\documentclass[12pt]{amsart}
\usepackage[margin=1in]{geometry}
\usepackage{amssymb,amsfonts,amsmath}
\usepackage{color}
\usepackage{soul}
\usepackage{enumerate}
\usepackage{mathrsfs}
\usepackage[bbgreekl]{mathbbol}
\usepackage{hyperref}
\usepackage[capitalise]{cleveref}
\usepackage{constants}

\newtheorem{thm}{Theorem}[section]
\newtheorem{cor}[thm]{Corollary}
\newtheorem{hyp}[thm]{Hypothesis}
\newtheorem{prop}[thm]{Proposition}
\newtheorem{lem}[thm]{Lemma}

\newtheorem{quest*}{Question}
\newtheorem{prob*}{Problem}

\theoremstyle{definition}

\theoremstyle{remark}
\newtheorem*{remark}{Remark}

\theoremstyle{remark}
\newtheorem*{remark1}{Remarks}

\numberwithin{equation}{section}

\crefname{figure}{Figure}{Figures}
\theoremstyle{plain}
\newtheorem*{thm*}{Theorem}
\crefname{thm}{Theorem}{Theorems}
\crefname{cor}{Corollary}{Corollarys}
\newtheorem*{cor*}{Corollary}
\crefname{cor*}{Corollary}{Corollarys}
\crefname{lem}{Lemma}{Lemmas}
\crefname{prop}{Proposition}{Propositions}
\crefname{conj}{Conjecture}{Conjectures}
\newtheorem*{conj*}{Conjecture}
\crefname{conj*}{Conjecture}{Conjectures}
\crefname{defn}{Definition}{Definitions}
\crefname{hyp}{Hypothesis}{Hypotheses}

\newcommand{\Z}{\mathbb{Z}}
\newcommand{\R}{\mathbb{R}}
\newcommand{\Q}{\mathbb{Q}}

\newcommand{\re}{\textup{Re}}
\newcommand{\im}{\textup{Im}}

\newcommand{\GL}{\mathrm{GL}}

\newcommand{\Ar}{\mathrm{Ar}}
\newcommand{\Sw}{\mathrm{Sw}}

\newcommand{\N}{\mathbb{N}}
\newcommand{\A}{\mathbb{A}}
\newcommand{\Run}{\underset{s=1}{\text{Res}}\;}

\renewcommand{\bar}{\overline}

\renewcommand{\epsilon}{\varepsilon}
\newcommand{\cE}{\mathcal{E}}

\newcommand{\cF}{\mathcal{F}}

\newcommand{\sL}{\mathscr{L}}

\renewcommand{\Im}{\mathrm{Im}}

\renewcommand{\pmod}[1]{\, (\mathrm{mod} {\, #1})}

\newcommand{\kp}{\mathfrak{p}}

\newcommand{\kq}{\mathfrak{q}}

\renewcommand{\Re}{\mathrm{Re}}

\renewcommand{\pmod}[1]{\left(\mathrm{mod}\,\,#1\right)}

\newconstantfamily{abcon}{symbol=c}

\makeatletter
\let\@wraptoccontribs\wraptoccontribs
\makeatother

\title{Zeros of Rankin-Selberg $L$-functions at the edge of the critical strip}
\author{Farrell Brumley}
\address{LAGA - Institut Galil\'ee, 99 avenue Jean Baptiste Cl\'ement, 93430 Villetaneuse, France}
\email{brumley@math.univ-paris13.fr}
\author{Jesse Thorner}
\address{Department of Mathematics, Stanford University, Stanford, CA 94305}
\email{jthorner@stanford.edu}
\author{Asif Zaman}
\address{Department of Mathematics, Stanford University, Stanford, CA 94305}
\email{aazaman@stanford.edu}

\address{King's College London, Department of Mathematics, Strand, London WC2R 2LS, England}
\email{colin.bushnell@kcl.ac.uk}

\address{Laboratoire de Math\'ematiques d'Orsay, Univ Paris-Sud, CNRS, Universit\'e Paris-Saclay, 91405 Orsay, France}
\email{Guy.Henniart@math.u-psud.fr}

\contrib[with an appendix by]{Colin J. Bushnell and Guy Henniart}

\begin{document}
\begin{abstract}
Let $\pi$ and $\pi_0$ be unitary cuspidal automorphic representations.  We prove log-free zero density estimates for Rankin-Selberg $L$-functions of the form $L(s,\pi\times\pi_0)$, where $\pi$ varies in a given family and $\pi_0$ is fixed.  These estimates are unconditional in many cases of interest; they hold in full generality assuming an average form of the generalized Ramanujan conjecture.  We consider applications of these estimates related to mass equidistribution for Hecke-Maass forms, the rarity of Landau-Siegel zeros of Rankin-Selberg $L$-functions, the Chebotarev density theorem, and $\ell$-torsion in class groups of number fields.
\end{abstract}

\maketitle

\setcounter{tocdepth}{1}
\tableofcontents

\section{Statement of the main results}


The generalized Riemann hypothesis (GRH) for Dirichlet $L$-functions implies that if $a$ and $q\geq 1$ are coprime integers, then there exists a prime\footnote{We write $f=O(g)$ or $f\ll g$ to mean that $|f|\leq c|g|$ for some absolute and effective constant $c>0$.  For a parameter $\nu$, we write $f=O_{\nu}(g)$ to mean that $c$ might depend on $\nu$ in an effective manner.  We write $f\asymp g$ to mean that $f=O(g)$ and $g=O(f)$, and similarly for $f\asymp_{\nu}g$ and $f\sim_{\nu}g$.} $p\ll(q\log q)^2$ such that $p\equiv a\pmod{q}$.  Linnik \cite{Linnik} unconditionally proved that the least such prime is $O(q^A)$, where $A>0$ is an absolute and effective constant; up to the quality of $A$, Linnik's result is commensurate with what GRH predicts.  Linnik's proof developed powerful results for the distribution of zeros of Dirichlet $L$-functions near the point $s=1$, including a log-free zero density estimate.  In this paper, we prove a flexible log-free zero density estimate for families of $L$-functions and consider the arithmetic consequences of such an estimate in several different settings.  We use this estimate to study mass equidistribution for Hecke-Maass forms, the rarity of Landau-Siegel zeros for Rankin-Selberg $L$-functions, the Chebotarev density theorem, and $\ell$-torsion in class groups of number fields.

In the spirit of Linnik's original result, Kowalski and Michel \cite[Theorem 5]{KM} proved a log-free zero density estimate for general families of automorphic $L$-functions in the conductor aspect.  To describe their result, let $\mathbb{A}_{\mathbb{Q}}$ be the ring of adeles over $\mathbb{Q}$, let $d\geq 1$ be a fixed integer, and let $\mathcal{A}(d)$ be the set of cuspidal automorphic representations of $\mathrm{GL}_d(\mathbb{A}_{\Q})$ with unitary central character.  We make the implicit assumption that the central character of each $\pi\in\mathcal{A}(d)$ is trivial on the positive reals; this discretizes $\mathcal{A}(d)$.  For each $\pi\in\mathcal{A}(d)$, let
\[
L(s,\pi)=\sum_{n\geq 1}\frac{a_{\pi}(n)}{n^s}=\prod_{p}\prod_{j=1}^{d}(1-\alpha_{j,\pi}(p)p^{-s})^{-1}
\]
be the standard $L$-function associated to $\pi$, where $p$ runs through the primes.  Consider a finite set $S(q)$  of distinct cuspidal automorphic representations $\pi\in\mathcal{A}(d)$ such that:
\begin{enumerate}
	\item There exists some $\delta>0$ (depending at most on $d$) such that for each $\pi\in S(q)$, each $1\leq j\leq d$, and each prime $p$, we have the bound $|\alpha_{j,\pi}(p)|\leq p^{1/4-\delta}$. 
	\item There exists a constant $A>0$ such that for all $\pi\in S(q)$, the conductor of $\pi$ is $O(q^A)$.
	\item There exists a constant $M>0$ such that $\# S(q)\ll q^{M}$.
	\item Each $\pi\in S(q)$ has the same component $\pi_{\infty}$ at the infinite place of $\mathbb{Q}$.
\end{enumerate}
 Note that the generalized Ramanujan conjecture (GRC) predicts that $|\alpha_{j,\pi}(p)|\leq 1$ for all primes $p$.
Define
\[
N_{\pi}(\sigma,T):=\#\{\rho=\beta+i\gamma\colon \sigma\leq\beta,~|\gamma|\leq T,~L(\rho,\pi)=0\}.
\]
With these conventions and hypotheses, Kowalski and Michel prove that there exists a constant $c=c(A,\delta,M)>M$ and a constant $B>0$ (depending on $S(q)$ but not $q$) such that
\begin{equation}
\label{eqn:KM}
\sum_{\pi\in S(q)}N_{\pi}(\sigma,T)\ll T^B q^{c\frac{1-\sigma}{2\sigma-1}},\qquad \frac{3}{4}<\sigma\leq 1,\quad T\geq 2.
\end{equation}
If $\sigma\geq 1-\frac{M}{c}$ and $T$ is sufficiently small with respect to $q$, then \eqref{eqn:KM} tells us that at most a vanishingly small proportion of low-lying zeros of the $L$-functions $L(s,\pi)$ with $\pi\in S(q)$ lie near $s=1$.  In many problems, such a result can serve as a powerful substitute for GRH.  Until now, \eqref{eqn:KM} appears to be the most flexible and robust zero density estimate for studying zeros of automorphic $L$-functions near $s=1$.

For a pair of automorphic representations $\pi\in\mathcal{A}(d)$ and $\pi_0\in\mathcal{A}(d_0)$, consider the associated Rankin-Selberg $L$-function
\begin{align*}
L(s,\pi\times\pi_0)=\sum_{n\geq 1}\frac{a_{\pi\times\pi_0}(n)}{n^s}=\prod_{p} \prod_{j=1}^{d} \prod_{j_0=1}^{d_0} (1-\alpha_{j,j_0,\pi\times\pi_0}(p)p^{-s})^{-1},
\end{align*}
where
\begin{equation}
\label{eqn:RS_split}
\{\alpha_{j,j_0,\pi\times\pi_0}(p)\colon 1\leq j\leq d,~1\leq j_0\leq d_0\}=\{\alpha_{j,\pi}(p)\alpha_{j_0,\pi_0}(p)\colon 1\leq j\leq d,~1\leq j_0\leq d_0\}
\end{equation}
for all but finitely many primes $p$.  In this paper, we establish log-free zero density estimates for families of Rankin-Selberg $L$-functions $L(s,\pi\times\pi_0)$, where $\pi$ varies and $\pi_0$ is fixed.  In order to make this precise, we define
\[
\mathcal{A}=\bigcup_{d\geq 1}\mathcal{A}(d),
\]
and we let $\mathcal{F}$ be a subset of $\mathcal{A}$.  We define
\begin{equation}
\label{eqn:family}
\mathcal{F}_m(Q) = \{\pi\in\mathcal{F}\colon C(\pi)\leq Q,~\pi\in\mathcal{F}\cap \mathcal{A}(d)\implies d\leq m\},
\end{equation}
where $C(\pi)$ is the analytic conductor of $\pi$ (see \eqref{eqn:analytic_conductor_def} for the definition).  We require an average version of GRC.

\begin{hyp}
\label{hyp}
Let $\pi\in\mathcal{A}(d)$.  For all $\epsilon>0$,
\[
\prod_p \sum_{r=0}^{\infty}\frac{\max_{1\leq j\leq d}|\alpha_{j,\pi}(p)|^{2r}}{p^{r(1+\epsilon)}}\ll_{d,\epsilon} C(\pi)^{\epsilon}.
\]
\end{hyp}
\begin{remark}
Indeed, if $\pi$ satisfies GRC, then \cref{hyp} follows with lots to spare.  Brumley \cite[Theorem 1 and Corollary 2]{Brumley_2} proved that each $\pi\in\mathcal{A}(d)$ satisfies \cref{hyp} when $d\leq 4$ and gave sufficient conditions (strictly weaker than assuming GRC in full) under which $\pi$ may satisfy \cref{hyp} when $d\geq 5$.
\end{remark}


\begin{thm}
	\label{thm:lfzde_RS}
	Let $\pi_0\in\mathcal{A}(m_0)$, and let $Q,T\geq 1$.  Let $\mathcal{F}_m(Q)$ be as in \eqref{eqn:family}, and suppose that $\pi_0$ and each $\pi\in\mathcal{F}_m(Q)$ satisfy \cref{hyp}.  If $1/2\leq\sigma\leq 1$, then
	\[
	\sum_{\pi\in\mathcal{F}_m(Q)}N_{\pi\times\pi_0}(\sigma,T)\ll_{m,m_0} (C(\pi_0)QT)^{10^7 (m_0 m)^4 (1-\sigma)}.
	\]
\end{thm}

\begin{remark}
When $\pi_0\in\mathcal{A}(1)$ is the trivial representation, whose corresponding $L$-function is the Riemann zeta function, \cref{thm:lfzde_RS} immediately recovers \eqref{eqn:KM} (up to the quality of the coefficient of $1-\sigma$) with the added benefit of a significantly improved dependence on $T$.  \cref{thm:lfzde_RS} is new for all other choices of $\pi_0$, even if one assumes GRC in full.
\end{remark}
\begin{remark}
For simplicity, we have made no attempt to optimize the exponent, but there is room for some noticeable improvement (especially if one assumes GRC).  Obtaining such a numerical improvement was big component of the work in \cite{TZ1} (see Theorem 3.2).
\end{remark}

Our proof of \cref{thm:lfzde_RS} in fact produces the upper bound
\begin{equation}
\label{eqn:family_dependence}
\sum_{\pi\in\mathcal{F}_m(Q)}N_{\pi\times\pi_0}(\sigma,T)\ll_{m,m_0} (C(\pi_0)QT\#\mathcal{F}_m(Q))^{10^6 (m_0 m)^3(1-\sigma)}
\end{equation}
(see \eqref{eqn:bound_without_family}).  However, \eqref{eqn:family_dependence} only becomes meaningful when there exists a constant $c_m>0$ (depending only on $m$) such that $\#\mathcal{F}_m(Q)\ll_{\mathcal{F},m}Q^{c_m}$.  The situation is the same as in \eqref{eqn:KM}, which is why Kowalski and Michel assume the bound $\#S(q)\ll q^M$.  A standard calculation for Dirichlet characters reveals that $\#\mathcal{F}_1(Q)\ll Q^{2}$, and the existence of some suitable $c_m>0$ for $m\geq 2$ follows from work of Michel and Venkatesh \cite[Section 2.6.5]{MV}. We expect that $\mathcal{F}_m(Q)\sim_{\mathcal{F},m}Q^{m+1}$ for all $m\geq 1$; Brumley and Mili{\'c}evi{\'c} \cite[Theorems 1.1 and 1.2]{BM} proved this claim (and much more) when $m=2$.  For $m\geq 3$, Brumley and Mili{\'c}evi{\'c} prove the claim when each $\pi\in\mathcal{F}_m(Q)$ corresponds to a Hecke-Maass newform.  We unconditionally prove:
\begin{thm}
	\label{thm:universal}
	For all $\epsilon>0$, we have the bound $\#\mathcal{F}_m(Q)\ll_{\epsilon,m}Q^{2m+\epsilon}$.
\end{thm}
The truth of \cref{thm:universal} follows immediately from \cref{thm:poly-growth}, which we prove in the appendix.  The bound in \cref{thm:universal} along with \eqref{eqn:family_dependence} produces \cref{thm:lfzde_RS}.

The bound \cref{thm:lfzde_RS} improves noticeably if $\pi_0$ satisfies GRC and there exists a primitive real Dirichlet character $\chi\pmod{q}$ with $q\leq 2Q$ such that $L(s,\chi)$ has real zero close to $s=1$.


\begin{thm}
\label{thm:Landau_Siegel}
Let $\pi_0\in\mathcal{A}(m_0)$ satisfy GRC, and let $Q,T\geq 1$.  Let $\mathcal{F}_m(Q)$ be as in \eqref{eqn:family}.  Suppose that each $\pi\in\mathcal{F}_m(Q)$ satisfies \cref{hyp}.  Let $\chi\pmod{q}$ be a real primitive Dirichlet character with $q\leq 2Q$.  If $1/2\leq\sigma\leq 1$, then
	\[
	\sum_{\pi\in\mathcal{F}_m(Q)}N_{\pi\times\pi_0}(\sigma,T)\ll_{m,m_0} \min\{1,(1-\beta_{\chi})\log(QT)\}\cdot(C(\pi_0)QT)^{10^7 (m_0 m)^4 (1-\sigma)},
	\]
	where $\beta_{\chi}$ denotes the largest real zero of the Dirichlet $L$-function $L(s,\chi)$, if it exists.
\end{thm}

Page's theorem \cite[Chapter 14]{Davenport} tells us that there exists an absolute and effective constant $\Cl[abcon]{c_zielgel}>0$ such that for every $Q\geq 3$, there exists at most one modulus $q\in(Q,2Q]$ and at most one primitive real character $\chi\pmod{q}$ such that $L(s,\chi)$ has a real zero $\beta_{\chi}$ with the property that $\beta_{\chi}\geq 1-\Cr{c_zielgel}/\log q$.  Moreover, such a zero $\beta_{\chi}$, which we call a Landau-Siegel zero, must be simple.  If a primitive real character $\chi\pmod{q}$ with $q\in(Q,2Q]$ has an associated Landau-Siegel zero $\beta_{\chi}$, then \cref{thm:Landau_Siegel} improves on \cref{thm:lfzde_RS}.  While it is well-known that Landau-Siegel zeros associated to real characters repel the zeros of Dirichlet $L$-functions from the point $s=1$, \cref{thm:Landau_Siegel} appears to be the first explicit instance in the literature where Landau-Siegel zeros associated to real characters repel zeros of high-degree $L$-functions.  This adds to the growing literature on interesting consequences of the existence of Landau-Siegel zeros of Dirichlet $L$-functions \cite{CI16,DM17,FI03,FI04,FI05,FI13,HB83}.

Our proof of \cref{thm:lfzde_RS}, which is noticeably different from that of \eqref{eqn:KM}, descends naturally from Gallagher's approach to log-free zero density estimates for Dirichlet $L$-functions \cite{Gallagher}.  Much like the classical approach to zero-free regions for $L$-functions, if $L(s,\pi\times\pi_0)$ has a zero $\rho_0$ such that $|\rho_0-(1+it)|\leq\epsilon$ for some small $\epsilon>0$, then high derivatives of $-L'/L(s,\pi\times\pi_0)$ near $s=1+\epsilon+it$ will be large; this is made quantitative via the lower bound for power sums due to S{\'o}s and Tur{\'a}n \cite{Turan}.  Moreover, one can show that if these derivatives are large, then the mean value of a certain Dirichlet polynomial roughly of the shape
\[
P(t,\pi\times\pi_0)=\sum_{A<p<B}\frac{a_{\pi\times\pi_0}(p)\log p}{p^{1+it}}
\]
must also be large when $t$ is close to $\im(\rho_0)$.    A new ``pre-sifted'' large sieve inequality (\cref{thm:large_sieve}) in the spirit of Duke and Kowalski \cite[Theorem 4]{DK}  shows that the mean value of $P(t,\pi\times\pi_0)$ cannot be large for too many $\pi\in\mathcal{F}_m(Q)$ simultaneously; \cref{thm:lfzde_RS} follows once this is made precise. The coefficients of $P(t,\pi\times\pi_0)$ are supported on large unramified primes, in which case $a_{\pi\times\pi_0}(p)=a_{\pi}(p)a_{\pi_0}(p)$ by means of \eqref{eqn:RS_split}; this decisive identity facilitates the averaging over $\pi\in\mathcal{F}_m(Q)$ while keeping $\pi_0$ fixed.  We prove \cref{thm:Landau_Siegel} similarly by simultaneously considering the twists $L(s,\pi\times\pi_0)$ and $L(s,\pi\times(\pi_0\otimes\chi))$ and exploiting the fact that if $\chi$ is a real primitive Dirichlet character with a Landau-Siegel zero, then $\chi$ behaves like the M{\"o}bius function.  This approach contrasts with the method of proof for \eqref{eqn:KM}, which uses mollification to detect zeros and a mean value theorem involving Selberg's pseudo-characters to show that the aggregate contributions from the zeros of each $L$-function is small.  It is unclear to the authors how one would modify the proof of \eqref{eqn:KM} to incorporate a twist by $\pi_0$ while maintaining a log-free estimate.

In \cite[Corollary 2.6]{ST}, Soundararajan and the first author establish the first unconditional log-free zero density estimate for each individual Rankin-Selberg $L$-function $L(s,\pi\times\pi_0)$ with an application to the weak subconvexity problem.  The proof of \cite[Corollary 2.6]{ST} relies on the same method of detecting zeros that we use here.  Unfortunately, the means by which the proofs in \cite{ST} avoid appealing to a weak form of GRC (such as \cref{hyp}) appears to be incompatible with the process of averaging over $\pi\in\mathcal{F}_m(Q)$.  In particular, \cref{hyp} appears to be indispensable in the proof of \cref{thm:large_sieve} unless $\#\mathcal{F}_m(Q)=1$, which is precisely the case considered in \cite{ST}.


\section{Arithmetic applications}
\label{sec:Landau_Siegel}


\subsection{Subconvexity and mass equdistribution}
\label{sec:Subconvexity}

Let $f$ be a Hecke-Maass newform for the congruence subgroup $\Gamma_0(q_f)\subset\mathrm{SL}_2(\Z)$ with Laplace eigenvalue $\lambda_f$ and trivial central character.  Define
\begin{equation}
	\label{eqn:maass_family}
\mathscr{G}(Q)=\{\textup{$f$: $q_f$ squarefree, $\lambda_f q_f\leq Q$}\}.
\end{equation}
Let $f_0$ denote a fixed Hecke-Maass newform, and consider the $L$-functions $L(s,f\times f)$ and $L(s,f\times f \times f_0)$ as $f\in\mathscr{G}(Q)$ varies.  Since $q_f$ is squarefree, the conductor of $f\times f$ is $q_f^2$.

The generalized Lindel{\"o}f hypothesis (which follows from GRH) predicts that for all $\epsilon>0$ and all $f\in\mathscr{G}(Q)$, we have the bounds I explicate the $t$ and $f_0$ dependence:
\[
L(1/2+it,f\times f)\ll_{\epsilon} ((|t|+1)^4\lambda_f q_f^2)^{\epsilon},\qquad L(1/2,f\times f \times f_0)\ll_{\epsilon}(\lambda_{f_0}^4 q_{f_0}^4\lambda_f^2 q_f^4)^{\epsilon}.
\]
The so-called convexity bounds
\[
L(1/2+it,f\times f)\ll ((|t|+1)^4\lambda_f q_f^{2})^{1/4},\qquad L(1/2,f\times f \times f_0)\ll (\lambda_{f_0}^4 q_{f_0}^4\lambda_f^{2}q_f^4)^{1/4}
\]
follow from the work of Heath-Brown \cite{HB}.  Subconvexity bounds of the shape 
\begin{equation}
\label{eqn:subconvexity}
L(1/2+it,f\times f)\ll ((|t|+1)^4\lambda_f q_f^{2})^{1/4-\delta},\qquad L(1/2,f\times f \times f_0)\ll (\lambda_{f_0}^4 q_{f_0}^4\lambda_f^{2}q_f^4)^{1/4-\delta}
\end{equation}
are not yet known; obtaining bounds of these sorts is a very active area of research which has some spectacular partial results (see \cite[Theorem 1.1]{Sound_weak}, for instance).

A standard calculation involving the approximate functional equation for Dirichlet $L$-functions and the large sieve shows that if $Q$ is large, then for all except at most a density zero subset of the moduli $q\leq Q$, we have the bound $L(1/2,\chi)\ll_{\epsilon}q^{\epsilon}$ for all primitive Dirichlet characters $\chi\pmod q$.  Similarly, a sufficiently strong analogue of the large sieve for automorphic forms will show that there exists a constant  $\delta>0$ such that \eqref{eqn:subconvexity} holds for almost all $f\in\mathscr{G}(Q)$.  The best candidate for such a large sieve is that of Duke and Kowalski \cite[Theorem 4]{DK}, but it falls short because the best unconditional bound toward GRC for Hecke-Maass newforms is not strong enough (though assuming GRC in full is not necessary).  However, a straightforward application of \cref{thm:lfzde_RS} yields such an average result.

\begin{thm}
	\label{thm:subconvexity}
	Let $\epsilon>0$, and let $\mathscr{G}(Q)$ be as in \eqref{eqn:maass_family}. For all except at most $O_{f_0}(Q^{\epsilon})$ of the Hecke-Maass forms $f\in \mathscr{G}(Q)$, the bounds in \eqref{eqn:subconvexity} hold simultaneously with $\delta=10^{-20}\epsilon$.
\end{thm}

\begin{remark}
It follows from recent work of Brumley and Mili{\'c}evi{\'c} \cite{BM} that
\begin{equation}
\label{eqn:size_maass_family}
\#\mathscr{G}(Q)\asymp Q^2.
\end{equation}
Thus \cref{thm:subconvexity} is nontrivial when $\epsilon$ is sufficiently small.  (In the discussion in \cite[Section 3]{BM}, one can replace $\Gamma_1(q)$ with $\Gamma_0(q)$ without loss, which yields \eqref{eqn:size_maass_family}.)
\end{remark}

Our interest in \eqref{eqn:subconvexity} is motivated by the quantum unique ergodicity conjecture.  Lindenstrauss \cite{Lindenstrauss} and Soundararajan \cite{Sound_QUE_Maass} proved that as $f$ traverses the Hecke-Mass forms with $q_f=1$ and $\lambda_f\to\infty$, the $L^2$ mass of $f$ equidistributes in $\Gamma_0(1)\setminus\mathbb{H}$ with respect to the standard hyperbolic measure.  This affirmatively resolved the quantum unique ergodicity conjecture of Rudnick and Sarnak \cite{RS} for the modular surface.  More specifically, let 
\begin{equation}
\label{eqn:def_mu_f}
\mu_f(\phi)=\int_{\Gamma_0(q_f)\setminus\mathbb{H}}|f(z)|^2\phi(z)\frac{dx dy}{y^2},\qquad \mu(\phi)=\int_{\Gamma_0(1)\setminus\mathbb{H}}\phi(z)\frac{dx dy}{y^2},
\end{equation}
where $\phi$ is a bounded measurable function on $\Gamma_0(1)\setminus\mathbb{H}$.  It is now known that as $f$ traverses the Hecke-Maass forms of eigenvalue $\lambda_f\to\infty$ with $q_f=1$,
\begin{equation}
\label{eqn:discrepancy}
D_f(\phi):=\frac{\mu_f(\phi)}{\mu_f(1)}-\frac{\mu(\phi)}{\mu(1)}=o_{\phi}(1).
\end{equation}
Unfortunately, the methods in \cite{Lindenstrauss,Sound_QUE_Maass} do not yield any information about the rate of convergence in \eqref{eqn:discrepancy}.  See \cite{HS,Nelson,NPS} for an unconditional proof of \eqref{eqn:discrepancy} with an effective rate of convergence as $f$ traverses the holomorphic cuspdial newforms of weight $k_f$ and level $q_f$ with $k_f q_f\to\infty$; this proof relies heavily on the fact that GRC is known for such newforms.  For work in the direction of establishing \eqref{eqn:discrepancy} for Hecke-Maass forms in $q_f$-aspect when $q_f$ is large and prime, see \cite{Nelson2}.

We consider the problem of proving that for all except at most a density zero subset of $f\in \mathscr{G}(Q)$, one has \eqref{eqn:discrepancy} with a power-saving rate of convergence in the hybrid $q_f$ and $\lambda_f$ aspects.  When $f$ traverses the even Hecke-Maass forms with $q_f=1$, this follows from Zhao's computation of the quantum variance of the modular surface \cite{Zhao}.  It is unclear to the authors whether one can adapt the proofs for the problem considered here.

Nelson \cite{Nelson} proved that for $f\in\mathscr{G}(Q)$ (given by \eqref{eqn:maass_family}), subconvexity bounds of the form \eqref{eqn:subconvexity} imply the bound
\begin{equation}
\label{eqn:power_saving_rate}
D_f(\phi)\ll_{\phi}(\lambda_f q_f^2)^{-\delta+o(1)}.
\end{equation}
(See Remarks 1.4 and 1.7 as well as Section 4 of \cite{Nelson}.)  Thus the next result follows immediately from \cref{thm:subconvexity} and the remark that follows it.

\begin{cor}
	\label{thm:QUE_1}
	Fix $\epsilon>0$, and let $\mathscr{G}(Q)$ be as in \eqref{eqn:maass_family}.  For all except at most $O_{\phi}(Q^{\epsilon})$ of the Hecke-Maass forms $f\in\mathscr{G}(Q)$, the bound \eqref{eqn:power_saving_rate} holds with $\delta=10^{-20}\epsilon$.
\end{cor}

\begin{remark}
By appealing to the extension of Watson's formula proved by Nelson, Pitale, and Saha (see \cite{NPS}) and the calculations in \cite[Page 11]{DK}, one can extend the definition of $\mathscr{G}(Q)$ to allow $q_f$ to be {\it any} integer at the cost of allowing the exceptional set to be of size $O_{\phi}(Q^{1/2+\epsilon})$ in \cref{thm:QUE_1}. The proof is entirely analogous.
\end{remark}

\subsection{Rarity of Landau-Siegel zeros} 
Let $\mathcal{F}_m(Q)$ be as in \eqref{eqn:family}, and let $\pi\in\mathcal{F}_m(Q)\cap\mathcal{A}(d)$.  While GRH predicts that $L(s,\pi)$ has no zero in the region $\re(s)>1/2$, at present we know that $L(s,\pi)$ has at most one zero in the region
\begin{equation}
\label{eqn:ZFR_cusp_form}
\re(s)\geq 1-\frac{\Cl[abcon]{ZFR_3}}{d^4\log(C(\pi)(|\im(s)|+3))}.
\end{equation}
(See \cite[Theorem 5.10]{IK}.) If $L(s,\pi)$ has a zero in this region, then $\pi$ is self-dual (so the Dirichlet coefficients of $L(s,\pi)$ are real), and the zero must be simple and real.  We call such a zero a Landau-Siegel zero.  Hoffstein and Ramakrishnan \cite[Theorem A]{Hoffstein} proved that such Landau-Siegel zeros are quite rare.  In particular, for some suitable effective constant $c(m)>0$, there is at most one $\pi\in\mathcal{F}_m(Q)$ such that $L(s,\pi)$ has a real zero $\beta$ satisfying $\beta>1-c(m)/\log Q$.  This generalizes Page's theorem for Dirichlet characters.  Moreover, it is known by the work of Hoffstein and Ramakrishnan \cite[Theorem C]{Hoffstein} and Banks \cite{Banks} that if $m=2$ or $3$, then no $\pi\in\mathcal{F}_m(Q)$ has an $L$-function possessing a Landau-Siegel zero.  The proof of \cite[Theorem A]{Hoffstein} relies crucially on the cuspidality of the $\pi\in\mathcal{F}_m(Q)$.

The situation for Rankin-Selberg $L$-functions is much more difficult.  Currently, an unconditional zero-free region (with at most one exceptional zero) roughly of the shape \eqref{eqn:ZFR_cusp_form} exists for $L(s,\pi\times\pi_0)$ when at least one of $\pi$ and $\pi_0$ is self-dual (see \cite{Humphries} for further discussion).  Since it is not known in general whether $L(s,\pi\times\pi_0)$ factors into a product of $L$-functions associated to cuspidal automorphic representations (though this is expected), it is unclear how to unconditionally generalize \cite[Theorem C]{Hoffstein} to establish the rarity of Landau-Siegel zeros for Rankin-Selberg $L$-functions.  Despite these setbacks, one can still show that few Rankin-Selberg $L$-functions have a Landau-Siegel zero.

\begin{thm}
\label{thm:page}
Assume the above notation.  Let $A>0$, and let $\mathcal{S}=\mathcal{S}(A,Q,T,\mathcal{F}_m(Q))$ be the set of all $\pi\in\mathcal{F}_m(Q)$ such that $L(s,\pi\times\pi_0)$ has a zero in the region
\[
s=\sigma+it,\qquad |t|\leq T,\qquad \sigma \geq 1-\frac{A}{10^7 (m_0 m)^4 \log(C(\pi_0)Q(T+2))}.
\]
\begin{enumerate}[(i)]
	\item Under the hypotheses of \cref{thm:lfzde_RS}, $|\mathcal{S}|=O_{m,m_0}(e^{A})$.
	\item Let $\chi\pmod{q}$ be a primitive real Dirichlet character modulo $q\leq 2Q$.  Under the hypotheses of \cref{thm:Landau_Siegel}, $|\mathcal{S}|=O_{m,m_0}(e^{A}\cdot\min\{1,(1-\beta_{\chi})\log(QT)\})$.
\end{enumerate}
\end{thm}

If there exists a sequence of primitive real characters $\chi\pmod{q}$ with $q\in(Q,2Q]$ such that $(1-\beta_{\chi})\log Q\to 0$ as $Q\to\infty$, then the size of the exceptional set in \cref{thm:page}(ii) is zero once $Q$ is sufficiently large relative to $T$.  Therefore, under \cref{hyp} for all cusp forms, the existence of a sequence of primitive real characters whose $L$-functions have a Landau-Siegel zero implies the nonexistence of Landau-Siegel zeros for all other Rankin-Selberg $L$-functions of comparable analytic conductor.  This provides an interesting companion to another result of Hoffstein and Ramakrishnan \cite[Theorem B]{Hoffstein} which roughly states that if all Rankin-Selberg $L$-functions factor into products of $L$-functions of cuspidal automorphic representations (as predicted by Langlands), then the only primitive $L$-functions over $\Q$ which could possibly admit a Landau-Siegel zero are those associated to primitive real Dirichlet characters.

Suppose that $\pi_0\not=\widetilde{\pi}$ for all $\pi\in\mathcal{F}_m(Q)$.  By setting $T=Q$ and $A=\log(C(\pi_0)Q)$, it follows readily from \cref{thm:page}(i) that apart from at most a few exceptional $\pi$ in $\mathcal{F}_m(Q)$, one can obtain strong approximations for $L(1,\pi\times\pi_0)$ as a short Euler product.  See \cite{CM,GS,Lamzouri} for further discussion and applications of such approximations.

\subsection{The Chebotarev density theorem in families} \label{subsec:CDT_intro}

Let $K$ be a number field of degree  $n=[K:\Q]$ with $D_K = |\mathrm{disc}(K/\Q)|$ and Galois closure $\widetilde{K}$ over $\Q$.  Let $G$ be isomorphic to the Galois group of $\widetilde{K}/\Q$, and let $C$ be a conjugacy class of $\mathrm{Gal}(\widetilde{K}/\Q)$.  Consider the prime counting function
\[
\pi_C(x,\widetilde{K}/\Q):=\#\Big\{p\leq x\colon p\nmid D_{\widetilde{K}},~\Big[\frac{\widetilde{K}/\Q}{p}\Big]=C\Big\},
\]
where the Artin symbol $[\frac{\widetilde{K}/\Q}{p}]$ denotes the conjugacy class of Frobenius automorphisms attached to the prime ideals of $\widetilde{K}$ which lie over $p$.  The Chebotarev density theorem states that as $x\to\infty$,
\[
\mathcal{E}_C(x,\widetilde{K}/\Q):=\Big|\pi_C(x,\widetilde{K}/\Q)-\frac{|C|}{|G|}\pi(x)\Big|=o\Big(\frac{|C|}{|G|}\frac{x}{\log x}\Big),
\]
where $\pi(x)$ is the number of rational primes up to $x$. It follows from the work of Lagarias and Odlyzko \cite[Theorem 1.1]{LO} that GRH for the Dedekind zeta function $\zeta_{\widetilde{K}}(s)$ implies
\begin{equation}
\label{eqn:LO}
\mathcal{E}_C(x,\widetilde{K}/\Q)\ll \frac{|C|}{|G|}x^{1/2}\log(D_{\widetilde{K}} x^{|G|})\quad \text{ for } x \geq (\log D_{\widetilde{K}})^2(\log\log D_{\widetilde{K}})^4.
\end{equation}
We need the $\log\log$ if we are to have an asymptotic.  The least prime in Chebotarev being of size $(\log D)^2$ is a consequence of smoothing.  Unconditionally, refining a result of Lagarias and Odlyzko \cite{LO}, it follows from work of Murty \cite[Section 4]{VKM} that 
 \begin{equation}
\label{eqn:CDT_1}
\mathcal{E}_C(x,\widetilde{K}/\Q)\ll \frac{|C|}{|G|}\Big(\frac{x^{\beta_{1}}}{\log x}+\frac{x}{\exp(\Cl[abcon]{LOzfr}(\log x)^{1/2}|G|^{-1/2})}\Big)\qquad \text{for } x\gg_{|G|} e^{\Cl[abcon]{LOrange} (\log D_{\widetilde{K}} )^2/|G| },
\end{equation}
where $\beta_{1}$ is a putative Landau-Siegel zero of $\zeta_{\widetilde{K}}(s)$. Recent work of the authors \cite{TZ3} shows for any $A > 1$, there exists $B = B(A) > 1$ such that
\begin{equation}
\label{eqn:CDT_1-1}
\mathcal{E}_C(x,\widetilde{K}/\Q)\ll_A \frac{|C|}{|G|}\Big(\frac{x^{\beta_{1}}}{\log x}+\frac{x}{(\log x)^A} \Big)\qquad \text{for } x\gg_{|G|,A} D_{\tilde{K}}^{B \log \log D_{\widetilde{K}}}.
\end{equation}
For large $x$, \eqref{eqn:CDT_1} remains the strongest upper bound for $\mathcal{E}_C$ and it is non-trivial in the absence of a Landau--Siegel zero. On the other hand, \eqref{eqn:CDT_1-1} exhibits a weaker estimate but for much smaller values of $x$. Nonetheless, even when ignoring the Landau--Siegel zero, both \eqref{eqn:CDT_1} and \eqref{eqn:CDT_1-1} fall far short of exhibiting non-trivial bounds for values of $x$ commensurate in size with \eqref{eqn:LO}. Even establishing such bounds for $x \geq D_{\widetilde{K}}^{o(1)}$ would be extremely desirable. 


Substantial progress has recently been made by Pierce, Turnage-Butterbaugh, and Wood \cite{PTW} when $K$ varies in certain families. They show that the ranges of $x$ in \eqref{eqn:CDT_1} and \eqref{eqn:CDT_1-1} can be significantly improved for most $K$.  We briefly summarize their results.  Let $G\in\{C_m,D_p,S_3,S_4,A_4\}$, where $C_m$ is a cyclic group of order $m\geq2$, $S_m$ is a symmetric group acting on $m\geq2$ elements, $D_p$ is a dihedral group of order $2p$ with $p$ an odd prime, and $A_4$ is an alternating group acting on 4 elements.  Let $\mathscr{F}(X) = \mathscr{F}(X;G,n,R_G)$ denote the set of number fields $K$ with $[K:\Q]=n$ and $D_K\leq X$ such that $\mathrm{Gal}(\widetilde{K}/\Q)\cong G$ and each $K$ satisfies a certain arithmetic restriction $R_G$ depending only on $G$.  In particular,
\begin{equation*}
	R_G=\begin{cases}
		\text{$K$ is totally ramified}&\mbox{if $G=C_m$,}\\[2mm]
		\text{$K$ has square-free absolute discriminant}&\mbox{if $G=S_m$,}\\[2mm]
		\text{Every prime $p$ that ramifies tamely in $K$ has its inertia group}\\
        \text{generated by an element in the conjugacy class of reflections}&\mbox{if $G=D_p$,}\\[2mm]
        \text{Every prime $p$ that ramifies tamely in $K$ has inertia group}\\
         \text{generated by an element in either $\{(1~2~3),(1~3~4),(1~4~2),(2~4~3)\}$} \\
         \text{or $\{(1~3~2),(1~4~3),(1~2~4),(2~3~4)\}$}&\mbox{if $G=A_4$.}
	\end{cases}
\end{equation*}
As demonstrated in \cite{PTW}, there exists some constant $a=a(G,n)\in(0,1]$ such that, for all choices of $G$, $n$, and $R_G$ under consideration,  $\#\mathscr{F}(X)\gg_{G,n} X^a$.

With this setup in mind, let $A\geq 2$ and $\eta > 0$.  Pierce, Turnage-Butterbaugh, and Wood \cite[Theorem 1.4]{PTW} proved that there exist effective constants $\alpha=\alpha(\eta, A,G,n)>0$ and $\epsilon = \epsilon(G,n)>0$ such that for all fields $K \in \mathscr{F}(X)$ with at most $O_{G,n}(X^{-\epsilon} \# \mathscr{F}(X))$ exceptions,  one has
\begin{equation}
\label{eqn:CDT_2}
\mathcal{E}_C(x,\widetilde{K}/\Q)\ll_A \begin{cases}
	\cfrac{|C|}{|G|}\cfrac{x}{(\log x)^A}&\mbox{if $e^{\alpha(\log\log D_{\widetilde{K}})^{\frac{5}{3}+\eta}} \leq x \ll_{|G|} e^{\Cr{LOrange} (\log D_{\widetilde{K}} )^2/|G| }$},\\[5mm]
	\cfrac{|C|}{|G|}\cfrac{x}{\exp(\Cr{LOzfr}(\log x)^{1/2}|G|^{-1/2})}
	 &\mbox{if $x \gg_{|G|}  e^{\Cr{LOrange} (\log D_{\widetilde{K}} )^2/|G| }$.} 
\end{cases}
\end{equation}
Notice \eqref{eqn:CDT_2} eliminates the Landau--Siegel zero and, most importantly,  goes beyond the range of $x$ in \eqref{eqn:CDT_1-1}. Somewhat surprisingly, when a Landau--Siegel zero does not exist, the estimate for $\cE_C$ in \eqref{eqn:CDT_2} surpasses neither  \eqref{eqn:CDT_1} nor \eqref{eqn:CDT_1-1} in their respective weaker ranges of $x$. (We have only collected their unconditional results; see \cite[Section 2]{PTW} for a discussion regarding degree $n$ $S_n$- and $A_n$-fields with $n\geq 5$.) 

The proofs in \cite{PTW} rely decisively on \eqref{eqn:KM}, and the $T$-dependence in \eqref{eqn:KM} inhibits their proof from achieving a result that is more commensurate with what GRH predicts in \eqref{eqn:LO}. Using \cref{thm:lfzde_RS} and \cref{thm:QuasiGRH}, we improve both the range of $x$ and quality of error term in \eqref{eqn:CDT_2}. In particular, we obtain a range much closer to what GRH predicts with a power savings error term for small values of $x$. 
\begin{thm}
\label{thm:CDT}
Let $G$ be isomorphic to one of	$C_m$, $S_3$, $S_4$, $D_p$, or $A_4$; let $C\subset G$ be a conjugacy class; and let $\mathscr{F}(X) = \mathscr{F}(X;G,n,R_G)$ be as above.  There exist small positive constants $\eta = \eta(G,n)$ and $\epsilon = \epsilon(G,n)$  such that, for all fields $K\in \mathscr{F}(X)$ with at most $O_{G,n}( X^{-\epsilon} \mathscr{F}(X) )$ exceptions, 
\[
\Big|\pi_C(x,\widetilde{K}/\Q)-\frac{|C|}{|G|}\pi(x)\Big|\ll \begin{cases}
	\cfrac{|C|}{|G|}x^{1-\eta}&\mbox{if $(\log D_{\widetilde{K}})^{ 2/\eta}\leq x< D_{\widetilde{K}}^{1/(24 \eta)}$,}\\[5mm] 
	\cfrac{|C|}{|G|}\cfrac{x}{\exp(\Cr{LOzfr}(\log x)^{1/2}|G|^{-1/2})}&\mbox{if $x\geq D_{\widetilde{K}}^{1/(24 \eta)}$.}
\end{cases}
\]
\end{thm}
\begin{remark}
For a more uniform version of the error term in \cref{thm:CDT}, see \eqref{eqn:asif_hates_pauses}.
\end{remark}

\subsection{Landau--Siegel zeros and torsion in class groups} Let us continue with the notation of \cref{subsec:CDT_intro}. 
Let $\mathrm{Cl}_K$ denote the ideal class group of a number field $K$. It is widely believed that if $\ell$ is a positive integer, then the $\ell$-torsion subgroup $\mathrm{Cl}_K[\ell]$ is of size $O_{\epsilon,n,\ell}(D_K^{\epsilon})$ for all $\epsilon>0$, while the trivial bound is $O_{\epsilon,\ell,n}(D_K^{1/2+\epsilon})$.  Ellenberg and Venkatesh \cite[Lemma 2.3 and Proposition 3.1]{EV} proved that if, for any $\epsilon>0$, one has
\begin{equation}
\label{eqn:EV_hypothesis}
\#\Big\{p\leq D_K^{\frac{1}{2\ell(n-1)} - \epsilon}\colon \textup{$p\nmid D_K$ and splits completely in $K$}\Big\}\gg_{\epsilon,n} D_K^{\frac{1}{2\ell(n-1)}-\epsilon},
\end{equation}
then
\begin{equation}
\label{eqn:EV}
|\mathrm{Cl}_K[\ell]|\ll_{\epsilon,n,\ell}D_K^{\frac{1}{2}-\frac{1}
{2\ell(n-1)}+\epsilon}.
\end{equation}
Since primes that split  completely in $\widetilde{K}$ also split completely in $K$, the hypothesis \eqref{eqn:EV_hypothesis} follows easily from \eqref{eqn:LO}, which is a consequence of GRH.  It is a straightforward consequence of \eqref{eqn:CDT_2} that for any positive integer $\ell$, all except at most a density zero subset of the fields $K\in\mathscr{F}(X;G,n,R_G)$ satisfy \eqref{eqn:EV_hypothesis}, and hence \eqref{eqn:EV}, {\it unconditionally}.  This provides the first nontrivial upper bounds for $|\mathrm{Cl}_K[\ell]|$, for all integers $\ell\geq 1$, applicable to infinite families of fields of arbitrarily large degree. This elegant application of \eqref{eqn:CDT_2} in \cite{PTW} was achieved by exhibiting large zero-free regions for $\zeta_{\widetilde{K}}(s)$ for most fields $K$ in a given family. 

We proceed in a complementary direction using the zero repulsion phenomenon of a Landau--Siegel zero. If the Dedekind zeta function of a quadratic subfield $\Q(\sqrt{d})$ has a Landau-Siegel zero, then \cref{thm:Landau_Siegel} implies that certain number fields $K$, whose Galois closure does not contain $\Q(\sqrt{d})$ as a subfield, possess GRH-quality bounds on $\ell$-torsion in their class groups.

\begin{thm}
\label{thm:ell-torsion} 
Let $K/\Q$ be a number field of degree $n$ with Galois closure $\widetilde{K}$ over $\Q$. Let $\ell \geq 1$ be a positive integer and $\epsilon > 0$ be arbitrary. Let $\chi$ be the real Dirichlet character modulo a fundamental discriminant $d$. Assume the following:
\begin{enumerate}[(i)]
	\item $\zeta_{\widetilde{K}}(s)$ is the $L$-function of an automorphic representation of $\mathrm{GL}_{[\widetilde{K}:\Q]}(\mathbb{A}_{\Q})$.
	\item $\Q(\sqrt{d})\cap\widetilde{K}=\Q$ and   $\log D_K \asymp_{n,\epsilon,\ell} \log |d|$. 
	\item The Dirichlet $L$-function $L(s,\chi)$ has a real zero $\beta_{\chi} = 1 - \eta_{\chi}/\log d$ with $\eta_{\chi}$ sufficiently small, depending only on $n, \epsilon,$ and $\ell$. 
\end{enumerate}  
Then 
\[
|\mathrm{Cl}_K[\ell]|\ll_{\epsilon,n,\ell}D_K^{\frac{1}{2}-\frac{1}
{2\ell(n-1)}+\epsilon}.
\]
\end{thm}
\begin{remark1} 
~
\begin{enumerate}
\item We emphasize that \cref{thm:ell-torsion} is a pointwise bound, whereas the bounds in \cite{PTW} hold as one averages over $K$.
\item Arthur and Clozel \cite[page 223]{AC} proved that $\zeta_{\widetilde{K}}(s)$ is automorphic over $\Q$ when $\mathrm{Gal}(\widetilde{K}/\Q)$ is solvable.  Therefore, by a well-known result of Feit and Thompson, hypothesis (i) on $\zeta_{\widetilde{K}}(s)$ is satisfied when $|\mathrm{Gal}(\widetilde{K}/\Q)|$ is odd.
\end{enumerate}
\end{remark1}


\section{Properties of $L$-functions}
\label{sec:L-functions}


We recall some standard facts about $L$-functions arising from cuspidal automorphic representations and their Rankin-Selberg convolutions.  Much of the material we present here can be found in \cite[Section 1]{Brumley}.  We refer the reader there for a more detailed overview.

\subsection{Standard $L$-functions}

Let $d\geq1$ be an integer, let $\mathbb{A}$ denote the ring of adeles over $\Q$, and let $\mathcal{A}(d)$ be the set of all cuspidal automorphic representations of $\mathrm{GL}_d(\mathbb{A})$ (up to equivalence).  We consider each $\pi=\otimes_{p}\pi_{p}\in\mathcal{A}(d)$ to be normalized so that $\pi$ has unitary central character which is trivial on the positive reals; here, $p$ ranges over the primes.  We write $\widetilde{\pi}\in\mathcal{A}(d)$ for the representation which is contragredient to $\pi$.

Let $\pi=\otimes_p \pi_p\in\mathcal{A}(d)$, and let $N_{\pi}$ denote the conductor of $\pi$.  The standard $L$-function $L(s,\pi)$ associated to $\pi$ is of the form
\[
L(s,\pi)=\prod_{p} L(s,\pi_p)=\sum_{n=1}^{\infty}\frac{a_{\pi}(n)}{n^s}.
\]
The Euler product and Dirichlet series converge absolutely when $\re(s)>1$.  For each $p$, the local factor $L(s,\pi_p)$ is given in the form
\[
L(s,\pi_p)=\prod_{j=1}^{d}\Big(1-\frac{\alpha_{j,\pi}(p)}{p^{s}}\Big)^{-1}=1+\sum_{j=1}^{\infty}\frac{a_{\pi}(p^j)}{p^{js}}
\]
for suitable complex numbers $\alpha_{j,\pi}(p)$. With this convention, we have $\alpha_{j,\pi}(p)\neq0$ for all $j$ whenever $p\nmid N_{\pi}$, and it might be the case that $\alpha_{j,\pi}(p)=0$ for some $j$ when $p\mid N_{\pi}$.  At the archimedean place of $\Q$, there are $d$ complex Langlands parameters $\mu_{\pi}(j)$ from which we define
\[
L(s,\pi_{\infty}) = N_{\pi}\pi^{-\frac{ds}{2}}\prod_{j=1}^{d}\Gamma\Big(\frac{s+\mu_{\pi}(j)}{2}\Big).
\]
By the work of Rudnick and Sarnak \cite[Proposition A.1]{RS} and Blomer and Brumley \cite[Corollary 1]{BB}, we know that there exists a constant
\begin{equation}
\label{eqn:delta_d}
\delta_d\in\Big[0,\frac{1}{2}-\frac{1}{d^2+1}\Big]
\end{equation}
such that
\begin{equation}
\label{eqn:LRS_finite}
	| \alpha_{j,\pi}(p)|\leq  p^{\delta_d}\qquad\textup{and}\qquad\re(\mu_{\pi}(j))\geq -\delta_d
\end{equation}
for all $j$ and $p$.  The generalized Selberg eigenvalue conjecture and GRC assert that $\delta_d=0$ for all $d\geq 1$.  For each $p$,
\[
\{\alpha_{j,\widetilde{\pi}}(p)\}=\{\overline{ \alpha_{j,\pi}(p)}\},\qquad \{\mu_{\widetilde{\pi}}(j)\}=\{\overline{\mu_{\pi}(j)}\}.
\]

Let $r_{\pi}$ denote the order of the pole of $L(s,\pi)$ at $s=1$ and $\kappa_{\pi}$ be the residue of $L(s,\pi)$ at $s=1$.  The completed $L$-function
\[
\Lambda(s,\pi) = (s(s-1))^{r_{\pi}}N_{\pi}^{s/2}L(s,\pi)L(s,\pi_{\infty})
\]
is an entire function of order 1, and there exists a complex number $W(\pi)$ of modulus 1 such that for all $s\in\mathbb{C}$,
\[
\Lambda(s,\pi)=W(\pi)\Lambda(1-s,\widetilde{\pi}).
\]
On one hand, $L(s,\pi)$ has a zero at each pole of $L(s,\pi_{\infty})$; we call such a zero a trivial zero.  On the other hand, since $\Lambda(s,\pi)$ is entire of order 1, it has a Hadamard factorization
\[
\Lambda(s,\pi)=e^{a_{\pi}+b_{\pi}s}\prod_{\rho}\Big(1-\frac{s}{\rho}\Big)e^{s/\rho},
\]
where $\rho$ runs through the so-called nontrivial zeros of $L(s,\pi)$.

Finally, we define the analytic conductor of $\pi$ to be
\begin{equation}
\label{eqn:analytic_conductor_def}
C(\pi,t)=N_{\pi}\prod_{j=1}^d(1+|it+\mu_{\pi}(j)|),\qquad C(\pi)=C(\pi,0).
\end{equation}

\subsection{Rankin-Selberg $L$-functions}

Let $\pi=\otimes_p \pi_{p}\in\mathcal{A}(d)$ and $\pi'=\otimes_p \pi_{p}'\in\mathcal{A}(d')$.  The Rankin-Selberg $L$-function $L(s,\pi\times\pi')$ associated to $\pi$ and $\pi'$ is of the form
\[
L(s,\pi\times\pi')=\prod_{p}L(s,\pi_p\times\pi_{p}')=\sum_{n=1}^{\infty}\frac{a_{\pi\times\pi'}(n)}{n^s}.
\]
The Euler product and Dirichlet series converge absolutely when $\re(s)>1$.  For each $p$, the local factor $L(s,\pi_p)$ is given in the form
\[
L(s,\pi_p\times\pi_{p}')=\prod_{j=1}^{d}\prod_{j'=1}^{d'}(1-\alpha_{j,j',\pi\times\pi'}(p) p^{-s})^{-1}
\]
for suitable complex numbers $\alpha_{j,j',\pi\times\pi'}(p)$.  With $\delta_d$ as in \eqref{eqn:delta_d}, we have the pointwise bound
\begin{equation}
\label{eqn:LRS_2}
|\alpha_{j,j',\pi\times\pi'}(p)|\leq  p^{\delta_{d}+\delta_{d'}}\leq p^{1-\frac{1}{d'd}}.
\end{equation}
If $p\nmid N_{\pi}N_{\pi'}$, then we have the equality of sets
\begin{equation}
\label{eqn:separate_dirichlet_coeffs}
\{\alpha_{j,j',\pi\times\pi'}(p)\}=\{ \alpha_{j,\pi}(p)\alpha_{j',\pi'}(p)\}.
\end{equation}
At the archimedean place of $\Q$, there are $d'd$ complex Langlands parameters $\mu_{\pi\times\pi'}(j,j')$ from which we define
\[
L(s,\pi_{\infty}\times\pi_{\infty}')=N_{\pi\times\pi'}\pi^{-\frac{d'ds}{2}}\prod_{j=1}^{d}\prod_{j'=1}^{d'}\Gamma\Big(\frac{s+\mu_{\pi\times\pi'}(j,j')}{2}\Big).
\]
These parameters satsify
\[
\{\mu_{\widetilde{\pi}\times\widetilde{\pi}'}(j,j')\}=\{\overline{\mu_{\pi\times\pi'}(j,j')}\}
\]
and satisfy the pointwise bound
\begin{equation}
\label{eqn:LRS_22}
\Re(\mu_{\pi\times\pi'}(j,j'))\geq-\delta_{d}-\delta_{d'}\geq -1+(d'd)^{-1}.
\end{equation}

Let $r_{\pi\times\pi'}$ be the order of the pole of $L(s,\pi\times\pi')$ at $s=1$ and let $\kappa_{\pi\times \pi'}$ be the residue of $L(s,\pi\times\pi')$ at $s=1$.  By our normalization for $\pi$ and $\pi'$, we have that $r_{\pi\times\pi'}=1$ if and only if $\pi=\widetilde{\pi}'$; otherwise, $r_{\pi\times\pi'}=0$ and hence $\kappa_{\pi \times \pi'} = 0$.  The function
\[
\Lambda(s,\pi\times\pi')=(s(s-1))^{r_{\pi\times\pi'}}N_{\pi\times\pi'}^{s/2}L(s,\pi\times\pi')L(s,\pi_{\infty}\times\pi_{\infty}')
\]
is entire of order 1, and there exists a complex number $W(\pi\times\pi')$ of modulus 1 such that $\Lambda(s,\pi\times\pi')$ satisfies the functional equation
\[
\Lambda(s,\pi\times\pi')=W(\pi\times\pi')\Lambda(1-s,\widetilde{\pi}\times\widetilde{\pi}').
\]
On one hand, $L(s,\pi\times\pi')$ has a zero at each pole of $L(s,\pi_{\infty}\times\pi_{\infty}')$; we call such a zero a trivial zero.  On the other hand, since $\Lambda(s,\pi\times\pi')$ is entire of order 1, it has a Hadamard factorization
\[
\Lambda(s,\pi\times\pi')=e^{a_{\pi\times\pi'}+b_{\pi\times\pi'}s}\prod_{\rho}\Big(1-\frac{s}{\rho}\Big)e^{s/\rho},
\]
where $\rho$ runs through the so-called nontrivial zeros of $L(s,\pi\times\pi')$.

As with $L(s,\pi)$, we define the analytic conductor of $\pi\otimes\pi'$ to be
\[
C(\pi\times\pi',t)=N_{\pi\times\pi'}\prod_{j=1}^{d}\prod_{j'=1}^{d'}(1+|it+\mu_{\pi\times\pi'}(j,j')|),\qquad C(\pi\times\pi')=C(\pi\times\pi',0).
\]
It will be important to be able to decouple the dependencies of $C(\pi\times\pi',t)$ on $\pi$, $\pi'$, and $t$.  To this end, we have the combined work of Bushnell and Henniart \cite[Theorem 1]{BH} and Brumley \cite[Lemma A.2]{Humphries} which yields
\begin{equation}
\label{eqn:BH}
C(\pi\times\pi',t)\leq C(\pi)^{d'}C(\pi')^{d}(1+|t|)^{d' d},\qquad C(\pi\times\pi')\leq C(\pi)^{d'}C(\pi')^{d}.
\end{equation}

\section{Detecting zeros of $L$-functions}


Let $\Lambda(n)$ be the von Mangoldt function, and define the numbers
\begin{equation}
\label{eqn:lambda_def}
\lambda_{\pi\times\pi'}(n)=\begin{cases}
	\sum_{j=1}^{d}\sum_{j'=1}^{d'}\alpha_{j,j',\pi\times\pi_0}(p)^k&\mbox{if $n=p^k$ for a prime $p$},\\
	0&\mbox{otherwise}
\end{cases}
\end{equation}
so that
\[
-\frac{L'}{L}(s,\pi\times\pi')=\sum_{n=1}^{\infty}\frac{\lambda_{\pi\times\pi'}(n)\Lambda(n)}{n^s}.
\]
It follows from the definition of $\lambda_{\pi\times\pi'}(n)$ that if $\gcd(n,N_{\pi}N_{\pi'})=1$, then $\lambda_{\pi\times\pi'}(n)=\lambda_{\pi}(n)\lambda_{\pi'}(n)$.  In particular, if $\gcd(n,N_{\pi})=1$, then $|\lambda_{\pi}(n)|^2=\lambda_{\pi\times\widetilde{\pi}}(n)$.  During the proof of \cite[Lemma A.1]{RS}, it is shown that $\lambda_{\pi\times\widetilde{\pi}}(n)\geq 0$ for all $n\geq 0$.  Brumley \cite[Appendix]{ST} proved that regardless of whether $\gcd(n,N_{\pi}N_{\pi'})>1$, we always have the inequality
\begin{equation}
\label{eqn:Brumley_Cauchy_Schwarz}
	|\lambda_{\pi\times\pi'}(n)|\leq \sqrt{\lambda_{\pi\times\widetilde{\pi}}(n)\lambda_{\pi'\times\widetilde{\pi}'}(n)}\leq\frac{\lambda_{\pi\times\widetilde{\pi}}(n)+\lambda_{\pi'\times\widetilde{\pi}'}(n)}{2}.
\end{equation}
If $\chi$ is a primitive Dirichlet character modulo $q$, then $\chi\overline{\chi}$ is the trivial character modulo $q$, hence
\begin{equation}
\label{eqn:Brumley_Cauchy_Schwarz_2}
	|\lambda_{\pi\times(\pi'\otimes\chi)}(n)|\leq \sqrt{\lambda_{\pi\times\widetilde{\pi}}(n)\lambda_{(\pi'\otimes\chi)\times(\widetilde{\pi}'\otimes\overline{\chi})}(n)}\leq\frac{\lambda_{\pi\times\widetilde{\pi}}(n)+\lambda_{\pi'\times\widetilde{\pi}'}(n)}{2}.
\end{equation}

The proof of \cref{thm:lfzde_RS} will use the following result on the detection of zeros near the line $\re(s)=1$.

\begin{prop}
\label{prop:zero_detection}
Let $\pi\in\mathcal{F}_m(Q)$ and $\pi_0\in\mathcal{A}(m_0)$; suppose that both $\pi$ and $\pi_0$ satisfy \cref{hyp}.  Let $\chi\pmod{q}$ be a real primitive Dirichlet character, and let $\beta_{\chi}\geq 1/2$ denote a real zero of $L(s,\chi)$ (if it exists).  Let
\begin{equation}
\label{eqn:eta_def}
\frac{1}{\log(C(\pi_0)qQT)}\leq \eta \leq \frac{1}{10^7 (m_0 m)^2}
\end{equation}
and
\begin{equation}
\label{eqn:K_def1}
K\geq 4000(m_0 m)^2\eta\log(C(\pi_0)qQT)+O_{m_0,m}(1).
\end{equation}
with a sufficiently large implied constant.  If $L(s,\pi\times\pi_0)$ has a zero $\rho_0$ (trivial or nontrivial) satisfying $|\rho_0-(1+i\tau)|\leq \eta$ and $\rho_0\neq\beta_{\chi}$, then
\begin{align*}
1\ll (200)^{4K}\Big[\eta^3 \int_{A_1}^{A_2}\Big|\sum_{A_1<p\leq u}\frac{\lambda_{\pi\times\pi_0}(p)\log p}{p^{1+i\tau}}(1+\chi(p)p^{\beta_{\chi}-1})\Big|^2\frac{du}{u}+r_{\pi\times\pi_0}\mathbf{1}(\tau)\min\Big\{1,\frac{1-\beta_{\chi}}{\eta}\Big\}\Big],
\end{align*}
where $A_1=\exp(K/(300\eta))$, $A_2=\exp(40K/\eta)$, and
\[
\mathbf{1}(\tau)=\begin{cases}
	1&\mbox{if $|\tau|\leq 200\eta$,}\\
	0&\mbox{otherwise.}
\end{cases}
\]
We follow the usual convention of dropping terms involving $\beta_{\chi}$ if $\beta_{\chi}$ does not exist.
\end{prop}

When $\pi,\pi_0\in\mathcal{A}(1)$ and $\pi_0$ is trivial, \cref{prop:zero_detection} reduces to a result of Weiss \cite[Proposition 4.2]{Weiss}; we follow Weiss's proof with the modifications which follow \cite{RJLOT,ST} to allow for more general choices of $\pi$ and $\pi_0$.  Relative to the ideas in \cite{RJLOT,ST,Weiss}, there are three novelties here.  First, we exploit the existence of an exceptional zero of a Dirichlet $L$-function in the zero-detection process for $L(s,\pi\times\pi_0)$, which generalizes \cite[Proposition 4.2]{Weiss}.  Second, we use \cref{hyp} for both $\pi$ and $\pi_0$ instead of assuming that at least one of $\pi$ and $\pi_0$ satisfies GRC as in \cite{RJLOT} so that, unlike the approach in \cite{ST}, the Dirichlet polynomial can be supported on primes.  Third, much like \cite[Section 4]{ST}, the proof here makes explicit some of the effective constants in \cite{RJLOT,Weiss}, which we believe makes the proof a bit easier to read.

\subsection{Preliminary estimates}

\begin{lem}
\label{lem:Lemma_3.5_RJLOT}
Let $\pi\in\mathcal{F}_m(Q)$, let $\pi_0\in\mathcal{A}( m_0 )$, and let $\chi$ be a primitive Dirichlet character modulo $q$. If $\eta>0$, then
	\[
	\sum_{n\geq 1}\frac{|\lambda_{\pi\times(\pi_0\otimes\chi)}(n)|\Lambda(n)}{n^{1+\eta}}\leq\frac{1}{\eta}+ \frac{m_0 m}{2}\log(C(\pi_0)Q)+O(( m_0 m)^2).
	\]
\end{lem}
\begin{proof}
Suppose $\pi\in\mathcal{A}(d)\cap\mathcal{F}_m(Q)$.  It follows from \eqref{eqn:Brumley_Cauchy_Schwarz_2} and the discussion in \cite{ST} which follows Lemma 2.3 that
\[
\sum_{n\geq 1}\frac{|\lambda_{\pi\times(\pi_0\otimes\chi)}(n)|\Lambda(n)}{n^{1+\eta}}\leq\frac{1}{\eta}+ \frac{1}{4}\log C(\pi\times\widetilde{\pi})+\frac{1}{4}\log C(\pi_0\times\widetilde{\pi}_0)+O((d m_0)^2),
\]
The desired result follows from \eqref{eqn:BH}, the bound $C(\pi)\leq Q$, and the bound $d\leq m$.
\end{proof}

%
%
%
%
%
%

\begin{lem}
\label{lem:Lemma_3.4_RJLOT}
	Let $\pi\in\mathcal{F}_m(Q)$, let $\pi_0\in\mathcal{A}( m_0 )$, and let $\chi$ be a primitive Dirichlet character modulo $q\leq 2Q$.  Let $n(\eta;s)$ denote the number of zeros $\rho$ of $L(s,\pi\times(\pi_0\otimes\chi))$ with $|s-\rho|\leq\eta$.  For all $\re(s)\geq 1$ and all $0<\eta<1/2$, we have the bound
	\[
	n(\eta;s)\leq 20 (m m_0)^2 \eta \log(C(\pi_0)Q) + 5m_0 m\eta\log(|\im(s)|+2)+O((m_0 m)^2).
	\]
\end{lem}
\begin{proof}
	It suffices to prove the result for $n(\eta;1+it)$ because $n(\eta;1+it)\geq n(\eta;\sigma+it)$ for any $\sigma\geq 1$.  For $\pi\in\mathcal{A}(d)\cap\mathcal{F}_m(Q)$, it follows from \cite[Lemma 3.1]{ST} that for such $s$,
	\[
	n(\eta;s)\leq 10  d m_0 \eta \log C(\pi\times(\pi_0\otimes\chi)) + 5  d m_0\eta\log(|\im(s)|+2)+O(( d m_0)^2).
	\]
	The result now follows from \eqref{eqn:BH} and the hypothesis that $q\leq 2Q$.
\end{proof}

\begin{lem}
	\label{lem:hyp}
	If $\pi\in\mathcal{A}(d)$ satisfies \cref{hyp}, $y>C(\pi)$, and $\eta$ is as in \eqref{eqn:eta_def}, then
	\begin{equation}
	\label{eqn:the_original_bound}
	\sum_{\substack{n\in[y,y^{12000}] \\ \textup{$n$ composite}}}\frac{\lambda_{\pi\times\widetilde{\pi}}(n)\Lambda(n)}{n^{1+\eta}}\ll_{d}y^{-\frac{1}{2(d^2+1)}-\eta}(\log y)^3.
	\end{equation}
\end{lem}
\begin{proof}
We first bound the contribution to the sum in \eqref{eqn:the_original_bound} from the $n$ which share a prime factor with $N_{\pi}$ separately.  Note that $O(\log y)$ primes divide $N_{\pi}$ as $y > C(\pi) \geq N_{\pi}$.  Thus by  \eqref{eqn:LRS_2} and  \eqref{eqn:lambda_def} applied to the ramified prime, we have
\begin{align*}
\sum_{\substack{n\in[y,y^{12000}] \\ \textup{$n$ composite} \\ (n,N_{\pi})>1}}\frac{\lambda_{\pi\times\widetilde{\pi}}(n)\Lambda(n)}{n^{1+\eta}}&\ll d^2(\log y)\sum_{2\leq r\leq 20000\log y}~\sum_{\substack{y^{1/r}\leq p\leq y^{12000/r} \\ p\mid N_{\pi}}}p^{-r\frac{2}{d^2+1}-r\eta}\\
&\ll d^2 (\log y)^2\sum_{2\leq r\leq 20000\log y}y^{-\frac{2}{d^2+1}-\eta}\ll d^2 y^{-\frac{2}{d^2+1}-\eta}(\log y)^3.
\end{align*}

If $p\nmid N_{\pi}$, then $\lambda_{\pi\times\widetilde{\pi}}(p^r)=|\lambda_{\pi}(p^r)|^2$.  From \eqref{eqn:lambda_def}, we see that
\[
|\lambda_{\pi}(p^r)|^2\leq d^2\max_{1\leq j\leq d}|\alpha_{j,\pi}(p)|^{2r}.
\]
Define
\[
\beta_p = p^{-1}\max_{1\leq j\leq d}|\alpha_{j,\pi}(p)|^2.
\]
Note that $\beta_p\leq p^{-2/(d^2+1)}$ by \eqref{eqn:LRS_finite}.  Thus the contribution to the sum in \eqref{eqn:the_original_bound} arising from the $n$ which are coprime to $N_{\pi}$ is
\[
\ll d^2(\log y)\sum_{r=2}^{\infty}~\sum_{y^{1/r}\leq p\leq y^{12000/r}}\beta_p^r p^{-r\eta}\ll d^2(\log y)\sum_{2\leq R\leq 20000\log y}~\sum_{y^{1/R}<p\leq y^{12000/R}}~\sum_{r=R}^{\infty}\beta_p^r p^{-r\eta},
\]
Subject to \cref{hyp}, we will prove that
\begin{equation}
\label{eqn:original_1}
S_R:=\sum_{y^{1/R}<p\leq y^{12000/R}}~\sum_{r=R}^{\infty}\beta_p^r p^{-r\eta}\ll_d y^{-\frac{1}{2(d^2+1)}-\eta}(\log y)
\end{equation}
uniformly for all $2\leq R\leq 20000\log y$, which suffices to prove the lemma.

The inner sum is geometric, so
\begin{align*}
S_R= \sum_{y^{1/R}<p\leq y^{12000/R}} \frac{(\beta_p p^{-\eta})^R}{1-\beta_p p^{-\eta}}\leq y^{-\frac{1}{d^2+1}-\eta}\sum_{y^{1/R}<p\leq y^{12000/R}}\frac{\beta_p}{1-\beta_p p^{-\eta}}.
\end{align*}
We decompose the sum according to whether $p$ is greater than $2^{d^2}$ (in which case $1-\beta_p p^{-\eta}\geq 1/2$) or not.  The contribution from the latter range to the sum is $O_d(1)$, so we have
\begin{align}
\label{eqn:original_22}
S_R &\leq y^{-\frac{1}{d^2+1}-\eta}\Big(2\sum_{y^{1/R}<p\leq y^{12000/R}}\beta_p+O_d(1)\Big)\notag\\
&\leq y^{-\frac{1}{d^2+1}-\eta}\Big(2y^{\frac{1}{2(d^2+1)}}\sum_{y^{1/R}<p\leq y^{12000/R}}\beta_p p^{-\frac{1}{12000(d^2+1)}}+O_d(1)\Big).
\end{align}
Note that $x\leq 2\log(x+1)$ for all $0\leq x\leq 5/2$. Thus, since $0\leq\beta_p p^{-\frac{1}{12000(d^2+1)}}\leq 1$ for all $p$, the final sum in \eqref{eqn:original_22} is bounded by
\begin{align*}
2\sum_{y^{1/R}<p\leq y^{12000/R}}\log(1+\beta_p p^{-\frac{1}{12000(d^2+1)}})&\leq 2\sum_{p}\log\Big(\sum_{r=0}^{\infty}\beta_p^r p^{-\frac{r}{12000(d^2+1)}}\Big)\\
&=2\log\Big(\prod_{p}\sum_{r=0}^{\infty}\frac{\max_{1\leq j\leq d}|\alpha_{j,\pi}(p)|^{2r}}{p^{r(1+\frac{1}{12000(d^2+1)})}}\Big)
\end{align*}
The above display is of size $2\epsilon\log y + O_d(1)$ by \cref{hyp}, which establishes \eqref{eqn:original_1}.
\end{proof}

\subsection{Proof of \cref{prop:zero_detection}}

Let $\pi\in\mathcal{A}(d)\cap\mathcal{F}_m(Q)$ and $\pi_0\in\mathcal{A}(m_0)$ and let $\chi$ be a primitive real character modulo $q\leq 2Q$.  Suppose that $L(s,\pi\times\pi_0)$ has a zero $\rho_0\neq\beta_{\chi}$ such that $|\rho_0-(1+i\tau)|\leq\eta$, where $\tau\in\R$, $|\tau|\leq T$, and
\begin{equation}
\label{eqn:eta_def1}
\frac{1}{\log(C(\pi_0)Q T)}\leq\eta\leq\frac{1}{10^7(m_0 m)^2}.
\end{equation}
We let $s=1+\eta+i\tau$.  We define
\[
F(z):=\frac{L'}{L}(z,\pi\times\pi_0)+\frac{L'}{L}(z+1-\beta_{\chi},\pi\times(\pi_0\otimes\chi))
\]
and
\[
G_k(z) := \frac{(-1)^k}{k!}F^{(k)}(z)+\frac{r_{\pi\times\pi_0}}{(z-1)^{(k+1)}}-\frac{r_{\pi\times(\pi_0\otimes\chi)}}{(z+1-2\beta_{\chi})^{k+1}}.
\]
We apply \cite[Equation 4.2]{ST} and \eqref{eqn:BH} to obtain
\begin{align}
\label{eqn:power_sum}
G_k(s)&=\sum_{\substack{L(\rho,\pi\times\pi_0)=0 \\ \rho\neq\beta_{\chi} \\ |s-\rho|\leq 200\eta}}\frac{1}{(s-\rho)^{k+1}}+\sum_{\substack{L(\rho,\pi\times(\pi_0\otimes\chi))=0 \\ \rho\neq\beta_{\chi} \\ |s-\rho|\leq 200\eta}}\frac{1}{(s-\rho)^{k+1}}\notag\\
&+O\Big(\frac{d m_0\eta(\log(C(\pi\times\pi_0)+\log C(\pi\times(\pi_0\otimes\chi)))}{(200\eta)^k}\Big)\notag\\
&=\sum_{\substack{L(\rho,\pi\times\pi_0)=0 \\ \rho\neq\beta_{\chi} \\ |s-\rho|\leq 200\eta}}\frac{1}{(s-\rho)^{k+1}}+\sum_{\substack{L(\rho,\pi\times(\pi_0\otimes\chi))=0 \\ \rho\neq\beta_{\chi} \\ |s-\rho|\leq 200\eta}}\frac{1}{(s-\rho)^{k+1}}+O\Big(\frac{(m_0 m)^2\eta\log(C(\pi_0)QT)}{(200\eta)^k}\Big)
\end{align}
\cref{lem:Lemma_3.4_RJLOT} tells us that the two sums in \eqref{eqn:power_sum} have, in total, at most $K$ terms for any
\begin{align}
\label{eqn:K_def}
K \geq 8000  (m_0 m)^2\eta\log(C(\pi_0)QT)+O((m_0 m)^2).
\end{align}
Just like \cite{RJLOT,ST,Weiss}, we rely on the following diophantine result due to S{\'o}s and Tur{\'a}n \cite{Turan}.
\begin{lem}
\label{lem:turan}
Let $z_1,\ldots,z_{\nu}\in\mathbb{C}$.  If $K\geq \nu$, then there exists an integer $k\in[K,2K]$ such that $|z_1^k+\cdots+z_m^k|\geq(\frac{1}{50}|z_1|)^k$.
\end{lem}
If $L(s,\pi\times\pi_0)$ has a zero $\rho_0$ (trivial or nontrivial) satisfying $|\rho_0-(1+i\tau)|\leq\eta$ with $\rho_0\neq\beta_{\chi}$, we can apply \cref{lem:turan} to the sums over zeros in \eqref{eqn:power_sum} and find that if the implied constant in \eqref{eqn:K_def} is sufficiently large, then for some $k\in[K,2K]$,
\begin{align*}
\eta^{k+1}|G_k(s)|\geq\frac{\eta^{k+1}}{(50|s-\rho_0|)^{k+1}}-O\Big(\frac{ (m_0 m)^2\eta\log(C(\pi_0)QT)}{(200)^{k}}\Big)\geq\frac{3}{4(100)^{k+1}}.
\end{align*}
It follows from the calculations on \cite[pages 80-81]{Weiss} that
\[
\eta^{k+1}\Big|\frac{1}{(s-1)^{k+1}}-\frac{1}{(s+1-2\beta_{\chi})^{k+1}}\Big|\leq \frac{1}{4(100)^{k+1}}+\mathbf{1}(\tau)\min\{1,16\cdot 2^k((1-\beta_{\chi})/\eta)^{1/2}\},
\]
Therefore, for some $k\in[K,2K]$ with $K$ given by \eqref{eqn:K_def}, we have the lower bound
\begin{equation}
\label{eqn:part2_lower}
\frac{\eta^{k+1}}{k!}|F^{(k)}(s)|+r_{\pi\times\pi_0}\mathbf{1}(\tau)\min\{1,16\cdot 2^{k}((1-\beta_{\chi})/\eta)^{1/2}\}\geq\frac{1}{2(100)^{k+1}}.
\end{equation}

On the other hand, since $\eta>0$, we can use the absolute convergence of the Dirichlet series which defines $F(s)$ to directly compute
\begin{equation}
\label{eqn:sum_11}
\frac{\eta^{k+1}}{k!}|F^{(k)}(s)|=\eta\Big|\sum_{n\geq 1}\frac{(\lambda_{\pi\times\pi_0}(n)+\lambda_{\pi\times(\pi_0\otimes\chi)}(n)n^{\beta_{\chi}-1})\Lambda(n)}{n^{1+i\tau}}j_k(\eta\log n)\Big|,
\end{equation}
where $j_k(u)=(k!)^{-1}u^k e^{-u}$.  Let $A_1=\exp(K/(300\eta))$ and $A_2=\exp(40K/\eta)$. Suppressing summands, we write the right hand side of \eqref{eqn:sum_11} as
\begin{equation}
\label{eqn:summand_suppressed}
\eta\sum_{n\geq 1}=\eta\Big(\sum_{n\notin[A_1,A_2]}+\sum_{\substack{ n\in[A_1,A_2]\\ \textup{$n$ composite}}}+\sum_{p\in[A_1,A_2]}\Big)
\end{equation}

First, we bound the contribution from $n\notin[A_1,A_2]$.  Since $k!\geq(k/e)^k$, we find from a small numerical calculation \cite[Proof of Lemma 4.3]{ST} that
\begin{equation}
\label{eqn:j_bounds}
j_k(\eta\log n)\leq (110)^{-k}n^{-\eta/2}\qquad\textup{if $n\notin[A_1,A_2]$.}
\end{equation}  By \eqref{eqn:j_bounds},
\begin{align*}
\Big|\eta\sum_{n\notin[A_1,A_2]}\Big|&\ll\eta(110)^{-k}\sum_{n\geq 1}\frac{(|\lambda_{\pi\times\pi_0}(n)|+|\lambda_{\pi\times(\pi_0\otimes\chi)}(n)|)\Lambda(n)}{n^{1+\eta/2}}.
\end{align*}
By \cref{lem:Lemma_3.5_RJLOT}, the above display is $\ll\eta(110)^{-k}(\eta^{-1}+ (m_0 m)^2\log(C(\pi_0)qQT))$.  Using \eqref{eqn:K_def}, we see that the contribution from $n\notin[A_1,A_2]$ is $O_{m_0,m}(k(110)^{-k})$.

Second, we bound the contribution from the composite $n\in[A_1,A_2]$.  Since $(\log u)^k\leq k!u$ for all $k\geq1$ and $u\geq 1$, we find that
\[
j_k(\eta\log n)=\frac{(\eta\log n)^k}{k!n^{\eta}}=\frac{1}{n^{\eta}}(110)^{-k}\frac{(\log n^{110\eta})^k}{k!}\leq \frac{1}{n^{\eta}}(110)^{-k}n^{110\eta}.
\]
This estimate and \eqref{eqn:Brumley_Cauchy_Schwarz_2} imply that
\begin{align*}
\Big|\eta\sum_{\substack{n\in[A_1,A_2] \\ \textup{$n$ composite}}}\Big|&\leq \eta(110)^{-k}\sum_{\substack{n\in[A_1,A_1^{12000}] \\ \textup{$n$ composite}}}\frac{(|\lambda_{\pi\times\pi_0}(n)|+|\lambda_{\pi\times(\pi_0\otimes\chi)}(n)|)\Lambda(n)}{n^{1+\eta}}n^{110\eta}\\
&\leq \eta(110)^{-k}\sum_{\substack{n\in[A_1,A_1^{12000}] \\ \textup{$n$ composite}}}\frac{(\lambda_{\pi\times\widetilde{\pi}}(n)+\lambda_{\pi_0\times\widetilde{\pi}'}(n))\Lambda(n)}{n^{1+\eta}} n^{110\eta}\\
&\leq\eta(110)^{-k}A_1^{1320000\eta}\sum_{\substack{n\in[A_1,A_1^{12000}] \\ \textup{$n$ composite}}}\frac{(\lambda_{\pi\times\widetilde{\pi}}(n)+\lambda_{\pi_0\times\widetilde{\pi}'}(n))\Lambda(n)}{n^{1+\eta}}.
\end{align*}
By \cref{lem:hyp} and \eqref{eqn:eta_def1}, the above display is
\[
\ll_{d, m_0 }\eta(110)^{-k}A_1^{1320000\eta-\frac{1}{2(( d m_0)^2+1)}-\eta}(\log A_1)^2\ll_{d, m_0 } \eta(110)^{-k}\ll_{m, m_0 }k(110)^{-k}.
\]

Finally, we estimate the contribution from the primes $p\in[A_1,A_2]$.  Summation by parts gives us the identity
\begin{align}
\label{eqn:sum_by_parts_turan}
\eta\sum_{p\in[A_1,A_2]}&=j_k(\eta\log A_2)\eta\sum_{p\in[A_1,A_2]}\frac{(\lambda_{\pi\times\pi_0}(p)+\lambda_{\pi\times(\pi_0\otimes\chi)}(p)p^{\beta_{\chi}-1})\Lambda(p)}{p^{1+i\tau}}\notag\\
&-\eta^2\int_{A_1}^{A_2}j_k'(\eta\log u)\sum_{p\in[A_1,u]}\frac{(\lambda_{\pi\times\pi_0}(p)+\lambda_{\pi\times(\pi_0\otimes\chi)}(p)p^{\beta_{\chi}-1})\Lambda(p)}{p^{1+i\tau}}\frac{du}{u}.
\end{align}
Much like the above calculations, we use \cref{lem:Lemma_3.5_RJLOT} to deduce that the sum over $p\in[A_1,A_2]$ in \eqref{eqn:sum_by_parts_turan} is
\[
\ll \eta(110)^{-k}A_2^{-\eta/2}\sum_{n<A_2}\frac{(|\lambda_{\pi\times\pi_0}(n)|+|\lambda_{\pi\times(\pi_0\otimes\chi)}(n)|)\Lambda(n)}{n}\ll \frac{k}{(110)^{k}}.
\]
Since $|j_k'(u)|=|j_{k-1}(u)-j_k(u)|\leq j_{k-1}(u)+j_k(u)\leq 1$, we find that
\[
\Big|\eta\sum_{p\in[A_1,A_2]}\Big|\leq \eta^2\int_{A_1}^{A_2}\Big|\sum_{p\in[A_1,u]}\frac{(\lambda_{\pi\times\pi_0}(p)+\lambda_{\pi\times(\pi_0\otimes\chi)}(p)p^{\beta_{\chi}-1})\Lambda(p)}{p^{1+i\tau}}\Big|\frac{du}{u}+O\Big(\frac{k}{(110)^{k}}\Big).
\]

If $K$ is given by \eqref{eqn:K_def}, then the condition $p\in[A_1,A_2]$ implies that $p\nmid N_{\pi}N_{\pi_0}q$.  Therefore, by \eqref{eqn:separate_dirichlet_coeffs} and \eqref{eqn:lambda_def},
\[
(\lambda_{\pi\times\pi_0}(p)+\lambda_{\pi\times(\pi_0\otimes\chi)}(p)p^{\beta_{\chi}-1})\Lambda(p)= \lambda_{\pi\times\pi_0}(p)(1+\chi(p)p^{\beta_{\chi}-1})\log p.
\] 
We collect our estimates for the three sums in \eqref{eqn:summand_suppressed} to find that for all $k\in[K,2K]$ with $K$ given by \eqref{eqn:K_def1},
\begin{align}
\label{eqn:sum_111}
	\frac{\eta^{k+1}}{k!}|F^{(k)}(s)|\leq \eta^2\int_{A_1}^{A_2}\Big|\sum_{A_1<p\leq u}\frac{\lambda_{\pi\times\pi_0}(p)\log p}{p^{1+i\tau}}(1+\chi(p)p^{\beta_{\chi}-1})\Big|\frac{du}{u}+O_{m_0, m}\Big(\frac{k}{(110)^k}\Big).
\end{align}

We enlarge $K$ according to \eqref{eqn:K_def1}, which we are free to do.  If $k\in[K,2K]$ and the implied constant in \eqref{eqn:K_def1} is sufficiently large, then $O_{m, m_0 }(k(110)^{-k})\leq \frac{1}{4}(100)^{-k-1}$.  Therefore, it follows from \eqref{eqn:part2_lower} and \eqref{eqn:sum_111} that if $L(s,\pi\times\pi_0)$ has a zero $\rho_0\neq\beta_{\chi}$ which satisfies $|\rho_0-(1+i\tau)|\leq\eta$, then with $K$ given by \eqref{eqn:K_def1}, we have the bound
\begin{align*}
1&\leq 4(100)^{2K+1}\eta^2\int_{A_1}^{A_2}\Big|\sum_{A_1<p\leq u}\frac{\lambda_{\pi\times\pi_0}(p)\log p}{p^{1+i\tau}}(1+\chi(p)p^{\beta_{\chi}-1})\Big|\frac{du}{u}\\
&+4r_{\pi\times\pi_0}(100)^{2K+1}\mathbf{1}(\tau)\min\{1,16\cdot 2^{2K}((1-\beta_{\chi})/\eta)^{1/2}\}.
\end{align*}
We square both sides and apply the Cauchy-Schwarz inequality to obtain the bound
\begin{align*}
1&\ll (100)^{4K}\eta^4\Big(\int_{A_1}^{A_2}\frac{du}{u}\Big)\Big(\int_{A_1}^{A_2}\Big|\sum_{A_1<p\leq u}\frac{\lambda_{\pi\times\pi_0}(p)\log p}{p^{1+i\tau}}(1+\chi(p)p^{\beta_{\chi}-1})\Big|^2\frac{du}{u}\Big)\\
&+r_{\pi\times\pi_0}(100)^{4K}\mathbf{1}(\tau)\min\{1, 2^{4K}\eta^{-1}(1-\beta_{\chi})\}.
\end{align*}
Since $\int_{A_1}^{A_2}u^{-1}du\ll K/\eta$, \cref{prop:zero_detection} follows.

\section{A new large sieve inequality}


We will generalize the large sieve for Dirichlet coefficients of automorphic representations due to Duke and Kowalski \cite[Theorem 4]{DK}.  As observed by Brumley \cite{Brumley_2}, one can adjust their proof to show that if $\mathcal{F}_m(Q)$ satisfies \cref{hyp} and $Q,x\geq 2$, then
\begin{equation}
\label{eqn:DK_large_sieve}
\sum_{\pi\in\mathcal{F}_m(Q)}\Big|\sum_{n\leq x}a_{\pi}(n)b(n)\Big|^2\ll_{\epsilon,m}(Q^{\epsilon}x+Q^{\frac{3}{2m}}x^{1-\frac{1}{m^2}}\#\mathcal{F}_m(Q))\sum_{n\leq x}|b(n)|^2,
\end{equation}
where $b(n)$ is any complex-valued function supported on the integers.  We require two modifications to \eqref{eqn:DK_large_sieve}.  First, we need to take sums over $n$ in intervals of length $x/T$, where $T$ is arbitrarily large.  Second, we need a variant of \eqref{eqn:DK_large_sieve} which applies with more sensitivity to sequences $b(n)$ supported on the primes.

We establish a ``pre-sifted'' large sieve inequality over short intervals for families of automorphic representations which satisfy \cref{hyp}.  We anticipate that this will be useful in contexts beyond this paper.  In what follows, we define $P^{-}(n)$ to be the least prime dividing a positive integer $n$; we set $P^{-}(1)=\infty$ by convention.

\begin{prop}
\label{thm:large_sieve}
Let $b(n)$ be a complex-valued function supported on the integers, and suppose that each $\pi\in\mathcal{F}_m(Q)$ (see  \eqref{eqn:family}) satisfies \cref{hyp}.  If $Q\geq 3$, $T\geq1$, $x>0$, and $z\gg_m Q^{6m}$ with a sufficiently large implied constant, then for every $\epsilon>0$,
\begin{align*}
\sum_{\pi\in\mathcal{F}_m(Q)}\Big|\sum_{\substack{x<n\leq x e^{1/T}\\ P^{-}(n)>z}}a_{\pi}(n)b(n)\Big|^2\ll_{\epsilon,m}\Big(\frac{x}{T\log z}+Q^{\frac{3}{2m}}T^{\frac{3}{4}}x^{1-\frac{1}{m^2}}z^{2+\epsilon}\#\mathcal{F}_m(Q)\Big)\sum_{\substack{x<n\leq x e^{1/T}  \\ P^{-}(n)>z }}|b(n)|^2.
\end{align*}
\end{prop}

\subsection{The na{\"i}ve Rankin-Selberg $L$-function}
\label{subsec:naive}

Let $\pi\in\mathcal{A}(d)$ and $\pi'\in\mathcal{A}(d')$.  For each prime $p$, define
\begin{equation}
\label{eqn:square-free_L-factor}	
L^{RS}(s,\pi_p\times\pi_{p}')=1+\sum_{j=1}^{\infty}\frac{a_{\pi}(p^j)a_{\pi'}(p^j)}{p^{js}}.
\end{equation}
We call the Dirichlet series
\begin{equation}
\label{eqn:square-free_L-function}
L^{RS}(s,\pi\times\pi'):=\sum_{\substack{n\geq 1 \\ (n,N_{\pi}N_{\pi'})=1}}\frac{a_{\pi}(n)a_{\pi'}(n)}{n^s}=\prod_{p\nmid N_{\pi}N_{\pi'}}L^{RS}(s,\pi_p\times\pi_{p}')
\end{equation}
the na{\"i}ve Rankin-Selberg $L$-function.  We access the Dirichlet coefficients of $L^{RS}(s,\pi\times\pi')$ by relating $L^{RS}(s,\pi\times\pi')$ to $L(s,\pi\times\pi')$.  In order to accomplish this, we use \cref{hyp} and the following result of Brumley (see the proof of \cite[Corollary 3]{Brumley_2}).

\begin{lem}[Brumley]
\label{cor:hyp}
	Suppose that  $\pi, \pi'\in\mathcal{F}_m(Q)$ satisfy \cref{hyp}.  For each prime $p$, define $H(s,\pi_p\times\pi_{p})$ by the equality
	\[
	L^{RS}(s,\pi_p\times\pi_{p}')=L(s,\pi_p\times\pi_{p}')H(s,\pi_p\times\pi_{p}').
	\]
	The Euler product $H(s,\pi\times\pi'):=\prod_{p\nmid N_{\pi}N_{\pi'}} H(s,\pi_p\times\pi_{p})$ converges absolutely for $\re(s)\geq 1-m^{-2}$; this yields the factorization
	\[
	L^{RS}(s,\pi\times\pi')=L(s,\pi\times\pi')H(s,\pi\times\pi')\prod_{p\mid N_{\pi}N_{\pi'}}L(s,\pi_p\times\pi_{p}')^{-1}
	\]
	in the region $\re(s)\geq 1-m^{-2}$.  Furthermore, $H(s,\pi\times\pi')\ll_{\epsilon,m}Q^{\epsilon}$ for $\re(s)\geq 1-m^{-2}$.
\end{lem}

\subsection{Preliminary estimates}


Let $\pi,\pi'\in\mathcal{F}_m(Q)$, and assume throughout this subsection that both $\pi$ and $\pi'$ satisfy \cref{hyp}.  Let
\[
g^{RS}_d(s,\pi\times\widetilde{\pi}'):=\prod_{p\mid d}(1-L^{RS}(s,\pi_p\times\widetilde{\pi}_{p}')^{-1}),
\]
and let $d\geq 1$ be a square-free integer which is coprime to $N_{\pi}N_{\pi'}$.  We will consider the Dirichlet series given by
\[
L_{d}^{RS}(s,\pi\times\widetilde{\pi}'):=\sum_{\substack{n\geq1 \\ d\mid  n \\ (n,N_{\pi}N_{\pi'})=1}}\frac{a_{\pi}(n)a_{\widetilde{\pi}'}( n)}{n^s}=L^{RS}(s,\pi\times\widetilde{\pi}')g_{d}^{RS}(s,\pi\times\widetilde{\pi}').
\]
A bound for $g^{RS}_d(s,\pi\times\widetilde{\pi}')$ follows readily from \eqref{eqn:LRS_2}.

\begin{lem}
\label{lem:g}
Let $d\geq 1$ be square-free and $\pi,\pi'\in\mathcal{F}_m(Q)$.  In the region $\sigma\geq 1-m^{-2}$, we have that $g^{RS}_d(s,\pi\times\widetilde{\pi}')\ll_{\epsilon,m}d^{\epsilon}$.  If $d\geq 2$, then $0\leq g^{RS}_d(1,\pi\times\widetilde{\pi})<1$.
\end{lem}
\begin{proof}
The fact that $0\leq g^{RS}_d(1,\pi\times\widetilde{\pi})<1$ for $d\geq 2$ follows immediately from \eqref{eqn:square-free_L-factor}.  The bound \eqref{eqn:LRS_2} yields $|L^{RS}(s,\pi_p\times\widetilde{\pi}'_p)^{-1}|\ll_m 1$ for $\re(s)\geq 1-m^{-2}$.  The lemma now follows from the well-known bound $\omega(n)\ll(\log\log n)^{-1}\log n$, where $\omega(n)$ is the number of distinct prime factors of $n$.
\end{proof}

We require some uniform estimates for $L_d^{RS}(s,\pi\times\widetilde{\pi}')$.
\begin{lem}
\label{lem:H}
Let $s=\sigma+it$, $\pi,\pi'\in\mathcal{F}_m(Q)$, and . For any squarefree integer $d \geq 1$ coprime to $N_{\pi}N_{\pi'}$, 
\[
|(\sigma-1)^{r(\pi\times\widetilde{\pi}')}L_{d}^{RS}(s,\pi\times\widetilde{\pi}')|,\ll_{\epsilon,m} d^{\epsilon}Q^{\frac{3}{2m}}(1+|t|)^{\frac{3}{4}}
\]
uniformly in the region $\sigma\geq 1-m^{-2}$.
\end{lem}
\begin{proof}
First, we establish the bound
\begin{equation}
\label{eqn:intermediate}
|(\sigma-1)^{r_{\pi\times\widetilde{\pi}'}}L(s,\pi\times\widetilde{\pi}')|\ll_{\epsilon,m}(Q^{2m}(1+|t|)^{m^2})^{\max\{\frac{1}{2}(1-\sigma),0\}+\epsilon},\qquad 1/2\leq\sigma\leq 3
\end{equation}
for every $\epsilon>0$. By the work of Li \cite[Theorem 2]{Li}, we know that for some constant $c_m>0$ depending at most on $m$,
\begin{equation}
\label{thm:Li_2}
	(\sigma-1)^{r_{\pi\times\widetilde{\pi}'}}|L(\sigma,\pi\times\widetilde{\pi}')|\ll\exp\Big(c_m\frac{\log C(\pi\times\widetilde{\pi}')}{\log\log C(\pi\times\widetilde{\pi}')}\Big)\ll_{\epsilon,m} C(\pi\times\widetilde{\pi}')^{\epsilon},\quad 1\leq\sigma\leq 3.
\end{equation}
By replacing $\widetilde{\pi}'$ with $\widetilde{\pi}'\otimes|\det|^{-it}$ in the proof of \eqref{thm:Li_2} (which does not change the proof substantially), we obtain
\begin{equation}
\label{thm:Li}
|(\sigma-1)^{r_{\pi\times\widetilde{\pi}'}}L(\sigma+it,\pi\times\widetilde{\pi}')|\ll_{\epsilon,m}C(\pi\times\widetilde{\pi}',t)^{\epsilon},\qquad 1\leq\sigma\leq 3.
\end{equation}
The refined version of the convexity bound for $L$-functions proved by Heath-Brown in \cite{HB} yields
\begin{equation}
\label{eqn:HB_1}
|L(1/2+it,\pi\times\widetilde{\pi}')|\ll_{m}|L(3/2+it,\pi\times\widetilde{\pi}')|^2 C(\pi\times\widetilde{\pi}',t)^{1/4}.
\end{equation}
Hence, by \eqref{thm:Li},
\begin{equation}
\label{eqn:HB_2}
|L(1/2+it,\pi\times\pi')|\ll_{\epsilon,m}C(\pi\times\pi',t)^{1/4+\epsilon}.	
\end{equation}
Thus \eqref{eqn:intermediate} follows from \eqref{eqn:BH}, \eqref{thm:Li},  \eqref{eqn:HB_2}, and an application of the Phragm{\'e}n-Lindel{\"o}f principle.

We see from \eqref{eqn:LRS_2} and the bound $\omega(n)\ll(\log\log n)^{-1}\log n$ that for every $\epsilon>0$, one has the bound
\[
\prod_{p\mid N_{\pi}N_{\pi'}}|L(s,\pi_p\times\widetilde{\pi}_p')^{-1}|\ll_{\epsilon,m}Q^{\epsilon},\qquad \re(s)\geq 1-m^{-2}.
\]
With this bound in hand, the lemma follows from \cref{cor:hyp}, \cref{lem:g}, and \eqref{eqn:intermediate}.
\end{proof}

Fix a smooth function $\phi$ whose support is a compact subset of $(-2,2)$.  Let
\[
\widehat{\phi}(s)=\int_{-\infty}^{\infty}\phi(y)e^{sy}dy.
\]
Thus $\widehat{\phi}(s)$ is entire, and integrating by parts several times yields the bound
\begin{equation}
\label{eqn:Laplace}
\widehat{\phi}(s)\ll_{\phi,k}\frac{e^{2|\re(s)|}}{|s|^k}.
\end{equation}
any integer $k\geq0$.  Let $T\geq 1$; by Fourier inversion, one has the identity
\[
T\phi(T\log x)=\frac{1}{2\pi i}\int_{c-i\infty}^{c+i\infty}\widehat{\phi}(s/T)x^{-s}ds
\]
for any $x>0$ and any $c\in\mathbb{R}$.

\begin{lem}
\label{lem:smooth_positive_sums}
Let $\pi,\pi'\in\mathcal{F}_m(Q)$ with $m\geq 2$.  Let $x>0$, $T\geq 1$, and $d\geq 1$ be a square-free integer which is coprime to $N_{\pi}N_{\pi'}$.
\begin{enumerate}
\item If $\phi$ is as above, then
\begin{align*}
\Big|\sum_{\substack{n\geq 1 \\ d\mid n \\ (n,N_{\pi}N_{\pi'})=1}}a_{\pi}(n)&a_{\widetilde{\pi}'}(n)\phi\Big(T\log\frac{n}{x}\Big)\\
&-g_d^{RS}(1,\pi\times\widetilde{\pi}')\kappa_{\pi\times\widetilde{\pi}'}H(1,\pi\times\widetilde{\pi}')\frac{x\widehat{\phi}(1/T)}{T}\prod_{p\mid N_{\pi}N_{\widetilde{\pi}'}}L^{RS}(1,\pi_p\times\widetilde{\pi}_{p}')^{-1}\Big|\\
&\qquad\qquad\qquad\qquad\qquad\qquad\qquad\qquad\ll_{\epsilon,m,\phi}d^{\epsilon}Q^{\frac{3}{2m}}T^{\frac{3}{4}}x^{1-\frac{1}{m^2}}.
\end{align*}
\item If $z\gg_{m} Q^{6m}$ with a sufficiently large implied constant, then
\[
\sum_{\substack{n\leq z \\ (n,N_{\pi})=1}}\frac{|a_{\pi}(n)|^2}{n}\geq \frac{\log z}{20}\kappa_{\pi\times\widetilde{\pi}}H(1,\pi\times\widetilde{\pi})\prod_{p\mid N_{\pi}}L^{RS}(1,\pi_p\times\widetilde{\pi}_{p})^{-1}+\frac{1}{2}.
\]
\end{enumerate}
In both results, the quantities involving $\pi \times \widetilde{\pi}'$ are   positive when $\pi=\pi'$. Otherwise, $\kappa_{\pi \times \widetilde{\pi}'} = 0$ whenever $\pi \neq \pi'$. 
\end{lem}
\begin{proof}
For Part 1, the quantity we want to estimate equals, by \cref{cor:hyp} and the above properties of $\phi$,
	\begin{align*}
\frac{1}{2\pi iT}\int_{1-m^{-2}-i\infty}^{1-m^{-2}+i\infty}L_d^{RS}(s,\pi\times\widetilde{\pi}')\widehat{\phi}(s/T)x^s ds.
	\end{align*}
	By \cref{lem:H} and \eqref{eqn:Laplace}, the integral in the above display is
	\begin{align*}
	&\ll_{\epsilon,m}\frac{x^{1-\frac{1}{m^2}}d^{\epsilon}Q^{\frac{3}{2m}}}{T}\int_{-\infty}^{\infty}\Big|\widehat{\phi}\Big(\frac{1-m^{-2}+it}{T}\Big)\Big|(1+|t|)^{\frac{3}{4}}dt\\
	&\ll_{\epsilon,m,\phi}\frac{x^{1-\frac{1}{m^2}}d^{\epsilon}Q^{\frac{3}{2m}}}{T}\int_{-\infty}^{\infty}\min\Big\{1,\frac{T^2}{(|t|+2)^2}\Big\}(1+|t|)^{\frac{3}{4}}dt\ll_{\epsilon,m,\phi}d^{\epsilon}Q^{\frac{3}{2m}}T^{\frac{3}{4}}x^{1-\frac{1}{m^2}}.
	\end{align*}

We proceed to Part 2.  Let
\[
\phi(t)=\begin{cases}
	\exp(16+t^{-1}(t+\frac{1}{2})^{-1})&\mbox{if $t\in(-\frac{1}{2},0)$,}\\
	0&\mbox{otherwise.}
\end{cases}
\]
Observe that if $z\geq 4$, then by \cref{cor:hyp}, \cref{lem:H}, and \eqref{eqn:Laplace},
\[
\sum_{\substack{n\geq 1 \\ (n,N_{\pi})=1}}\frac{|a_{\pi}(n)|^2}{n}\phi\Big(\log\frac{n}{z}\Big)-\kappa_{\pi\times\widetilde{\pi}}H(1,\pi\times\widetilde{\pi})\widehat{\phi}(0)\prod_{p\mid N_{\pi}N_{\widetilde{\pi}'}}L(1,\pi_p\times\widetilde{\pi}_{p})^{-1}
\]
equals
\[
\frac{1}{2\pi i}\int_{-m^{-2}-i\infty}^{-m^{-2}+i\infty}L_d^{RS}(s+1,\pi\times\widetilde{\pi})\widehat{\phi}(s)z^s ds \ll_{m}Q^{\frac{3}{2m}}z^{-\frac{1}{m^2}}.
\]
The intervals $[2^{-j}e^{-1/2}z,2^{-j}z]$ and $[2^{-j-1} e^{-1/2}z,2^{-j-1}z]$ are disjoint for all integers $0\leq j\leq \frac{\log z}{\log 4}$, so
\begin{align*}
\sum_{\substack{n\leq z \\ (n,N_{\pi})=1}}\frac{|a_{\pi}(n)|^2}{n}&\geq1+\sum_{j\leq \frac{\log z}{\log 4}}\sum_{\substack{n\geq 1 \\ (n,N_{\pi})=1}}\frac{|a_{\pi}(n)|^2}{n}\phi\Big(\log\frac{n}{z/2^j}\Big)\\
&=1+\left\lfloor\frac{\log z}{\log 4}\right\rfloor\kappa_{\pi\times\widetilde{\pi}}H(1,\pi\times\widetilde{\pi})\widehat{\phi}(0)\prod_{p\mid N_{\pi}}L^{RS}(1,\pi_p\times\widetilde{\pi}_p)^{-1}+O_{m}(Q^{\frac{3}{2m}}z^{-\frac{1}{2m^2}}).
\end{align*}
Since $\widehat{\phi}(0)\geq 1/10$, the result follows once $z\gg_{m} Q^{6m}$.

For Parts 1 and 2, note that $\kappa_{\pi\times\widetilde{\pi}'}>0$ if and only if $\pi=\pi'$.  The same holds for $H(1,\pi\times\widetilde{\pi}')$ by appealing to the Euler product definition of $H$ in \cref{cor:hyp} and fact that $a_{\pi\times\widetilde{\pi}}(n)\geq0$ for all $n\geq 1$ (see \cite[Lemma a]{Hoffstein}).   From \eqref{eqn:square-free_L-factor} and the fact that $\bar{a}_{\pi}(n)=a_{\widetilde{\pi}}(n)$, we have that $\prod_{p\mid N_{\pi}}L^{RS}(1,\pi_p\times\widetilde{\pi}_p)^{-1}>0$, and the lemma follows.
	\end{proof}

\subsection{Proof of \cref{thm:large_sieve}}

It suffices to prove the bound
\begin{align}
\label{eqn:dual_sum}
\sum_{\substack{x<n\leq xe^{1/T} \\ P^{-}(n)>z}}\Big|\sum_{\pi\in\mathcal{F}_m(Q)}a_{\pi}(n)b_{\pi}\Big|^2\ll_{\epsilon,m} \Big(\frac{x/T}{\log z}+T^{\frac{3}{4}}Q^{\frac{3}{2m}}x^{1-\frac{1}{m^2}}z^{2+\epsilon}\#\mathcal{F}_m(Q)\Big)\sum_{\pi\in\mathcal{F}_m(Q)}|b_{\pi}|^2
\end{align}
for any sequence of complex numbers $\{b_{\pi}\}_{\pi\in\mathcal{F}_m(Q)}$ with the convention that $a_{\pi}(n)=0$ when $(n,N_{\pi})>1$, where $Q\geq 3$, $T\geq 1$, $x>0$, and $z\gg_m Q^{6m}$.  Indeed, with \eqref{eqn:dual_sum} in hand, it follows from a standard application of the duality principle that
\begin{align*}
\sum_{\pi\in\mathcal{F}_m(Q)}\Big|\sum_{\substack{x<n\leq x e^{1/T}\\ P^{-}(n)>z}}a_{\pi}(n)b(n)\Big|^2\ll_{\epsilon,m}\Big(\frac{x}{T\log z}+Q^{\frac{3}{2m}}T^{\frac{3}{4}}x^{1-\frac{1}{m^2}}z^{2+\epsilon}\#\mathcal{F}_m(Q)\Big)\sum_{\substack{x<n\leq x e^{1/T}  \\ P^{-}(n)>z }}|b(n)|^2,
\end{align*}
again with the convention that $a_{\pi}(n)=0$ when $(n,N_{\pi})>1$ and with $Q$, $T$, $x$, and $z$ as before.  Since $z>Q$ by hypothesis and $Q\geq N_{\pi}$ for all $\pi\in\mathcal{F}_m(Q)$, the condition $P^{-}(n)>z$ implies that $(n,N_{\pi})=1$ for all $n\in(x,xe^{1/T}]$.  The proposition now follows.

To bound \eqref{eqn:dual_sum}, we choose a compactly supported, infinitely differentiable function $\phi$ such that $\phi(t)\geq 1$ for $t\in[0,1]$ and $\phi(t)\geq0$ otherwise.  Then $\phi(T\log\frac{n}{x})$ is a pointwise upper bound for the indicator function of the interval $(x,xe^{1/T}]$.  If $w_z$ is any function such that $w_z(n)\geq 1$ if $P^{-}(n)>z$ and $w_z(n)\geq 0$ otherwise, then
\begin{equation}
\label{eqn:dual_sum_2.0}
\sum_{\substack{x<n\leq xe^{1/T} \\ P^{-}(n)>z}}\Big|\sum_{\pi\in\mathcal{F}_m(Q)}a_{\pi}(n)b_{\pi}\Big|^2\leq\sum_{n\geq1}\Big|\sum_{\pi\in\mathcal{F}_m(Q)}a_{\pi}(n)b_{\pi}\Big|^2  w_{z}(n)\phi\Big(T\log\frac{n}{x}\Big).
\end{equation}
We expand the square, swap the order of summation, and apply \cite[Lemma 1]{DK} so that the righthand side of \eqref{eqn:dual_sum_2.0} equals
\begin{align}
\label{eqn:DK_3}
&\sum_{\pi,\pi'\in\mathcal{F}_m(Q)}b_{\pi}\overline{b_{\pi'}}\sum_{n\geq 1}a_{\pi}(n)\overline{a_{\pi'}(n)} w_{z}(n)\phi\Big(T\log\frac{n}{x}\Big)\notag\\
&\leq \Big(\max_{\pi\in\mathcal{F}_m(Q)}\sum_{\pi'\in\mathcal{F}_m(Q)}\Big|\sum_{n\geq1}a_{\pi}(n)a_{\widetilde{\pi}'}(n) w_{z}(n)\phi\Big(T\log\frac{n}{x}\Big)\Big|\Big)\sum_{\pi\in\mathcal{F}_m(Q)}|b_{\pi}|^2.
\end{align}

We now choose $ w_{z}(n)$ as in the Selberg sieve.  Let $\pi_1\in\mathcal{F}_m(Q)$ be a representation which achieves the maximum in \eqref{eqn:DK_3}.  Let $g(d)=g_d(1,\pi_1\times\widetilde{\pi}_1)$, and define
\[
P(z)=\prod_{\substack{p<z \\ g(p)\neq 0 \\ p\nmid N_{\pi_1}}}p,\qquad \mathcal{D}_z = \{d\colon d\leq z,~d\mid P(z)\}.
\]
Let $\rho_{d}$ be a real-valued function such that
\begin{enumerate}
	\item $\rho_{1}=1$,
	\item $\rho_{d}=0$ unless $d\in\mathcal{D}_z$,
	\item $|\rho_{d}|\leq 1$ for all $d$.
\end{enumerate}
For integers $a,b\geq1$, let $[a,b]$ and $(a,b)$ denote the least common multiple and greatest common divisor of $a$ and $b$, respectively.  Conditions (1) and (2) in our choice of $\rho_{d}$ imply that if $P^{-}(n)>z$, then the statement $d\mid n$ implies that $d=1$ or $\rho_{d}=0$.  Therefore, we choose
\[
 w_{z}(n) = \Big(\sum_{d\mid(n,P(z))}\rho_{d}\Big)^2.
\]
Upon expanding the square and swapping the order of summation, \eqref{eqn:DK_3} equals
\begin{equation*}
\label{eqn:DK_4}
\sum_{\pi'\in\mathcal{F}_m(Q)}\Big|\sum_{d_1,d_2\in\mathcal{D}_z}\rho_{d_1}\rho_{d_2}\sum_{\substack{n\geq 1 \\ [d_1,d_2]\mid n} }a_{\pi_1}(n)a_{\widetilde{\pi}'}(n)\phi\Big(T\log\frac{n}{x}\Big)\Big|\sum_{\pi\in\mathcal{F}_m(Q)}|b_{\pi}|^2.
\end{equation*}
Our convention that $a_{\pi}(n)=0$ if $(n,N_{\pi})>1$ for each $\pi\in\mathcal{F}_m(Q)$ means that the above display is bounded by
\begin{align}
\label{eqn:DK_4}
&\Big|\sum_{d_1,d_2\in\mathcal{D}_z}\rho_{d_1}\rho_{d_2}\sum_{\substack{n\geq 1 \\ [d_1,d_2]\mid n \\ (n,N_{\pi_1})=1} }|a_{\pi_1}(n)|^2\phi\Big(T\log\frac{n}{x}\Big)\Big|\sum_{\pi\in\mathcal{F}_m(Q)}|b_{\pi}|^2\notag\\
&+\max_{\substack{\pi'\in\mathcal{F}_m(Q) \\ \pi_1\not=\pi'}}\Big|\sum_{d_1,d_2\in\mathcal{D}_z}\rho_{d_1}\rho_{d_2}\sum_{\substack{n\geq 1 \\ [d_1,d_2]\mid n \\ (n,N_{\pi_1}N_{\pi'})=1} }a_{\pi_1}(n)a_{\widetilde{\pi}'}(n)\phi\Big(T\log\frac{n}{x}\Big)\Big|\#\mathcal{F}_m(Q)\sum_{\pi\in\mathcal{F}_m(Q)}|b_{\pi}|^2.
\end{align}
We use \cref{lem:smooth_positive_sums} along with condition (3) to bound \eqref{eqn:DK_4} by
\begin{align}
\label{eqn:DK_5}
\Big(\kappa_{\pi_1\times\widetilde{\pi}_1}H(1,\pi_1\times\widetilde{\pi}_1)x\frac{\widehat{\phi}(1/T)}{T}\prod_{p\mid N_{\pi_1}}L^{RS}(1,&(\pi_{1})_{p}\times(\widetilde{\pi}_1)_{p})^{-1}\Big|\sum_{d_1,d_2\in\mathcal{D}_z}\rho_{d_1}\rho_{d_2}g([d_1,d_2])\Big|\notag\\
&+O_{\epsilon,m}(T^{\frac{3}{4}}Q^{\frac{3}{2m}}x^{1-\frac{1}{m^2}}z^{2+\epsilon}\#\mathcal{F}_{m}(Q))\Big)\sum_{\pi\in\mathcal{F}_m(Q)}|b_{\pi}|^2.
\end{align}

By proceeding as in the formulation of the Selberg sieve in \cite[Theorem 7.1]{FI}, we find that there exists a choice of $\rho_d$ satisfying conditions (1)-(3) such that
\begin{align*}
\sum_{d_1,d_2\in\mathcal{D}_z}\rho_{d_1}\rho_{d_2}g([d_1,d_2])=\sum_{\substack{d\leq z^2 \\ d\mid P(z)}}\prod_{p\mid d}\frac{g(p)}{1-g(p)}\leq\Big(\sum_{\substack{n\leq z \\ (n,N_{\pi_1})=1}}\frac{|a_{\pi_1}(n)|^2}{n}\Big)^{-1}.
\end{align*}
Note that $\kappa_{\pi\times\widetilde{\pi}}$, $H(1,\pi_1\times\widetilde{\pi}_1)$, and $\prod_{p\mid N_{\pi_1}}L^{RS}(1,(\pi_1)_p\times(\widetilde{\pi}_1)_p)^{-1}$ are each positive.  Therefore, by \cref{lem:smooth_positive_sums} and the upper bound $\widehat{\phi}(1/T)\ll 1$ from \eqref{eqn:Laplace}, we have that if $z\gg_m Q^{6m}$ with a sufficiently large implied constant, then the main term in \eqref{eqn:DK_5} is bounded by
\begin{align*}
x\cdot\frac{\kappa_{\pi\times\widetilde{\pi}}H(1,\pi_1\times\widetilde{\pi}_1)\prod_{p\mid N_{\pi_1}}L^{RS}(1,(\pi_1)_p\times(\widetilde{\pi}_1)_p)^{-1}}{\cfrac{\log z}{20}\cdot\kappa_{\pi\times\widetilde{\pi}}H(1,\pi_1\times\widetilde{\pi}_1)\prod_{p\mid N_{\pi_1}}L^{RS}(1,(\pi_1)_p\times(\widetilde{\pi}_1)_p)^{-1}+\cfrac{1}{2}}\frac{\widehat{\phi}(1/T)}{T}\ll\frac{x}{T\log z}.
\end{align*}
This establishes the bound \eqref{eqn:dual_sum}, thus concluding the proof of \cref{thm:large_sieve}.

\subsection{Mean values of Dirichlet polynomials}

Using \cref{thm:large_sieve}, we bound the mean value of the Dirichlet polynomial appearing as the integrand in  \cref{prop:zero_detection}.

\begin{prop}
	\label{prop:MVT}
	Suppose that each $\mathcal{F}_m(Q)$ satisfies \cref{hyp}, and let $\pi_0\in\mathcal{A}(m_0)$.  Let $Q\geq 3$, $T\geq 1$, and $y\geq c_m (C(\pi_0)QT\#\mathcal{F}_m(Q))^{32(m_0 m)^3}$, where $c_m>0$ is a sufficiently large constant depending at most on $m$.  For any $u\in[y,y^{12000}]$,
	\begin{align*}
	\sum_{\pi\in\mathcal{F}_m(Q)}\int_{-T}^{T}\Big|\sum_{y<p\leq u}\frac{\lambda_{\pi\times\pi_0}(p)\log p}{p^{1+it}}(1+\chi(p)p^{\beta_{\chi}-1})\Big|^2 dt\ll_m \sum_{y<p\leq u}\frac{\lambda_{\pi_0\times\widetilde{\pi}_0}(p)|1+\chi(p)p^{\beta_{\chi}-1}|^2 \log p}{p}.
	\end{align*}
\end{prop}
\begin{proof}
A result of Gallagher \cite[Theorem 1]{Gallagher} states that for any sequence of complex numbers $a_n$ and any $T\geq1$, we have
\[
\int_{-T}^{T}\Big|\sum_{n\geq1}a_n n^{-it}\Big|^2 dt\ll T^2\int_0^{\infty}\Big|\sum_{x<n\leq xe^{1/T}}a_n\Big|^2\frac{dx}{x}.
\]
Assume $z \geq c_m Q^{6m}$ with $c_m$ sufficiently large. If $b(n)$ is as in \cref{thm:large_sieve}, then the above result with $a_n=b(n)a_{\pi}(n)$ yields the bound
\[
\sum_{\pi\in\mathcal{F}_m(Q)}\int_{-T}^{T}\Big|\sum_{\substack{n\geq1 \\ P^{-}(n)>z}}b(n)a_{\pi}(n) n^{-it}\Big|^2 dt\ll T^2\int_0^{\infty}\sum_{\pi\in\mathcal{F}_m(Q)}\Big|\sum_{\substack{x<n\leq xe^{1/T} \\ P^{-}(n)>z}}b(n)a_{\pi}(n)\Big|^2\frac{dx}{x}.
\]
We apply \cref{thm:large_sieve} and bound the right hand side of the above display by
\begin{align*}
&\ll_{\epsilon,m} T^2\int_0^{\infty}\Big(\frac{x}{T\log z}+Q^{\frac{3}{2m}}T^{\frac{3}{4}}x^{1-\frac{1}{m^2}}z^{2+\epsilon}\#\mathcal{F}_m(Q)\Big)\sum_{\substack{x<n\leq xe^{1/T}\\ P^{-}(n)>z }}|b(n)|^2\frac{dx}{x}\\
&\ll_{\epsilon,m}\sum_{\substack{n\geq1 \\ P^{-}(n)>z }}|b(n)|^2 n\Big(\frac{1}{\log z}+n^{-\frac{1}{m^2}}Q^{\frac{2}{3m}}T^{\frac{7}{4}}z^{2+\epsilon}\#\mathcal{F}_m(Q)\Big)\\
&\ll_{m}\frac{1}{\log z}\sum_{\substack{n\geq1 \\ P^{-}(n)>z }}|b(n)|^2 n(1+n^{-\frac{1}{m^2}}Q^{\frac{2}{3m}}T^{\frac{7}{4}}z^{3}\#\mathcal{F}_m(Q)).
\end{align*}
Choose $y$ such that $y\geq c_m(C(\pi_0)QT\#\mathcal{F}_m(Q))^{32(m_0 m)^3}$ and $y>z$, and choose $b(n)$ to be supported on the primes $p>y$.  Then the above display is
\[
\ll_{m}(1+Q^{\frac{3}{2m}}T^{\frac{7}{4}}y^{-\frac{1}{m^2}}z^{3}\#\mathcal{F}_m(Q))\frac{1}{\log z}\sum_{p>y}|b(p)|^2 p.
\]
Choose $z = y^{1/(5m^2)}$ so, by our assumption on $y$, we indeed have $z \geq c_m Q^{6m}$. It follows that
\begin{equation}
\label{eqn:MV_T}
\sum_{\pi\in\mathcal{F}_m(Q)}\int_{-T}^{T}\Big|\sum_{p>y}b(p)a_{\pi}(p)p^{-it}\Big|^2 dt\ll_m\frac{1}{\log y}\sum_{p>y}|b(p)|^2 p.
\end{equation}
Now, select
\[
b(p)=\begin{cases}
	a_{\pi_0}(p)(1+\chi(p)p^{\beta_{\chi}-1})\frac{\log p}{p}&\mbox{if $y< p\leq u$,}\\
	0&\mbox{otherwise.}
\end{cases}
\]
Since $y>C(\pi)C(\pi_0)$ for any $\pi\in\mathcal{F}_m(Q)$, we have by \eqref{eqn:separate_dirichlet_coeffs} that
\[
a_{\pi}(p)a_{\pi_0}(p)=\lambda_{\pi\times\pi_0}(p),\qquad |a_{\pi_0}(p)|^2=\lambda_{\pi_0\times\widetilde{\pi}_0}(p)
\]
for every $p>y$.  Therefore, we may conclude from \eqref{eqn:MV_T} that if $u\in[y,y^{12000}]$, then
\begin{align*}
&\sum_{\pi\in\mathcal{F}_m(Q)}\int_{-T}^{T}\Big|\sum_{y<p\leq u}\frac{\lambda_{\pi\times\pi_0}(p)\log p}{p^{1+it}}(1+\chi(p)p^{\beta_{\chi}-1})\Big|^2 dt\\
&\ll_m\frac{1}{\log y}\sum_{y<p\leq u}\frac{\lambda_{\pi_0\times\widetilde{\pi}_0}(p) |1+\chi(p)p^{\beta_{\chi}-1}|^2 (\log p)^2}{p}\\
&\ll_m \sum_{y<p\leq u}\frac{\lambda_{\pi_0\times\widetilde{\pi}_0}(p) |1+\chi(p)p^{\beta_{\chi}-1}|^2 \log p}{p},
\end{align*}
as desired.
\end{proof}


\section{Proof of \cref{thm:lfzde_RS,thm:Landau_Siegel} and the rarity of Landau--Siegel zeros}

We now begin the proofs of \cref{thm:lfzde_RS,thm:Landau_Siegel}, both of which use \cref{prop:zero_detection,prop:MVT}.  \cref{thm:page} will follow as a straightforward consequence of \cref{thm:lfzde_RS,thm:Landau_Siegel}.  The proofs of \cref{thm:lfzde_RS,thm:Landau_Siegel} run parallel for the most part and deviate only at the very end.  By \cite[Theorem 5.8]{IK}, we have that
\[
\sum_{\pi\in\mathcal{F}_m(Q)}N_{\pi\times\pi_0}(0,T)\ll_{m,m_0} \#\mathcal{F}_m(Q) T\log(C(\pi_0)Q T).
\]
Note that the left hand sides of \cref{thm:lfzde_RS,thm:Landau_Siegel} are decreasing functions of $\sigma$, and the right hand sides are $O_{m,m_0}(1)$ when $1-\sigma\leq 1/\log(C(\pi_0)QT)$ (this uses the polynomial bound in  \cref{thm:universal}).  Thus it suffices to prove \cref{thm:lfzde_RS} when $\sigma=1-\eta/2$, with $\eta$ as in \eqref{eqn:eta_def}.

Let $\eta$ satisfy \eqref{eqn:eta_def}, and let $\tau\in\R$ satisfy $|\tau|\leq T$.  In order to simultaneously satisfy \cref{prop:zero_detection,prop:MVT}, we choose
\begin{equation}
\label{eqn:K_def2}
K =  9600(m_0 m)^3\eta\log(C(\pi_0)qQT\#\mathcal{F}_m(Q)) + O_{m_0,m}(1),
\end{equation}
where the implied constant is sufficiently large.  \cref{lem:Lemma_3.4_RJLOT} implies that there are $\ll (m_0 m)^2\log(C(\pi_0)QT)$ zeros of $L(s,\pi\times\pi_0)$ satisfying $|\rho-(1+i\tau)|\leq\eta$.  Thus
\begin{align*}
	&\frac{\#\{\rho=\beta+i\gamma\colon \beta\geq 1-\eta/2,~|\gamma-\tau|\leq\eta/2\}}{(m_0 m)^2\log(C(\pi_0)QT)}\\
	&\ll (200)^{4K}\Big[\eta^3\int_{A_1}^{A_2}\Big|\sum_{A_1<p\leq u}\frac{\lambda_{\pi\times\pi_0}(p)\log p}{p^{1+i\tau}}(1+\chi(p)p^{\beta_{\chi}-1})\Big|^2\frac{du}{u}+r_{\pi\times\pi_0}\mathbf{1}(\tau)\min\Big\{1,\frac{1-\beta_{\chi}}{\eta}\Big\}\Big].
\end{align*}
We integrate both sides over $|\tau|\leq T$ and use the bound $(m_0 m)^2\log(C(\pi_0)QT)\ll K/\eta$ to conclude that $N_{\pi\times\pi_0}(1-\eta/2,T)$ is
\begin{align*}
	\ll (201)^{4K}\Big[\eta^2 \int_{-T}^{T}\int_{A_1}^{A_2}\Big|\sum_{A_1\leq p\leq u}\frac{\lambda_{\pi\times\pi_0}(p)\log p}{p^{1+i\tau}}(1+\chi(p)p^{\beta_{\chi}-1})\Big|^2\frac{dud\tau}{u}+r_{\pi\times\pi_0}\min\Big\{1,\frac{1-\beta_{\chi}}{\eta}\Big\}\Big].
\end{align*}
The last line of the above display is
\begin{align*}
	&\ll (201)^{4K}\Big[\eta^2 \log\frac{A_2}{A_1}\max_{u\in[A_1,A_2]}\int_{-T}^{T}\Big|\sum_{A_1\leq p\leq u}\frac{\lambda_{\pi\times\pi_0}(p)\log p}{p^{1+i\tau}}(1+\chi(p)p^{\beta_{\chi}-1})\Big|^2 d\tau+r_{\pi\times\pi_0}\min\Big\{1,\frac{1-\beta_{\chi}}{\eta}\Big\}\Big]\\
	&\ll (202)^{4K}\Big[\eta\max_{u\in[A_1,A_2]}\int_{-T}^{T}\Big|\sum_{A_1\leq p\leq u}\frac{\lambda_{\pi\times\pi_0}(p)\log p}{p^{1+i\tau}}(1+\chi(p)p^{\beta_{\chi}-1})\Big|^2 d\tau+r_{\pi\times\pi_0}\min\Big\{1,\frac{1-\beta_{\chi}}{\eta}\Big\}\Big].
\end{align*}
Since $r_{\pi\times\pi_0}=1$ for at most one $\pi\in\mathcal{F}_m(Q)$ and $r_{\pi\times\pi_0}=0$ otherwise, we can sum the above display over all $\pi\in\mathcal{F}_m(Q)$ and apply \cref{prop:MVT} to obtain
\begin{align}
\label{eqn:star_point}
&\sum_{\pi\in\mathcal{F}_m(Q)}N_{\pi\times\pi_0}(1-\eta/2,T)\notag\\
&\ll_m (202)^{4K}\Big[\eta \sum_{A_1<p\leq A_2}\frac{\lambda_{\pi_0\times\widetilde{\pi}_0}(p)|1+\chi(p)p^{\beta_{\chi}-1}|^2 \log p}{p}+\min\Big\{1,\frac{1-\beta_{\chi}}{\eta}\Big\}\Big].
\end{align}

\begin{proof}[Proof of \cref{thm:lfzde_RS}]
Since $|1+\chi(p)p^{\beta_{\chi}-1}|\leq 2$, it follows from \cref{lem:Lemma_3.5_RJLOT} (with $\pi=\pi_0$ and $\eta=1/\log A_2$) that
\begin{align*}
\eta\sum_{A_1<p\leq A_2}\frac{\lambda_{\pi_0\times\widetilde{\pi}_0}(p)|1+\chi(p)p^{\beta_{\chi}-1}|^2 \log p}{p}+\min\Big\{1,\frac{1-\beta_{\chi}}{\eta}\Big\}\ll_{m_0}\eta(\log A_2 + \log C(\pi_0))\ll_{m_0} K.
\end{align*}
It follows that \eqref{eqn:star_point} is $\ll (203)^{4K}$.  Unraveling our choice of $K$ in \eqref{eqn:K_def2}, recalling our hypothesis that $q\leq 2Q$, and recalling that $\eta=2(1-\sigma)$, we find that \eqref{eqn:star_point} is
\begin{equation}
\label{eqn:bound_without_family}
\ll_{m_0,m}(C(\pi_0)qQT\#\mathcal{F}_m(Q))^{2.5\cdot 10^5(m_0 m)^3\eta}\ll_{m_0,m}(C(\pi_0)QT\#\mathcal{F}_m(Q))^{10^6(m_0 m)^3(1-\sigma)}.
\end{equation}
\cref{thm:lfzde_RS} now follows from \cref{thm:universal}.
\end{proof}

\begin{proof}[Proof of \cref{thm:Landau_Siegel}]
It follows from Weiss's arguments during his proof of \cite[Theorem 4.3]{Weiss} that
\[
\eta\sum_{A_1<p\leq A_2}\frac{|1+\chi(p)p^{\beta_{\chi}-1}|^2\log p}{p}\ll K\min\Big\{1,\frac{K}{\eta}(1-\beta_{\chi})\Big\}\ll K^2\min\Big\{1,\frac{1-\beta_{\chi}}{\eta}\Big\}.
\]
(see also \cite[Section 6]{Bombieri2}). Therefore, if $\pi_0$ satisfies GRC, then $|\lambda_{\pi_0\times\widetilde{\pi}_0}(p)|\leq m_0^2$ and
\[
\eta\sum_{A_1<p\leq A_2}\frac{\lambda_{\pi_0\times\widetilde{\pi}_0}(p)|1+\chi(p)p^{\beta_{\chi}-1}|^2 \log p}{p}+\min\Big\{1,\frac{1-\beta_{\chi}}{\eta}\Big\}\ll_{m_0}K^2\min\Big\{1,\frac{1-\beta_{\chi}}{\eta}\Big\}
\]
Incorporating \eqref{eqn:eta_def}, we insert this bound into \eqref{eqn:star_point} to complete the proof.
\end{proof}

\begin{proof}[Proof of \cref{thm:page}]
	Parts 1 and 2 follow from \cref{thm:lfzde_RS,thm:Landau_Siegel}, respectively, by choosing $\sigma=1-\frac{A}{10^7 (m_0 m)^4\log(C(\pi_0)QT)}$.
\end{proof}


\section{Subconvexity and mass equdistribution}

\noindent
\begin{proof}[Proof of \cref{thm:subconvexity}]
Recall the notation and setup of \cref{sec:Subconvexity}, especially the definition of $\mathscr{G}(Q)$ in \eqref{eqn:maass_family}.  To each $f\in\mathscr{G}(Q)$, there corresponds a cuspidal automorphic representation $\pi_f\in\mathcal{A}(2)$ with trivial central character.  Let $\mathcal{F}$ denote the set of all such $\pi_f$, and define $\mathcal{F}_2(Q)$ according to \eqref{eqn:family}.   Since $L(s,f)=L(s,\pi_f)$, it suffices for us to work with $\mathcal{F}_2(Q)$ instead of $\mathscr{G}(Q)$.  We denote by $\pi_0\in\mathcal{A}(2)$ the representation corresponding to $f_0$.

For $\pi\in\mathcal{A}(2)$, let $\mathrm{Ad}^2\pi$ denote the adjoint square lift of $\pi$; if $\pi\in\mathcal{F}_2(Q)$, then $\mathrm{Ad}^2\pi\in\mathcal{A}(3)$ and $C(\mathrm{Ad}^2\pi)\asymp \lambda_f q_f^2\leq Q^2$.  If $\pi\in\mathcal{F}_2(Q)$ and $\pi_0\in\mathcal{A}(2)$, then it follows from the uniform bound $|\alpha_{j,\pi}(p)|,|\alpha_{j,\pi_0}(p)|\leq p^{7/64}$ that both $L(3/2,\mathrm{Ad}^2\pi)$ and $L(3/2,\mathrm{Ad}^2\pi\times\pi_0)$ are defined by absolutely convergent sums which are bounded independently of $\pi$ and $\pi_0$.  (The bound $|\alpha_{j,\pi}(p)|\leq p^{7/64}$ was proved by Kim and Sarnak \cite[Appendix]{Kim} when $p$ is unramified; the ramified case was handled by Blomer and Brumley \cite{BB2}.)  Theorem 1.1 of \cite{ST} now implies that for any $0\leq\delta<1/2$, we have the bounds
\begin{equation}
\label{eqn:subconvex_2}
\log|L(1/2,\mathrm{Ad}^2\pi\times\pi_0)|\leq\Big(\frac{1}{4}-\frac{\delta}{10^9}\Big)\log C(\mathrm{Ad}^2\pi\times\pi_0)+\frac{\delta}{10^7}N_{\mathrm{Ad}^2\pi\times\pi_0}(1-\delta,6)+O(1)
\end{equation}
and
\begin{equation}
\label{eqn:subconvex_1}
\log|L(1/2,\mathrm{Ad}^2\pi)|\leq\Big(\frac{1}{4}-\frac{\delta}{10^9}\Big)\log C(\mathrm{Ad}^2\pi)+\frac{\delta}{10^7} N_{\mathrm{Ad}^2\pi}(1-\delta,6)+O(1)
\end{equation}

The $L$-function associated to the the isobaric representation $\mathrm{Ad}^2\pi\boxplus \mathrm{Ad}^2\pi\otimes\pi_0$ equals $L(s,\mathrm{Ad}^2 \pi\otimes \pi_0)L(s,\mathrm{Ad}^2 \pi)$.  Thus
\begin{equation*}
\sum_{\pi\in\mathcal{F}_2(Q)}N_{\mathrm{Ad}^2\pi\boxplus \mathrm{Ad}^2\pi\otimes\pi_0}(\sigma,T) = \sum_{\pi\in\mathcal{F}_2(Q)}N_{\mathrm{Ad}^2 \pi \otimes\pi_0}(\sigma,T) + \sum_{\pi\in\mathcal{F}_2(Q)}N_{\mathrm{Ad}^2 \pi}(\sigma,T).
\end{equation*}
By the definition of $\mathscr{G}(Q)$, each $\pi_f\in\mathcal{F}_2(Q)$ has squarefree conductor and trivial central character; therefore, it follows from the multiplicity one theorem for $\mathrm{SL}(2)$ proved by Ramakrishnan \cite[Theorem 4.2 and Corollary 4.3]{Ramakrishnan} that if $\pi,\pi'\in\mathcal{F}_2(Q)$ and $\mathrm{Ad}^2\pi=\mathrm{Ad}^2\pi'$, then $\pi'=\pi$.  Therefore, if we let $\mathcal{G}_3(Q^2)$ be the image of $\mathcal{F}_2(Q)$ in $\mathcal{A}(3)$ under the adjoint square lift, then the map $\mathrm{Ad}^2\colon \mathcal{F}_2(Q)\to\mathcal{G}_3(Q^2)$ is bijective.  Therefore,
\begin{equation*}
\sum_{\pi\in\mathcal{F}_2(Q)}N_{\mathrm{Ad}^2 \pi \otimes\pi_0}(\sigma,T) + \sum_{\pi\in\mathcal{F}_2(Q)}N_{\mathrm{Ad}^2 \pi}(\sigma,T) =  \sum_{\pi\in\mathcal{G}_3(Q^2)}N_{\pi \otimes\pi_0}(\sigma,T) + \sum_{\pi\in\mathcal{G}_3(Q^2)}N_{\pi}(\sigma,T).
\end{equation*}

By \eqref{eqn:size_maass_family}, we have that $\#\mathcal{F}_2(Q)\ll Q^2$. By the above discussion, this also means that $\#\mathcal{G}_3(Q^2)\ll Q^2$.  Thus \cref{thm:lfzde_RS} implies that for every $\epsilon>0$,
\[
\sum_{\pi\in\mathcal{G}_3(Q^2)}N_{\pi \otimes\pi_0}\Big(1-\frac{\epsilon}{10^{11}},6\Big) + \sum_{\pi\in\mathcal{G}_3(Q^2)}N_{\pi}\Big(1-\frac{\epsilon}{10^{11}},6\Big)\ll_{\pi_0} Q^{\epsilon}.
\]
Therefore, the number of $\pi\in\mathcal{F}_2(Q)$ such that $L(s,\mathrm{Ad}^2\pi\otimes\pi_0)L(s,\mathrm{Ad}^2\pi)$ has a zero in the region $\re(s)\geq 1-\frac{\epsilon}{10^{11}}$ and $|\im(s)|\leq 6$ is $O_{\pi_0}(Q^{\epsilon})$.  For each of the remaining $\pi\in\mathcal{F}_2(Q)$, it follows that both $N_{\mathrm{Ad}^2\pi\times\pi_0}(1-\frac{\epsilon}{10^{11}},6)$ and $N_{\mathrm{Ad}^2\pi}(1-\frac{\epsilon}{10^{11}},6)$ equal zero.  By \eqref{eqn:BH}, equations \eqref{eqn:subconvex_2} and \eqref{eqn:subconvex_1} now read
\begin{align*}
\log|L(1/2,\mathrm{Ad}^2\pi\times\pi_0)|&\leq\Big(\frac{1}{4}-10^{-20}\epsilon\Big)\log C(\mathrm{Ad}^2\pi\times\pi_0)+O(1)
\end{align*}
and
\[
\log|L(1/2,\mathrm{Ad}^2\pi)|\leq\Big(\frac{1}{4}-10^{-20}\epsilon\Big)\log C(\mathrm{Ad}^2\pi)+O(1),
\]
respectively.  Twisting $\mathrm{Ad}^2\pi$ by $|\det|^{-it}$, we see from \eqref{eqn:BH} that
\begin{align*}
\log|L(1/2+it,\mathrm{Ad}^2\pi)|&\leq\Big(\frac{1}{4}-10^{-20}\epsilon\Big)\log C(\mathrm{Ad}^2\pi,t)+O(1).
\end{align*}
for any $t\in\mathbb{R}$.  Observe the factorizations $L(1/2+it,f\times f)=\zeta(1/2+it)L(1/2+it,\mathrm{Ad}^2 \pi_f)$ and $L(1/2,f\times f\times f_0)=L(1/2,\pi_{0})L(1/2,\mathrm{Ad}^2 \pi_f \times \pi_{0})$ (where $\zeta(s)$ denotes the Riemann zeta function).  Since $C(\mathrm{Ad}^2\pi_f,t)\ll\lambda_f q_f^2(2+|t|)^{3}$ and $C(\mathrm{Ad}^2\pi_f\times\pi_0)\ll \lambda_f^2 q_f^4 \lambda_{f_0}^3 q_{f_0}^3$ via \eqref{eqn:BH}, the proof of \cref{thm:subconvexity} is complete.
\end{proof}

\section{The Chebotarev density theorem in families}
The goal of this section is to prove \cref{thm:CDT}. Let $L/\Q$ be a Galois extension of number fields.  We begin by establishing a flexible variant of the Chebotarev density theorem. Given any zero-free region for the Dedekind zeta function $\zeta_L(s)$, we would like to compute an asymptotic expression for $\pi_C(x,L/\Q)$ with an error term depending on the zero-free region in an explicit form.  

\begin{prop} \label{prop:FlexiError} Let $L/\Q$ be a Galois extension of number fields with Galois group $G$. Let $\Delta : [3,\infty) \rightarrow (0,\infty)$ be a function such that $\zeta_L(s)/\zeta_{\Q}(s) \neq 0$ in the region
\[
\Re(s) > 1- \Delta(|\Im(s)| +3). 
\]
Define 
	\begin{equation}
	\eta(x) =  \inf_{t \geq 3} \Big[ \Delta(t) \log x + \log t \Big].
	\label{def:eta_function}
	\end{equation}
	Let $C$ be a conjugacy class of $G$, and suppose there exists an abelian subgroup $H$ of $G$ such that $H \cap C$ is non-empty and $\zeta_{L^H}(s)/\zeta_{\Q}(s)$ is entire, where $L^H$ is the subfield of $L$ fixed by $H$.  For $x \geq (\log D_L)^4$,
	\[
	\Big| \pi_C(x,L/\Q) - \frac{|C|}{|G|} \pi(x) \Big| \ll  \frac{|C|}{|G|} \frac{x e^{-\frac{1}{8}\eta(x)} }{\log x} \log D_L  + \frac{|C|}{|G|}\frac{x^{3/4}}{\log x}. 
	\]
\end{prop}
\begin{remark1}
\begin{enumerate}
\item The existence of this abelian subgroup $H$ is a mild condition for our purposes.  In the special case $C = \{1\}$, one can take $H = \{1\}$ and this follows unconditionally from the Aramata--Brauer theorem as $L^H = L$ is Galois over $\Q$. For an arbitrary conjugacy class $C$, one can take $H = \langle g \rangle$ to be the cyclic subgroup generated by some element $g \in C$ in which case this assumption follows easily from the strong Artin conjecture for $\zeta_L(s)$ over $\Q$. The strong Artin conjecture is known for all examples under consideration in \cref{thm:CDT}. 
\item An analogous result holds for any Galois extension $L/F$ with $\pi(x)$ replaced by the number of prime ideals of $F$ up to $x$ and $\zeta_{\Q}(s)$ replaced by $\zeta_F(s)$. We restrict to $F=\Q$ for simplicity and with \cref{thm:CDT} in mind.
\end{enumerate}
\end{remark1}

\begin{proof} For the proof, we will borrow heavily from results recorded in \cite{TZ3} and will therefore remain consistent with the notation therein.  Let $g \in H \cap C$ be arbitrary and set $C_H = \{g\}$. Let $K = L^H$ be the fixed field of $L$ by $H$. Select $f( \, \cdot \, ) = f(\, \cdot \, ; x, \ell, \epsilon)$ in \cite[Lemma 2.2]{TZ3} with
	\begin{equation}
	\epsilon = \min\{\tfrac{1}{8}, 8e^{-\eta(x)/4}\} + x^{-1/4}, \qquad \ell = 2. 
	\label{eqn:AlternateZFR_ParameterChoice}
	\end{equation}
	Note that $0 \leq f(t) \leq 1$ for all $t \in \R$, $f(t)$ supported in $t \in [\frac{1}{2}-\epsilon, 1+\epsilon]$ and $f(t) \equiv 1$ for $t \in [\frac{1}{2},1]$.  Its Laplace transform $F(z) = \int_0^{\infty} f(t) e^{-zt} dt$ is entire and satisfies many properties recorded in \cite[Lemma 2.2]{TZ3}. Consider the weighted prime sum $\widetilde{\psi}_{C_H}(x,f) = \widetilde{\psi}_{C_H}(x, L/L^{H}; f)$ given by \cite[Equation 2.9]{TZ3} or, equivalently, 
	\[
	\widetilde{\psi}_{C_H}(x;f) = \frac{|C_H|}{|H|} \sum_{\chi \in \widehat{H}} \bar{\chi}(C_H)  \frac{\log x}{2\pi i} \int_{2-i\infty}^{2+i\infty} - \frac{L'}{L}(s,\chi,L/L^{H}) F(-s\log x) ds, 
	\]
	where $\chi$ runs over all the (Hecke) characters of the dual group $\widehat{H}$. 
	By \cite[Lemma 4.3]{TZ3}, the bound $\epsilon \geq x^{-1/4}$ from \eqref{eqn:AlternateZFR_ParameterChoice}, and the bounds $n_L \ll \log D_L \leq x^{1/4}$, it follows that
	\begin{equation}
	\frac{|H|}{|C_H|} \frac{\widetilde{\psi}_{C_H}(x;f)}{\log x} = F(-\log x) -  \sum_{\chi \in \widehat{H}} \bar{\chi}(C_H) \sum_{\rho_{\chi}}  F(-\rho_{\chi} \log x)  + O\Big( \frac{x^{1/2}}{\log x} \Big),  
	\label{eqn:AlternateZFR_smooth}
	\end{equation}
	where $\rho_{\chi}$ runs over all  non-trivial zeros of the Hecke $L$-functions $L(s,\chi,L/L^{H})$. Note
	\[
	\zeta_L(s) = \zeta_K(s) \prod_{\substack{\chi \in \widehat{H} \\ \chi \neq 1}} L(s,\chi,L/L^{H})
	\] 
	and, by assumption, $\zeta_{L^H}(s)/\zeta_{\Q}(s)$ is entire. The   zeros of $\zeta_{\Q}(s)$ therefore contribute only to the zeros of the trivial character $\chi = 1$ in \eqref{eqn:AlternateZFR_smooth}. From these observations, it follows that 
	\begin{equation}
	F(-\log x) -  \sum_{\chi \in \widehat{H}} \bar{\chi}(C_H) \sum_{\rho_{\chi}}  F(-\rho_{\chi} \log x) = S(x) + O\Big( \sum_{\substack{ \rho \\ \frac{\zeta_L}{\zeta_\Q}(\rho) =0}} |F(-\rho \log x)|\Big),
	\label{eqn:AlternateZFR_DecomposeZeros}
	\end{equation}
	where 
	\[
	S(x) = F(-\log x) - \sum_{\substack{ \rho \\ \zeta_{\Q}(\rho) = 0}} F(-\rho \log x).
	\] 
	By standard arguments using Mellin inversion, one can verify that 
	\begin{equation}
	 (\log x) S(x) = \sum_{n \geq 1} \Lambda(n) f\Big( \frac{\log n}{\log x}\Big) +(\log x) F(0) - \frac{\log x}{2\pi i} \int_{-1/2-i\infty}^{-1/2+i\infty} -\frac{\zeta_{\Q}'}{\zeta_{\Q}}(s) F(-s\log x) ds. 	\label{eqn:AlternateZFR_RationalPrimes}
	\end{equation}
	By \cite[Lemma 2.2(iv)]{TZ3}, $|F(0)| \ll 1$. From the properties of $f$ described immediately following \eqref{eqn:AlternateZFR_ParameterChoice} and the prime number theorem, 
	\[
	\sum_{n \geq 1} \Lambda(n) f\Big( \frac{\log n}{\log x}\Big) =  \sum_{n \leq x} \Lambda(n)  + O( \epsilon x + x^{1/2}). 
	\]
	For $\Re(s) = -1/2$, we have
	\[
	-\frac{\zeta_{\Q}'}{\zeta_{\Q}}(s) \ll \log(|\Im(s)|+3), 
	\qquad 
	(\log x) |F(-s\log x)| \ll \epsilon^{-2} x^{-1/4} \ll x^{1/4},
	\]
	which follow from \cite[Lemmas 2.2(vi) and 2.5]{TZ3} and \eqref{eqn:AlternateZFR_ParameterChoice}. Combining all of these observations with \eqref{eqn:AlternateZFR_RationalPrimes} and noting $\epsilon \ll e^{-\eta(x)/4} + x^{-1/4}$ by \eqref{eqn:AlternateZFR_ParameterChoice}, it follows that
	\begin{equation}
	(\log x) S(x) = \sum_{n \leq x} \Lambda(n)  + O(x e^{-\eta(x)/4} + x^{3/4}). 
	\label{eqn:AlternateZFR_mainterm}
	\end{equation}
	All that remains is to consider the error term in \eqref{eqn:AlternateZFR_DecomposeZeros}. By \cite[Lemma 4.4]{TZ3} and the assumption $\log D_L \leq x^{1/4}$, the zeros $\rho$ with $|\rho| \leq 1/4$ have negligible contribution; namely, 
	\[
	\sum_{\substack{ \rho \\ \frac{\zeta_L}{\zeta_\Q}(\rho) =0}} |F(-\rho \log x)| = \sum_{\substack{ |\rho| \geq 1/4 \\ \frac{\zeta_L}{\zeta_\Q}(\rho) =0}} |F(-\rho \log x)| + O( x^{1/2} )  \qquad \text{for $x \geq 3$}. 
	\]
	Write $\rho = \beta+i\gamma$ for each non-trivial zero $\rho$. By \eqref{def:eta_function}, one can see that $\dfrac{x^{-(1-\beta)}}{(|\gamma|+3)} \leq e^{-\eta(x)}$. Thus, \cite[Lemma 2.2(iv)]{TZ3} and \eqref{eqn:AlternateZFR_ParameterChoice} imply that, for $|\rho| \geq 1/4$,
	\[
	(\log x) |F(-\rho \log x)| \ll \frac{x^{\beta}}{(|\gamma|+3)} \cdot \frac{\epsilon^{-2}}{(|\gamma|+3)^2} \ll xe^{-\eta(x)} \cdot \frac{e^{\eta(x)/2}}{(|\gamma|+3)^2} . 
	\] 
	Summing over all such  zeros, it follows that
	\begin{equation*}
	\begin{aligned}
		\sum_{\substack{ |\rho| \geq 1/4 \\ \frac{\zeta_L}{\zeta_\Q}(\rho) =0}} |F(-\rho \log x)| &
		  \ll \frac{x e^{-\eta(x)/2}}{\log x}
		 \sum_{\substack{\frac{\zeta_L}{\zeta_\Q}(\rho) = 0}} \frac{1}{(|\gamma|+3)^2}.
	\end{aligned}
	\end{equation*}
	Applying a standard estimate for the zeros of the Dedekind zeta function \cite[Lemma 2.5]{TZ3} and Minkowski's bound $n_L \ll \log D_L$, we see that the above expression is 
	\begin{equation}
	\ll \frac{xe^{-\eta(x)/2}}{\log x} \sum_{T=1}^{\infty}  \sum_{\substack{\frac{\zeta_L}{\zeta_{\Q}}(\rho) = 0 \\ T-1 \leq |\gamma| < T}} \frac{\log D_L + n_L\log(T+3)}{T^2}  \ll  \frac{xe^{-\eta(x)/2} \log D_L}{\log x}. 
	\label{eqn:AlternateZFR_zeros}
	\end{equation}
	Substituting \eqref{eqn:AlternateZFR_zeros},   \eqref{eqn:AlternateZFR_mainterm}, and \eqref{eqn:AlternateZFR_DecomposeZeros} into \eqref{eqn:AlternateZFR_smooth}, we conclude that
	\[
	\frac{|H|}{|C_H|}\widetilde{\psi}_{C_H}(x;f) =  \sum_{n \leq x} \Lambda(n) + O\Big(x e^{-\eta(x)/4} \log D_L + x^{3/4} \Big) \quad \text{for $x \geq (\log D_L)^4$. }
	\]
	 Via \cite[Lemma 2.3]{TZ3}, we may replace $\widetilde{\psi}_{C_H}(x;f)$ by the usual prime counting function $\psi_{C_H}(x)$ given by \cite[Equation 2.1]{TZ3} at the cost of $O(\epsilon x + x^{1/2})$. From \eqref{eqn:AlternateZFR_ParameterChoice}, this cost is absorbed into the existing error term in the above expression. By partial summation (see \cite[Lemma 2.1 and Equation 5.3]{TZ3}), it therefore follows that
	\begin{equation*}
	\begin{aligned}
	\frac{|H|}{|C_H|} \pi_{C_H}(x) & = \pi(x)  + O\Big( (\log D_L) \frac{x}{\log x}  \sup_{\sqrt{x} \leq y \leq x}\big( e^{-\eta(y)/4}\big) +	\frac{x^{3/4}}{\log x} + \log D_L\Big). 
	\end{aligned}
	\end{equation*}
	By \eqref{def:eta_function}, one can verify that $\eta(y)$ is an increasing function of $y$ and also $\eta(x^{1/2}) \geq \frac{1}{2} \eta(x)$. With these observations and the assumption  $\log D_L \leq x^{1/4}$, we conclude that 
	\begin{equation*}
	\begin{aligned}
	  \pi_{C_H}(x) = \frac{|C_H|}{|H|} \pi(x)   + O\Big(   \frac{|C_H|}{|H|} \frac{x}{\log x}  e^{-\eta(x)/8}  \log D_L +  \frac{|C_H|}{|H|}  \frac{x^{3/4}}{\log x} \Big).
	\end{aligned}
	\end{equation*}
	\cref{prop:FlexiError} now follows by an application of \cite[Lemma 5.2]{TZ3} from class field theory. To absorb the arising secondary error term, we again use that $n_L \ll \log D_L \leq x^{1/4}$. 
\end{proof}

First, we record a classical zero-free region for the Dedekind zeta function; for a proof, see \cite[Lemma 2.3]{LMO} for instance. 

\begin{lem} \label{lem:DedekindZeta_ZFR}
	The Dedekind zeta function $\zeta_L(s)$ has at most one simple real zero in the region
	\[
	\Re(s) > 1 - \frac{\Cl[abcon]{DedekindZetaZFR}}{\log D_L + n_L \log(|\Im(s)|+3)}.
	\]
\end{lem}

Assuming a strong zero-free region for the Dedekind zeta function,  we arrive at a natural form of the Chebotarev density theorem

\begin{thm} \label{thm:QuasiGRH}
Let $L/\Q$ be a Galois extension of number fields with Galois group $G$ and $L \neq \Q$. Let $C$ be a conjugacy class of $G$ satisfying the hypotheses of \cref{prop:FlexiError}. Let $0 < \delta \leq 1/2$ and $T \geq (\log D_L)^{24}$ be arbitrary. Assume $\zeta_L(s)/\zeta_{\Q}(s)$ has no zeros in the region
\begin{equation}
	\Re(s) > 1-\delta, \qquad |\Im(s)| \leq T. 
	\label{eqn:CDT_ZFR}
\end{equation}	
For $x \geq (\log D_L)^{16/\delta}$, 
\[
\Big|\pi_C(x,L/\Q) - \frac{|C|}{|G|} \pi(x)\Big| \ll \frac{|C|}{|G|}\frac{x}{\log x}\Big( x^{-\delta/8} + T^{-\frac{1}{24}} e^{- \frac{1}{24}\sqrt{\Cr{DedekindZetaZFR}(\log x)/n_L}} + T^{-\frac{1}{24}} e^{- \frac{1}{24} \frac{\Cr{DedekindZetaZFR} \log x}{\log D_L}} \Big). 
\]
\end{thm}

\begin{proof} By \cref{prop:FlexiError} and \cref{lem:DedekindZeta_ZFR}, it remains to compute $\eta(x)$ for
\[
\Delta(t) = \begin{cases} 
		 \delta , & 3 \leq t \leq T, \\
		\Cr{DedekindZetaZFR}(\log D_L + n_L \log t)^{-1}, &  t > T.
 \end{cases}
\]
Define $\eta(x) = \min\{ \eta_1(x), \eta_2(x) \}$, where
\[
\eta_1(x) = \inf_{3 \leq t \leq T}( \delta \log x +   \log t)  \qquad\textup{and}\qquad
\eta_2(x) = \inf_{t \geq T}(\frac{\Cr{DedekindZetaZFR} \log x}{\log D_L + n_L \log t} +  \log t).
\]
If $\eta(x) = \eta_1(x)$, then $\eta(x) \geq \delta \log x$. Otherwise, we may assume $\eta(x) = \eta_2(x)$. Arguing as in \cite[Lemma 4.6]{TZ3}, the expression $\frac{\Cr{DedekindZetaZFR} \log x}{\log D_L + n_L u} + u$ is positive for $u \geq 0$ and is globally minimized in this interval at $u = \max\{0,u_0\}$ where  $u_0 = \frac{(\Cr{DedekindZetaZFR} \log x)^{1/2}}{n_L^{1/2}} - \frac{\log D_L}{n_L}$. Therefore,
\[
\eta(x) = \eta_2(x) \geq \min\Big\{ \frac{\Cr{DedekindZetaZFR} \log x}{\log D_L}, \sqrt{\frac{\Cr{DedekindZetaZFR} \log x}{n_L}} \Big\}. 
\]
Since one always has the lower bound $\eta_2(x) \geq \log T \geq 24 \log D_L$, we see in all cases that
\begin{equation*}
\begin{aligned}
e^{-\eta(x)/8} & \leq e^{-\eta_1(x)/8} + e^{-\eta_2(x)/8} \\
	& \leq x^{-\delta/8} + e^{-\eta_2(x)/24} T^{-1/24} (\log D_L)^{-1} \\
	& \leq (\log D_L)^{-1} \big( x^{-\delta/16} + T^{-\frac{1}{24}} e^{-\frac{\Cr{DedekindZetaZFR} \log x}{24 \log D_L}} + T^{-\frac{1}{24}} e^{-\frac{1}{24}\sqrt{\Cr{DedekindZetaZFR} (\log x)/n_L}}) 	
\end{aligned} 
\end{equation*}
because $x \geq (\log D_L)^{16/\delta}$. This estimate, along with \cref{prop:FlexiError,lem:DedekindZeta_ZFR}, yields the result.  
\end{proof}

We conclude this section with the proof of \cref{thm:CDT}. 
\begin{proof}[Proof of \cref{thm:CDT}]
Recall $\mathscr{F}(X) = \mathscr{F}(X; n, G, R_G)$ is a family of number fields over $\Q$ whose Galois closure has Galois group is isomorphic to $G$, where $G$ is a fixed transitive subgroup of $S_n$ equal to one of $C_n, S_3, S_4, D_p$ or $A_4$. 
Let $K \in \mathscr{F}(X)$ and recall $\widetilde{K}/\Q$ is the Galois closure of $K$ over $\Q$. For $\Re(s) > 1$, 
\begin{equation}
\zeta_{\widetilde{K}}(s) = \zeta_{\Q}(s) \prod_{\rho\neq 1} L(s,\rho,\widetilde{K}/\Q)^{\dim \rho},
\label{eqn:ArtinFactorization} 
\end{equation}
where $\rho$ runs over the \textit{non-trivial} irreducible Artin representations of $G$. In all cases under consideration, the strong Artin conjecture is known for all of the non-trivial Artin representations $\rho$ of $G$. That is, $L(s,\rho,\widetilde{K}/\Q) = L(s,\pi)$ for some cuspidal automorphic representation $\pi = \pi_{\rho}$ of $\mathrm{GL}_d(\mathbb{A}_{\Q})$ with $d$ equal to the degree of $\rho$. Observe that $d$ is bounded by $m$, where $m = m(G)$ is the maximum degree of the irreducible representations of $G$. The map 
\begin{equation}
\rho \mapsto \pi_{\rho}
\label{eqn:StrongArtinMap}
\end{equation}
has image $\mathscr{A}(X) = \mathscr{A}(X;G,n,R_G)$, the set of automorphic representations $\pi$ obtained this way from $\mathscr{F}(X)$. 

Let $M(X) = M(X;G, n, R_G)$ be the maximum size of the fibres of the map in \eqref{eqn:StrongArtinMap}. As shown in \cite{PTW}, 
\begin{equation}
M(X) = \max_{F \neq \Q} \#\{ K \in \mathscr{F}(X) : \Q \subset F \subseteq \widetilde{K} \}, 
\label{eqn:MaxFibre}
\end{equation}
where the maximum runs over all number fields $F \neq \Q$. Since our notation differs with theirs, we explain \eqref{eqn:MaxFibre} for the sake of clarity. Fix some $\pi \in \mathscr{A}(X)$.  By a result of Kl\"{u}ners and Nicolae 
\cite[Theorem 5]{KlunersNicolae} refined by Pierce--Turnage-Butterbaugh--Wood \cite[Lemma 7.4]{PTW}, it follows that\footnote{Here we crucially use that the base field is $\Q$.}  $L(s,\rho_1,\widetilde{K}_1/\Q) = L(s,\rho_2,\widetilde{K}_2/\Q) = L(s,\pi)$ if and only if
\[
\widetilde{K}_1^{\ker(\rho_1)} = \widetilde{K}_2^{\ker(\rho_2)} = F
\]
for some number field $F$. Note that $F \neq \Q$ since the representations $\rho_1, \rho_2$ are non-trivial.  Hence, the size of the fibre above $\pi \in \mathscr{A}(X)$ in  \eqref{eqn:StrongArtinMap} equals $\#\{ K \in \mathscr{F}(X) : \Q \subset F \subseteq \widetilde{K} \}$ for some number field $F \neq \Q$, implicitly depending on $\pi$. This implies \eqref{eqn:MaxFibre}.

In light of \eqref{eqn:MaxFibre}, it follows from \cite[Proposition 7.9]{PTW} and \cite[Theorem 7.1]{PTW} that there exists a sufficiently small $\epsilon = \epsilon(n,G) > 0$ such that
\begin{equation}
M(X) \ll_{n,G,\epsilon} X^{-2\epsilon} \#\mathscr{F}(X).
\label{eqn:DiscMultiplicity}
\end{equation}
This result is one of the key innovations of \cite{PTW}.

Now, we verify the assumptions of \cref{thm:lfzde_RS} with $\pi_0\in\mathcal{A}(1)$ taken to be the trivial representation. Take $m = m(G)$ to be the maximum degree of the irreducible representations of $G$, $Q = X^{|G|/2}$, and $\cF_m(Q) = \mathscr{A}(X)$. By \eqref{eqn:ArtinFactorization} and  \eqref{eqn:StrongArtinMap},  each $\pi \in \cF_m(Q)$ satisfies 
\[
\deg(\pi) \leq m \qquad \text{ and } \qquad 	 C(\pi) \leq D_{\widetilde{K}}\text{ for some $K \in \mathscr{F}(X)$}.
\] 
Since $D_{\widetilde{K}} \leq D_K^{|G|/2} \leq X^{|G|/2} = Q$ for any $K \in \mathscr{F}(X)$, we indeed have that $C(\pi) \leq Q$ for every $\pi \in \cF_m(Q)$. Moreover, $\pi \in \cF_m(Q)$ satisfies GRC (and hence \cref{hyp}) since it corresponds to an Artin representation via \eqref{eqn:StrongArtinMap}. Thus, by \cref{thm:lfzde_RS}, it follows that  
\begin{equation}
\sum_{\pi \in \mathscr{A}(X)} N_{\pi}(1-\delta,T) \ll_{n,G} ( X^{|G|/2} T)^{10^7 m^4 \delta } 
\label{eqn:CDT_LFZDE} 
\end{equation}
uniformly for $T \geq 1$ and $0 < \delta < 1/2$. For $\epsilon \in (0,1)$ arbitrary, select 
\[
T = Q (\log Q)^{24}, \qquad  \delta = \frac{\epsilon}{10^8 |G| m^4} 
\]
Thus, by \eqref{eqn:CDT_LFZDE} and our definition of $Q = X^{|G|/2}$, for all except at most $O_{n,G,\epsilon}(X^{\epsilon})$ automorphic representations $\pi \in \mathscr{A}(X)$, the $L$-function $L(s,\pi)$ is zero-free in the region 
\begin{equation}
\Re(s) > 1-\delta, \qquad |\Im(s)| \leq Q(\log Q)^{24}. 
\label{eqn:CDT_strongZFR}
\end{equation}
Each exceptional $\pi$ corresponds to at most $M(X)$ exceptional fields $K \in \mathscr{F}(X)$.  Throwing out all of these exceptional fields, it follows by \eqref{eqn:DiscMultiplicity} that $\zeta_{\widetilde{K}}(s)/\zeta_{\Q}(s)$ is zero-free in the region \eqref{eqn:CDT_strongZFR} for all $K \in \mathscr{F}(X)$ with at most $O_{n,G,\epsilon}(X^{-\epsilon} \#\mathscr{F}(X) )$  exceptions.

Now, let $K \in \mathscr{F}(X)$ be a non-exceptional field. By \cref{thm:QuasiGRH}, we have that 
\begin{equation}
\label{eqn:asif_hates_pauses}
\Big| \pi_C(x,\widetilde{K}/\Q) - \frac{|C|}{|G|} \pi(x) \Big| \ll \frac{|C|}{|G|}\frac{x}{\log x} E(x) \qquad \text{ for $x \geq (\log D_{\widetilde{K}})^{16/\delta}$,}
\end{equation}
 where
\[
E(x) = x^{-\delta/8} + D_{\widetilde{K}}^{-\frac{1}{24}} \exp\Big[ - \frac{1}{24} \Big( \frac{\Cr{DedekindZetaZFR} \log x}{|G|} \Big)^{1/2}  \Big] + D_{\widetilde{K}}^{-\frac{1}{24}} \exp\Big[- \frac{1}{24} \frac{\Cr{DedekindZetaZFR} \log x}{\log D_{\widetilde{K}}}\Big].
\] 
Note we used that $D_{\widetilde{K}} \leq Q$ to express $E(x)$ in terms of $D_{\widetilde{K}}$ instead of $Q$. Choose $\eta = \delta/8$.  For $(\log D_{\widetilde{K}})^{2/\eta} \le x \leq D_{\widetilde{K}}^{(24\eta)^{-1}}$, one can directly verify that $E(x) \ll x^{-\delta/8} = x^{-\eta}$.  If $(24\eta)^{-1} \log D_{\widetilde{K}} \leq \log x \leq \Cr{DedekindZetaZFR}^{-1} |G| (\log D_{\widetilde{K}})^2$  then one can verify that $E(x) \ll D_{\widetilde{K}}^{-1/24} \ll e^{-\frac{1}{24} \sqrt{ \Cr{DedekindZetaZFR} (\log x)/|G|}}$.  Finally, if $\log x \geq \Cr{DedekindZetaZFR}^{-1} |G| (\log D_{\widetilde{K}})^2$ then one can verify that 
\[
E(x) \ll e^{-\frac{1}{24} \sqrt{ \Cr{DedekindZetaZFR} (\log x)/|G|}} + e^{-\frac{ \Cr{DedekindZetaZFR} \log x }{24 \log D_{\widetilde{K}}} } \ll e^{-\frac{1}{24} \sqrt{ \Cr{DedekindZetaZFR} (\log x)/|G|}}. 
\]
This completes the proof of \cref{thm:CDT}. 
\end{proof}

\section{Landau--Siegel zeros and torsion in class groups}
This section is dedicated to the proof of \cref{thm:ell-torsion}. The first ingredient is a lemma due to Ellenberg--Venkatesh \cite[Lemma 2.3]{EV}. It establishes a connection between the existence of small split primes and bounds for the class group. 
\begin{lem}[Ellenberg--Venkatesh] \label{lem:SmallSplitPrimes}
	Let $K/\Q$ be a number field of degree $n$ and let $\ell \geq 1$ be a positive integer. Set $0 < \delta < \frac{1}{2\ell(n-1)}$ and suppose there exists $M$ rational primes $p \leq D_K^{\delta}$ which are unramified and split completely in $K$. For any $\epsilon > 0$,
	\[
	|\mathrm{Cl}_K[\ell]| \ll_{\epsilon, \ell, n} D_K^{\frac{1}{2}+\epsilon} M^{-1}.
	\] 
\end{lem}

To make use of \cref{lem:SmallSplitPrimes}, we require a proposition relating low-lying zero free regions to the existence of small primes with a given  splitting behaviour.

\begin{prop} \label{prop:ell-torsion_primes}
	Let $L/\Q$ be a Galois extension of number fields and let $0 < \epsilon < \delta/2$ be arbitrary. Suppose $\zeta_L(s)$ has no zeros in the region
	\begin{equation}
		\Re(s) > 1 - \frac{H_{\delta,\epsilon}}{\log D_L}, \qquad |\Im(s)| \leq 1, 
		\label{eqn:ZFR_low}
	\end{equation}
	where $H_{\delta,\epsilon} \geq 1$ is sufficiently large. Then, for any conjugacy class $C \subseteq G$ , 
	\[
	\pi_C( D_L^{\delta}, L/\Q) \geq 	\frac{\epsilon}{8 \delta} \frac{|C|}{|G|}  D_L^{\delta - \epsilon} + O_{\delta,\epsilon}\Big(	\frac{|C|}{|G|}  D_L^{\delta - \epsilon} (\log D_L)^{-3} \Big).  
	\] 
\end{prop}
\begin{proof}
	This essentially follows from the arguments found in \cite{Zaman1}. We will outline the proof here and borrow heavily from \cite{Zaman1}, so we will remain as consistent as possible with the notation therein. In particular, set $\sL = \log D_L$.  Select $f$ as in \cite[Lemma 2.6]{Zaman1} with $\ell = 2, B = \delta,$ and $A = \epsilon/4$. Then 
	\begin{itemize}
		\item 	$0 \leq f(t) \leq A^{-1} \leq 4 \epsilon^{-1}$ for all $t \in \R$
		\item The support of $f$ is contained in $[B-2\ell A, B] = [\delta - \epsilon, \delta]$. 
		\item The Laplace transform $F(z) = \int_0^{\infty} f(t) e^{-zt}dt$ is entire and given by 
			\[
			F(z) = e^{-(B-2\ell A) z} \Big( \frac{1-e^{-Az}}{Az} \Big)^{2\ell} = e^{-(\delta-\epsilon) z} \Big( \frac{1-e^{-\epsilon z/4}}{ \epsilon z/4} \Big)^4. 
			\]
		\item For $s = \sigma + it \in \R$ with $\sigma < 1$ amd $t \in \R$, we have:
			\[
			|F((1-s) \sL)| \ll_{\epsilon} e^{-(\delta-\epsilon) (1-\sigma) \sL} \min\{1, |(1-s)\sL|^{-4}\}. 
			\]
			Furthermore, $F(0) = 1$. 
	\end{itemize}
	We will use these properties frequently and often without mention. Define
	\[
	S = \sum_{\substack{p \text{ prime} \\ p \nmid D_L}} \frac{\log p}{p} f\Big( \frac{\log p}{\sL} \Big) \mathbf{1}_C(p), 
	\]
	where, for primes $p$ unramified in $L$,  $\mathbf{1}_C(p) = 1$ if $[\tfrac{L/\Q}{p}] = C$ and $0$ otherwise. By the properties of $f$, one can verify that
	\begin{equation}
	S \leq  \frac{\delta \sL}{e^{(\delta-\epsilon) \sL}} \cdot 4 \epsilon^{-1} \cdot \sum_{\substack{p \leq D_L^{\delta} \\ p \nmid D_L}}  \mathbf{1}_C(p) \leq (4 \delta \epsilon^{-1} D_L^{-\delta+\epsilon} \log D_L)  \cdot  \pi_C(D_L^{\delta}, L/\Q). 
	\label{eqn:ZFR_low_naturalprimesum}
	\end{equation}
	Now, from the proof of \cite[Lemma 4.1]{Zaman1}, we have that
	\[
	\sL^{-1} S =  \sum_{\psi} \bar{\psi}(C) \frac{1}{2\pi i} \int_{2-i\infty}^{2+i\infty} -\frac{L'}{L}(s,\psi,L/\Q) F((1-s)\sL) ds + O_{\delta,\epsilon}(\sL^2 e^{-\delta \sL/4} ),
	\]
	where $\psi$ runs over the irreducible Artin characters of $\mathrm{Gal}(L/\Q)$. Using standard class field theory arguments (see \cite[Section 4.2]{Zaman1}), one can shift the contour as in \cite[Lemma 4.2]{Zaman1}  with $T_{\star} = 1$. This yields
	\begin{equation}
	\frac{|G|}{|C|} \sL^{-1} S = 1 + O_{\delta,\epsilon}\Big( \sum_{\substack{ |\Im(\rho)| \leq 1}} |F((1-\rho)\sL)| + \sL^{-3} \Big),
	\label{eqn:ZFR_low_primesum}
	\end{equation}
	where $\rho$ runs over all non-trivial zeros of $\zeta_L(s)$ satisfying $|\Im(\rho)| \leq 1$.  We apply \cite[Lemma 4.3]{Zaman1}  (with $J = 1, T_1 = 1,$ and $R_1 = H_{\delta,\epsilon}$ in their notation) to deduce that
	\[
	\sum_{\substack{ |\Im(\rho)| \leq 1}} |F((1-\rho)\sL)|  = \sum_{\substack{ |\Im(\rho)| \leq 1 \\ \Re(\rho) > 1 - \frac{H_{\delta,\epsilon}}{\log D_L}}} |F((1-\rho)\sL)|  + O_{\delta,\epsilon}\big( e^{- \delta H_{\delta,\epsilon}/2 } \big). 
	\] 
	By assumption \eqref{eqn:ZFR_low}, the remaining sum over zeros is empty. Combining these estimates with \eqref{eqn:ZFR_low_primesum} implies that
	\[
	S = \frac{|C|}{|G|} \log D_L \Big(1+ O_{\delta,\epsilon}\big( e^{-\delta H_{\delta,\epsilon}/2}  + (\log D_L)^{-3} \big)  \Big) \geq \frac{1}{2} \frac{|C|}{|G|} \log D_L \Big(1 + O_{\delta,\epsilon}\big( (\log D_L)^{-3} \big) \Big), 
	\]
	since $H_{\delta,\epsilon}$ is sufficiently large. Substituting this lower bound into \eqref{eqn:ZFR_low_naturalprimesum} yields the result. 
\end{proof}

\noindent
\begin{proof}[Proof of \cref{thm:ell-torsion}] Recall $K$ is a number field of degree $n$ and $\widetilde{K}$ is its Galois closure over $\Q$. By assumption, $\zeta_{\widetilde{K}}(s) = L(s,\pi)$ for some automorphic representation $\pi$ of $\mathrm{GL}_{m}(\mathbb{A}_{\Q})$ with $m = [\widetilde{K}:\Q] \leq n!$. Clearly, $L(s,\pi)$ satisfies GRC. Let 
\[
Q = \max\{D_{\widetilde{K}}, 2d\}, \qquad 
0 < \epsilon < \frac{1}{4 \ell(n-1)} \qquad
\delta  = \Big(\frac{1}{2\ell(n-1)} - \epsilon\Big)\frac{\log D_K}{\log D_{\widetilde{K}}},
\]
and let $H_{\delta,\epsilon} \geq 1$ be sufficiently large. From the estimate $D_K^{|G|/n} \leq D_{\widetilde{K}} \leq D_K^{|G|/2}$, one can see that $\delta < 1$ and $\delta$ is bounded away from zero {\it uniformly} in terms of $n, \ell,$ and $\epsilon$. Thus, when a quantity depends on $\delta$ (such as $H_{\delta,\epsilon}$), we may replace this dependence with $n, \ell,$ and $\epsilon$. In particular, we may treat $\delta$ as independent of $D_K$ and $D_{\widetilde{K}}$. 

We note that while $\zeta_{\widetilde{K}}(s)$ does not directly satisfy the hypotheses of \cref{thm:Landau_Siegel} (as $\zeta_{\widetilde{K}}(s)$ is not the $L$-function of a cuspidal automorphic representation $\mathbb{Q}$, and it only conjecturally factors into a product of such $L$-functions), it has an analytic continuation and functional equation just as described in \cref{sec:L-functions}.  The only part of the proof of \cref{thm:Landau_Siegel} which relies on cuspidality is in the use of \cref{thm:large_sieve}.  However, since we are considering a single $L$-function here instead of several, the use of \cref{thm:large_sieve} can be replaced with the field-uniform analogue of the Brun-Titchmarsh theorem proved in \cite[Proposition 2]{GMP}.

We apply \cref{thm:Landau_Siegel} with $\pi_0$ trivial,  $\cF_m(Q) = \{ \pi\}$, $T = 1$, and $\sigma = 1 - \frac{H_{\delta,\epsilon}}{\log D_{\widetilde{K}} }$ to deduce that
\[
N_{\pi}(1 - \tfrac{H_{\delta,\epsilon}}{\log D_{\widetilde{K}} }, 1 ) \ll_n \big( (1-\beta_{\chi}) \log Q \big) Q^{10^7 m^4 H_{\delta,\epsilon} / \log D_{\widetilde{K}} }. 
\]  
Since we have $\log Q  \asymp_{n,\epsilon,\ell} \log D_{\widetilde{K}} \asymp_{n,\epsilon,\ell} \log d$, the bound $m \leq n!$ and $\beta_{\chi} = 1 - \frac{\eta_{\chi}}{\log d}$,  we have that $N_{\pi}(1 - \tfrac{H_{\delta,\epsilon}}{\log D_{\widetilde{K}} }, 100 ) \ll_{n, \ell, \epsilon}  \eta_{\chi}$. As $\eta_{\chi}$ sufficiently small depending only on $n, \ell, \epsilon$, it follows that $\zeta_{\widetilde{K}}(s)$ has no zeros in the region
\[
\Re(s) > 1 - \frac{H_{\delta,\epsilon}}{\log D_{\widetilde{K}}}, \qquad |\Im(s)| \leq 1. 
\]
Thus, by \cref{prop:ell-torsion_primes}, there are $M$ rational primes $p \leq D_{\widetilde{K}}^{\delta} = D_K^{\frac{1}{2\ell(n-1)} - \epsilon}$ which split completely in $\widetilde{K}$ with 
\[
M \gg_{\epsilon,n,\ell} D_{\widetilde{K}}^{\delta-\epsilon} = D_K^{\frac{1}{2\ell(n-1)}-2\epsilon}, 
\] 
provided $D_K$ is sufficiently large depending on $\epsilon, n,$ and $\ell$. The result now follows from an application of \cref{lem:SmallSplitPrimes} and rescaling $\epsilon$ appropriately. 
\end{proof}

\appendix
\label{sec:appendix}

\section{Explicit upper bound on the universal family for $\GL_n$}

Let $F$ be a number field of degree $d$ over $\Q$ and discriminant $D$ and let $n\geq 1$ be an integer. Let $\mathscr{A}_{\rm cusp}$ denote the set of unitary cuspidal automorphic representations $\pi$ of $\GL_n(\A_F)$, with normalized central character, ordered by analytic conductor $C(\pi)$.  We recall that $C(\pi)=D^n N_\pi k_{\kp}i$, where $N_\pi={\rm Norm}(\kq_{\pi})$ is the arithmetic conductor and $k_{\kp}i$ the archimedean conductor, as in \eqref{eqn:analytic_conductor_def}. Note the factor of the discriminant, which arises naturally in the functional equation for the standard $L$-function of $\pi$.

For $Q\geq 1$ let
\[
\mathscr{F}(Q)=\{\pi\in\mathscr{A}_{\rm cusp}: C(\pi)\leq Q\}.
\]
We present an argument, due to Venkatesh \cite{Venk} and based on results in \cite{Brumley}, to deduce a polynomial upper bound on the cardinality $|\mathscr{F}(Q)|$. We can in fact make this polynomial bound explicit, using subsequent refinements of \textit{loc. cit.}, as in the following

\begin{thm}\label{thm:poly-growth} We have, for all fixed $\epsilon>0$, $|\mathscr{F}(Q)|\ll_{d,n,\epsilon} (D^{-n^2}Q^{2n})^{1+\epsilon}$.
\end{thm}

\begin{remark}
As we consider $d$ and $n$ as being fixed, we shall henceforth systematically suppress the dependence of implied constants on $n$ and $d$ in the notation.
\end{remark}
\begin{remark}
The expected value of the exponent of $Q$ in Theorem \ref{thm:poly-growth} is $n+1$, and indeed this was shown (with an asymptotic) in \cite{BM}, with one caveat: for $n\geq 3$ the authors restrict to the subfamily of $\mathscr{F}(Q)$ consisting of Maass forms. This restriction is fortunate, in a way, since it provides an occasion for this appendix, which has sat for a long time in a drawer (or inbox) and whose methods are quite different. While Theorem \ref{thm:poly-growth} says nothing about existence, and the upper bound is not sharp, we believe that the proof itself is of sufficient interest to merit circulation. 
\end{remark}

\begin{remark}
The results of \cite{BM} make no claim of uniformity in the number field $F$. (In fact one should note the difference in the notational conventions between these two papers: in \cite{BM}, the analytic conductor, denoted by $Q(\pi)$ there, does not include the factor of the discriminant.) The upper bound in Theorem \ref{thm:poly-growth} is, however, uniform in $D$, making this perhaps the most novel aspect of the result.
\end{remark}

The proof of Theorem \ref{thm:poly-growth} combines two ingredients: Rankin-Selberg theory and sphere packing bounds in large dimensions. It is natural to ask what effect assuming standard conjectures on these $L$-functions would have on the quality of the resulting bound. For example, a similar argument to the one we present here was used in \cite{FKM} to count $\ell$-adic sheaves of bounded complexity. In that article, Deligne's proof of the Riemann hypothesis over finite fields is used  to show that certain trace functions form a quasi-orthogonal system with small enough angular separation to deduce a polynomial upper bound. We show that the exponent $2n$ can be improved to $n+1$ under standard conjectures, demonstrating the strength of the method of proof. 

\begin{thm}\label{under-RR}
Denote by $\mathscr{F}_\chi(Q)$ the subfamily of $\mathscr{F}(Q)$ having fixed central character $\chi$. Assume the Ramanujan conjecture and the Riemann hypothesis for Rankin-Selberg $L$-functions. Then $|\mathscr{F}_\chi(Q)|\ll_\epsilon (D^{-n^2/2}Q^n)^{1+\epsilon}$, and
\[
|\mathscr{F}(Q)|\ll_\epsilon (DQ)^\epsilon D^{-n^2/2-n}Q^{n+1}.
\]
\end{thm}

\begin{remark}
Note that, by the results in \cite{BM}, the exponent of $Q$ in Theorem \ref{under-RR} is sharp, up to the $\epsilon$. Moreover, the $D$ dependence here and that of the \textit{main term} of the asymptotic given in \cite{BM} are in agreement.
\end{remark}
\begin{remark}\label{BH-assumptions}
The method of proof of Theorems \ref{thm:poly-growth} and \ref{under-RR} is sensitive to any loss of information incurred in the application of the Bushnell-Henniart bounds \cite{BH}. Recall that the main result in \textit{loc. cit.} provides upper bounds for the Rankin-Selberg Artin exponent $\Ar(\pi_v\times\tilde\pi_v')$ at finite places $v$ in terms the standard Artin exponents $\Ar(\pi_v)$ and $\Ar(\tilde\pi'_v)$, and the integers $n,n'$, where $\pi_v$ and $\pi_v'$ are smooth irreducible representations of $\GL_n(F_v)$ and $\GL_{n'}(F_v)$, respectively.

While the bounds in \textit{loc. cit.} are sharp in general, we apply them under additional hypotheses on $\pi_v$ and $\pi_v'$. Namely, in the course of the proof, we assume that 
\begin{enumerate}
\item\label{dim} the dimensions $n=n'$ are the same, 
\item\label{cond-exp} the Artin exponents $a=\Ar(\pi)=\Ar(\pi')$ are the same, 
\item\label{cent-char} the central characters are the same, say equal to $\chi$.
\end{enumerate}
Under the assumptions \eqref{dim} and \eqref{cond-exp} above, Theorem 1 in \cite{BH} establishes the sharp bound $\Ar(\pi_v\times\tilde\pi_v')\leq (2n-1)a$. In Theorem \ref{thm:improvedBH} of Appendix \ref{AppendixB}, Bushnell and Henniart show that, under the additional assumption of \eqref{cent-char}, this bound can be improved to $\Ar(\pi_v\times\tilde\pi_v')\leq (2n-2)a$.

This improved bound is an ingredient in the explicit exponents given in Theorems \ref{thm:poly-growth} and \ref{under-RR}. Without this improvement, the unconditional bound in Theorem \ref{thm:poly-growth} would have an additional factor of $Q$, and the conditional bound in Theorem \ref{under-RR} would have an additional factor of $Q^{1/2}$.
\end{remark}

\begin{remark}
The method of proof of Theorems \ref{thm:poly-growth} and \ref{under-RR} requires fixing certain representation theoretic data, of combinatorial nature. This data encodes the dimensional blocks of the inducing supercuspidal representations in the Bernstein-Zelevinsky classification, as well as the partition of these blocks according to the underlying twist equivalency classes. See \S\ref{sec:RS-local} for more details. After bounding the size of the subfamily associated with such data, one then sums over the finite number of such choices.

This decomposition allows one to prove, in principle, refined bounds for the cardinality of these subfamilies, since the Bushnell-Henniart bounds \cite{BH} can often be improved under such assumptions. For example, if the combinatorial data that one takes is ``trivial'', in the sense that it corresponds to $\pi_v$ and $\pi_v'$ supercuspidal on $\GL_n$, then (keeping the assumptions \eqref{cond-exp} and \eqref{cent-char} of the previous remark) one can use the bound $\Ar(\pi\times\tilde\pi')\leq na$ of \cite[Corollary C]{BH2017}, which is, in general, far better than the general bound of $(2n-2)a$ cited above. In this way one can show that, under Ramanujan and Riemann as in Theorem \ref{under-RR}, the subfamily of $\mathscr{F}(Q)$ consisting of $\pi$ which 
\begin{enumerate} 
\item are supercuspidal at all the places at which they ramify,
\item have  archimedean component lying in some fixed compact of the unitary dual, 
\end{enumerate}
has cardinality $O(Q^{\frac{n}{2}+2})$ (ignoring the discriminant dependence). This bound is surprisingly strong, and no trace formula was used to derive it. We have not found this type of interplay between conductor dropping phenomenon and improved bounds on dimension counts of automorphic forms elsewhere in the literature.

\end{remark}

\subsection{Idea of proof}\label{sec:outline}
We present here the basic argument to prove Theorem \ref{thm:poly-growth}. We shall later need to modify the presentation to obtain the best possible exponent.

Let $\kq$ be an integral ideal of $\mathcal{O}_F$. Let $\chi$ be a character of $\A_f^\times$ of conductor $\kq$, where $\A_f$ is the ring of finite adeles. Let
\[
\mathscr{A}_{\rm cusp}(\kq,\chi)=\{\pi\in\mathscr{A}_{\rm cusp}: \kq_{\pi_f}=\kq, \; \chi_{\pi_f}=\chi\}, \qquad \mathscr{F}_{\kq,\chi}(Q)=\mathscr{F}(Q)\cap \mathscr{A}_{\rm cusp}(\kq,\chi).
\]
Here, $\kq_{\pi_f}$ is the conductor of $\pi_f$ and $\chi_{\pi_f}$ is the central character of $\pi_f$. Then
\begin{equation}\label{LevelSum}
|\mathscr{F}(Q)|=\sum_{{\rm Norm}(\kq)\leq Q/D^n}\sum_{\substack{\chi\\ \text{cond} \;\kq}}  |\mathscr{F}_{\kq,\chi}(Q)|.
\end{equation}
The argument we sketch below provides a bound on $|\mathscr{F}_{\kq,\chi}(Q)|$ of the form
\[
O_\epsilon((D^{-n^2} {\rm Norm}(\kq)^{-2} Q^{2n+n^2})^{1+\epsilon}).
\]
Executing double sum over pairs $(\kq,\chi)$, this would produce a bound of $O_\epsilon((D^{-n^2}Q^{2n+n^2})^{1+\epsilon})$. We will later show how to remove the $n^2$ to establish Theorem \ref{thm:poly-growth}, as well as the sharp conditional bounds in Theorem \ref{under-RR}.

\subsubsection{Mapping $\mathscr{A}_{\rm cusp}(\kq,\chi)$ to a Hermitian space} 

We begin by describing a way to map $\mathscr{A}_{\rm cusp}(\kq,\chi)$ to a Hermitian space, whose inner product can be understood in terms of Rankin-Selberg $L$-functions. The reader is encouraged to read ahead to the next subsection describing the Dirichlet coefficients of these $L$-functions for a motivation of the following constructions. 

Recall that a partition $\mu=(\mu_i)$ is a sequence of non-increasing non-negative integers $\mu_1\geq\mu_2\geq\cdots $ with only finitely many non-zero entries. Write $\mathcal{P}$ for the set of all partitions. The \textit{length} of $\mu\in\mathcal{P}$, denoted $\ell(\mu)$, is the number of its non-zero entries. Write
\[
\mathcal{P}_\ell=\{\mu:\mu_1\geq\mu_{2}\geq\cdots\geq\mu_{\ell}\geq 0\}
\]
for the partitions of length at most $\ell$. Finally, for $\mu=(\mu_i)\in\mathcal{P}$, write $|\mu|=\sum_i \mu_i$. For an integer $r$, let $\mathcal{P}_\ell(r)=\{\mu\in\mathcal{P}_\ell: |\mu|=r\}$; this is the empty set when $r$ is negative.

Let $S$ be a finite set of finite places. Let $I^S$ denote the set of integral ideals of $\mathcal{O}_F$ supported outside of $S$. When $S$ is empty we abbreviate this to $I$ for the set of all integral ideals. Given an $\mathfrak{n}=\prod_{\kp} \kp^{r_{\kp}}\in I$ we write $\mathscr{P}_{n-1}(\mathfrak{n})$ for the set of sequences $\bbmu=(\mu_\kp)_{\kp}$ of partitions such that $\mu_{\kp}\in\mathcal{P}_{n-1}(r_{\kp})$. A $\mathcal{P}_{n-1}$\textit{-decorated prime-to-$S$ ideal} is a pair $(\mathfrak{n},\bbmu)$, where $\mathfrak{n}\in I^S$ and $\bbmu\in\mathscr{P}_{n-1}(\mathfrak{n})$. Let $\mathscr{I}^S$ denote the set of $\mathcal{P}_{n-1}$-decorated prime-to-$S$ ideals. We have a map $\mathscr{I}^S\rightarrow I^S$, $(\mathfrak{n},\bbmu)\mapsto \mathfrak{n}$, where we forget the decoration and take the underlying ideal. Observe that several $(\mathfrak{n},\bbmu)$ can have the same underlying ideal $\mathfrak{n}$. We shall sometimes write $\tilde{\mathfrak{n}}$ for a $\mathcal{P}_{n-1}$-decorated ideal with underlying ideal $\mathfrak{n}$.

For a parameter $X>1$, let $\mathscr{I}^S(X)=\{\tilde{\mathfrak{n}}\in\mathscr{I}^S: {\rm Norm}(\mathfrak{n})\leq  X\}$; this is the set of pairs $(\mathfrak{n},\bbmu)$ with $\mathrm{Norm}(\mathfrak{n})\leq X$ and $\bbmu\in\mathscr{P}_{n-1}(\mathfrak{n})$. Let $V^S(X)$ be the vector space of complex valued functions on $\mathscr{I}^S(X)$. Endow $V^S(X)$ with the standard scalar product
\[
\langle f,g\rangle=\sum_{\tilde{\mathfrak{n}}\in\mathscr{I}^S(X)} f( \tilde{\mathfrak{n}})\overline{g(\tilde{\mathfrak{n}})}.
\]
For an integral ideal $\kq$ of $\mathcal{O}_F$, with support $S$, we shall map $\mathscr{A}_{\rm cusp}(\kq,\chi)$ to $V^S(X)$ in the following way.

For a partition $\mu\in\mathcal{P}_{n-1}$ let $s_\mu$ denote the associated Schur function in $n$ variables. If $\mu$ is the zero partition, then $s_\mu$ is identically $1$. For $(\mathfrak{n},\bbmu)\in\mathscr{I}^S$ pose
\begin{equation}\label{prime2S}
a_\pi(\mathfrak{n},\bbmu)=\prod_{\kp\notin S} s_{\mu_{\kp}}(A_\pi(\kp)),\qquad \text{where}\quad A_{\pi}(\kp)=(\alpha_{1,\pi}(\kp),\ldots,\alpha_{n,\pi}(\kp)).
\end{equation}
We note that if $n=2$ and $\pi_f$ has trivial central character, then there is no decoration $\bbmu$ and the \eqref{prime2S} just recovers the Hecke eigenvalue of $\pi$ at $\mathfrak{n}$. In fact, more generally, when $n\geq 2$ and $\pi_f$ has trivial central character, if we take $\bbmu=(\mu_{\kp})_{\kp}$ to satisfy $\mu_{\kp}=(r_{\kp},0,\ldots )$, then we once again recover the Hecke eigenvalue at $\mathfrak{n}=\prod_{\kp}\kp^{r_{\kp}}$.

Let $f:\R\rightarrow\R$ be a non-negative smooth function supported in $\left[\frac12,1\right]$ and having Lebesque integral $1$. Write
\[
F_X^S(\mathfrak{n})=\sum_{(\mathfrak{m},S)=1}f({\rm Norm}(\mathfrak{nm}^n)/X).
\]
For every $\pi\in\mathscr{A}_{cusp}(\kq,\chi)$ we define a vector $\mathbf{v}_\pi^S\in V^S(X)$ by the rule
\[
\mathbf{v}_\pi^S: (\mathfrak{n},\bbmu)\mapsto \sqrt{F_X^S(\mathfrak{n})}a_\pi(\mathfrak{n},\bbmu).
\]

Note that for ${\rm Norm}(\mathfrak{n})>X$ we have $F_X^S(\mathfrak{n})=0$; in this way the function $\tilde{\mathfrak{n}}\mapsto \mathbf{v}_\pi^S(\tilde{\mathfrak{n}})$ can indeed be viewed as an element of $V^S(X)$. 

\subsubsection{Relation to Rankin-Selberg $L$-functions}
We now recall the description of the Rankin-Selberg Dirichlet coefficients. This will clarify the choice of map $\pi\mapsto \mathbf{v}_\pi^S$ and the inner product we put on $V^S(X)$.

Let $\pi,\pi'\in\mathscr{A}_{\rm cusp}(\kq,\chi)$. The prime-to-$S$ part of the Rankin-Selberg $L$-function is defined, for ${\rm Re}(s)>1$, by the Euler product
\[
L^S(s,\pi\times\tilde\pi')=\prod_{\kp\notin S}\prod_{j=1}^n \prod_{j'=1}^{n}(1-\alpha_{j,\pi}(\kp)\overline{\alpha_{j',\pi'}(\kp)}{\rm Norm}(\kp)^{-s})^{-1}.
\]
We write $a_{\pi\times\tilde\pi'}(\mathfrak{n})$, for $(\mathfrak{n},S)=1$, for the Dirichlet coefficients of $L^S(s,\pi\times\tilde\pi')$, so that
\[
L^S(s,\pi\times\tilde\pi')=\sum_{(\mathfrak{n},S)=1}a_{\pi\times\tilde\pi'}(\mathfrak{n})\mathrm{Norm}(\mathfrak{n})^{-s}.
\]
Cauchy's identity shows that
\[
a_{\pi\times\tilde\pi'}(\kp^r)=\sum_{\mu\in\mathcal{P}_n(r)}s_\mu(A_\pi(\kp))s_\mu(A_{\tilde\pi'}(\kp)).
\]
Following the exposition in \cite[\S 2]{Brumley_2}, for a partition $\mu=(\mu_1,\ldots ,\mu_{n-1},k,0,\ldots )\in\mathcal{P}_n$ we let $\widehat\mu=(\mu_1-k,\ldots ,\mu_{n-1}-k,0,\ldots)\in\mathcal{P}_{n-1}$. Then $s_\mu(A_\pi(\kp))=\chi^k(\varpi_{\kp})s_{\hat\mu}(A_\pi(\kp))$. Now, for any pair $(\mu,k)$, where $\mu\in\mathcal{P}_{n-1}$ and $k\geq 0$, there is a unique $\lambda\in\mathcal{P}_n$ such that $|\lambda|=|\mu|+kn$ and $\widehat{\lambda}=\mu$ (add $k$ to each of the first $n$ entries of $\mu$). Applying this we get
\begin{equation}\label{unram-RS-coeff}
a_{\pi\times\tilde\pi'}(\kp^r)=\sum_{k\geq 0}\sum_{\mu\in\mathcal{P}_{n-1}(r-nk)}s_{\mu}(A_\pi(\kp))s_{\mu}(A_{\tilde\pi'}(\kp)).
\end{equation}
The sum on $k$ is finite, going up to the integer part of $r/n$. Note that, in the above expression, we have used the fact that $\chi_{\pi_f}=\chi_{\pi_f'}=\chi$; this explains why we've decomposed according to central character in \eqref{LevelSum}. Thus, for $(\mathfrak{n},S)=1$, we have
\begin{align}\label{n-prime2S}
a_{\pi\times\tilde\pi'}(\mathfrak{n})&=\prod_{\kp^{r_{\kp}}\|\mathfrak{n}}\sum_{k_{\kp}\geq 0}\sum_{\mu_{\kp}\in\mathcal{P}_{n-1}(r_{\kp}-nk_{\kp})}s_{\mu_{\kp}}(A_\pi(\kp))s_{\mu_{\kp}}(A_{\tilde\pi'}(\kp))\nonumber\\
&=\sum_{\substack{(\mathfrak{m},S)=1\\\mathfrak{m}^n\mid\mathfrak{n}}}\sum_{\bbmu\in\mathscr{P}_{n-1}(\mathfrak{n}/\mathfrak{m}^n)}a_\pi(\mathfrak{n}/\mathfrak{m}^n,\bbmu)a_{\pi'}(\mathfrak{n}/\mathfrak{m}^n,\bbmu).
\end{align}

We now consider the smooth sum of coefficients 
\[
S(X)=\sum_{(\mathfrak{a},S)=1}a_{\pi\times\tilde\pi'}(\mathfrak{a})f({\rm Norm}(\mathfrak{a})/X).
\]
We have
\begin{align*}
S(X)&=\sum_{(\mathfrak{a},S)=1}f({\rm Norm}(\mathfrak{a})/X)\sum_{\substack{(\mathfrak{m},S)=1\\\mathfrak{m}^n\mid\mathfrak{a}}}\sum_{\bbmu\in\mathscr{P}_{n-1}(\mathfrak{a}/\mathfrak{m}^n)}a_\pi(\mathfrak{a}/\mathfrak{m}^n,\bbmu)\overline{a_{\pi'}(\mathfrak{a}/\mathfrak{m}^n,\bbmu)}\\
&=\sum_{(\mathfrak{n},S)=1}\sum_{(\mathfrak{m},S)=1}f({\rm Norm}(\mathfrak{nm}^n)/X)\sum_{\bbmu\in\mathscr{P}_{n-1}(\mathfrak{n})}a_\pi(\mathfrak{n},\bbmu)\overline{a_{\pi'}(\mathfrak{n},\bbmu)}\\
&=\sum_{(\mathfrak{n},\bbmu)\in\mathscr{I}^S}F_X^S(\mathfrak{n})a_\pi(\mathfrak{n},\bbmu)\overline{a_{\pi'}(\mathfrak{n},\bbmu)}\\
&=\sum_{\tilde{\mathfrak{n}}\in\mathscr{I}^S(X)}\mathbf{v}_\pi^S(\tilde{\mathfrak{n}})\overline{\mathbf{v}_\pi^S(\tilde{\mathfrak{n}})}.
\end{align*}
We recognize this as $\langle \mathbf{v}_\pi^S,\mathbf{v}_{\pi'}^S\rangle$. On the other hand, if we let
\[
\hat{f}(s)=\int_0^\infty f(x)x^s \frac{dx}{x}
\]
be the Mellin transform of $f$, then by the Mellin inversion formula one has
\[
S(X)=\frac{1}{2\pi i}\int_{(2)}L^S(s,\pi\times\tilde\pi')\hat{f}(s)X^sds.
\]
This allows us to read off the orthogonality properties of $\mathbf{v}_\pi^S$ and $\mathbf{v}_{\pi'}^S$ in terms of the analytic information of $L^S(s,\pi\times\tilde\pi')$.

\subsubsection{Strategy of proof}
Let
\[
\mathbf{u}_\pi^S=\frac{\mathbf{v}_\pi^S}{\langle \mathbf{v}_\pi^S,\mathbf{v}_\pi^S\rangle^{1/2}}
\]
be the projection of the vector $\mathbf{v}_\pi^S$ to the unit sphere in $V^S(X)$. The idea behind the proof of Theorem \ref{thm:poly-growth} is to show that, for $X$ large relative to $Q$,
\begin{enumerate}
\item\label{step1} the map $\mathscr{F}_{\kq,\chi}(Q)\rightarrow V^S$ given by $\pi\mapsto \mathbf{v}_\pi^S$ is injective;
\item\label{step2} when $\pi,\pi'\in\mathscr{F}_{\kq,\chi}(Q)$ are distinct, the vectors $\mathbf{u}_\pi^S$ and $\mathbf{u}_{\pi'}^S$ are quasi-orthogonal;
\item\label{step3} there cannot be too many such quasi-orthogonal vectors.
\end{enumerate}
Moreover, each of these steps will be seen to be quantifiable, polynomially in $Q$.

There is only one problem with this approach: we have thrown out the information at ramified primes. While this allows for a simpler presentation, the price to pay is a weaker bound in Theorem \ref{thm:poly-growth}. Indeed one obtains in this way the exponent $2n+n^2+\epsilon$ in the parameter $Q$, with or without assuming the Ramanujan conjecture and the Riemann hypothesis. See Remark \ref{rem:n-square} for more details on the source of this loss by a power of $n^2$.

To obtain the unconditional bound of Theorem \ref{thm:poly-growth} (as well as the conditional bound of Theorem \ref{under-RR}, which is sharp up to $\epsilon$), we shall need to take into account the information at ramified primes. To adapt the above argument along these lines, one must explicate the Rankin-Selberg coefficients at ramified primes, which has been done by Brumley in \cite[Appendix]{ST}. In particular, we shall see in \S\ref{sec:RS-local} that the ``combinatorial distance to supercuspidal'' of $\pi_S=\otimes_{\kp\in S}\pi_{\kp}$ governs the shape of the ramified Rankin-Selberg coefficients. Then, in \S\ref{sec:comb-decomp}, we further decompose $\mathscr{F}_{\kq,\chi}(Q)$ according to this data. After an appropriate enrichening of the space $V^S(X)$ to take into account this information, we then execute the above three steps.

\subsection{Rankin-Selberg theory}\label{sec:RS-local}

We now recall some of the basic local and global properties of the Rankin-Selberg $L$-function that we shall need in the proof of Theorem \ref{thm:poly-growth}.

\subsubsection{Induction data}\label{induction-data}
Let $v$ be a finite place of $F$ associated with a prime ideal $\kp$ of $\mathcal{O}_F$. Let $q_v$ be the cardinality of the residue field. Let $\pi_v$ be an irreducible unitary generic representation of $\GL_n(F_v)$.

Recall that by the Bernstein-Zelevinsky description of admissible dual, we may associate with $\pi_v$ (see \cite[\S A.2]{ST}) the following \textit{combinatorial} data:
\begin{enumerate}
\item[(C1)] a standard Levi subgroup $M\simeq \GL_{n_1}\times\cdots\times\GL_{n_r}$ of $\GL_n$;
\item[(C2)] a partition $\underline{J}=[J_1,\ldots ,J_A]$ of the set $\{1,\ldots ,r\}$;
\item[(C3)] an integer vector $\mathbf{d}=(d_1,\ldots ,d_r)\in \N^r$, where $d_j\mid n_j$, such that $m_j=n_j/d_j$ is constant (say equal to $m_a$) along $j\in J_a$;
\item[(C4)] an integer vector $\mathbf{e}=(e_1,\ldots ,e_A)\in \N^A$, where each $e_a\mid n$;
\end{enumerate}
the following \textit{analytic data}:
\begin{enumerate}
\item[(A1)] real numbers $\sigma_1\geq \cdots \geq \sigma_r$;
\item[(A2)] real numbers $t_1,\ldots ,t_r$;
\end{enumerate}
encoded in the complex numbers
\[
s_j=\sigma_j+it_j\qquad\text{and}\qquad z_j=q_v^{-s_j-n_j/2};
\]
as well as  the following \textit{arithmetic} data:
\begin{enumerate}
\item[(SC)] a set $\{\varrho_1,\ldots ,\varrho_A\}$ of pairwise twist-inequivalent unitary supercuspidal representations $\varrho_a$ of $\GL_{m_a}(F_v)$ having torsion number $e_a$.
\end{enumerate}

\subsubsection{Rankin-Selberg local factors}\label{RS-local-factor}
The local Rankin-Selberg $L$-factor can be expressed using the above combinatorial and analytic data. (The epsilon factor, on the other hand, encodes the arithmetic information contained in the choice of supercuspidal representations on each block. We do not define the epsilon factors here, but they are used implicitly in Appendix \ref{AppendixB}.)  We let $\mathrm{Comb}_v=\{(M,\underline{J}, \mathbf{d}, \mathbf{e})\}$ denote the collection of combinatorial data C1, C2, C3, C4. Let $\pi_v$ and $\pi_v'$ both have the same combinatorial type $(M,\underline{J}, \mathbf{d}, \mathbf{e})$. Let $z_j, z_j'$ denote their respective analytic data. 

By \cite[\S A.2, Example 1]{ST} we have
\begin{align}
L(s,\pi_v\times\tilde\pi_v')&=\prod_{a=1}^A\; \prod_{j,k\in J_a}\prod_{\nu=1}^{\min(n_j,n_k)}\left(1-(q_v^\nu z_j\overline{z_k'})^{e_a}q_v^{-e_as}\right)^{-1}\nonumber\\
&=\prod_{a=1}^A\prod_{\nu=1}^n \prod_{j,k\in J_a^\nu}\left(1-(q_v^\nu z_j\overline{z_k'})^{e_a}q_v^{-e_as}\right)^{-1},\label{RS-L-factor}
\end{align}
where $J_a^\nu=\{j\in J_a: n_j\geq\nu\}$. We expand the expression \eqref{RS-L-factor} into the local Dirichlet series, which we again denote by $a_{\pi\times\pi'}(\kp^r)$. We shall now describe these in terms of the analytic data $z_j$, similarly to the unramified setting of \S\ref{sec:outline}.

We now furthermore assume that the central characters of $\pi_v$ and $\pi'_v$ coincide. We fix $a$ and $\nu$ in \eqref{RS-L-factor} and expand the product over $j$ and $k$. We obtain
\[
\prod_{j,k\in J_a^\nu}(1-(q_v^\nu z_j\overline{z_k'})^{e_a}X^{e_a})^{-1}=\sum_{r\geq 0}a_{\pi\times\tilde\pi'}(\kp^{e_ar};\nu,a)X^{e_ar}.
\]
Cauchy's identity will once again allow us to describe these coefficients $a_{\pi\times\tilde\pi'}(\kp^{e_ar};\nu,a)$ as a combinatorial expression in terms of the local roots. With this in mind, we let $A_\pi(\kp;a,\nu)$ denote the set of parameters $q_v^{\nu/2}z_j$, for $j\in J_a^\nu$, completed to a size $n$ multiset by adding $n-|J_a^\nu|$ remaining zeros. For an integer $e\geq 1$ we write $A_\pi^e(\kp;a,\nu)$ for the set of $e$-th powers of the parameters in $A_\pi(\kp;a,\nu)$. We may then evaluate the Schur functions in $n$ variables on $A_\pi^e(\kp,a,\nu)$. Reasoning as in \eqref{unram-RS-coeff}, we find
\begin{equation}\label{lambda-nu-a}
a_{\pi\times\tilde\pi'}(\kp^{e_ar};\nu,a)=\sum_{k\geq 0}\sum_{\mu\in\mathcal{P}_{n-1}(r-nk)}s_{\mu}(A_\pi^{e_a}(\kp;a,\nu))s_{\mu}(A_{\tilde\pi'}^{e_a}(\kp;a,\nu)).
\end{equation}
Multiplying out $\nu$ and $a$ in \eqref{RS-L-factor}, we deduce that $a_{\pi\times\tilde\pi'}(\kp^r)$ is the complete homogeneous polynomial of degree $r$ in the coefficients $a_{\pi\times\tilde\pi'}(\kp^{e_af};\nu,a)$.

\begin{remark}
The combinatorial data $M=T$, $\underline{J}=\{1,\ldots ,n\}$, $\mathbf{d}=(1,\ldots ,1)$, and $\mathbf{e}=1$ corresponds to representations $\pi_v$ which are, up to a character twist, unramified. In this case, the coefficient $a_{\pi\times\pi'}(\kp^r;\nu,1)$ is zero for all $\nu>1$, since all $n_j=1$. Thus $a_{\pi\times\pi'}(\kp^r)=a_{\pi\times\pi'}(\kp^r;1,1)$. Note that when $\pi_v$ is unramified, $A_\pi(\kp;1,1)$ is the set of the Satake parameters $A_\pi(\kp)$, and \eqref{lambda-nu-a} recovers \eqref{unram-RS-coeff}. 
\end{remark}

\begin{remark}\label{rem:local-roots}
In \cite[(A.6)]{ST}, it is shown that when $\pi_v$ and $\pi_v'$ are irreducible unitary generic representations of $\GL_n(F_v)$ and $\GL_m(F_v)$, respectively, then
\begin{equation}\label{RS-L-factor2}
L(s,\pi_v\times\pi_v')=\prod_{(a,b)\in\Delta}\; \prod_{j\in J_a}\prod_{k\in K_b}\prod_{\nu=1}^{\min(n_j,n_k')}(1-(q_v^\nu z_jz_k')^{e_{\ell(a,b)}}q_v^{-e_{\ell(a,b)}s})^{-1}.
\end{equation}
The expression \eqref{RS-L-factor} is a special case of this, when both $\pi_v$ and $\pi_v'$ have the same combinatorial type. See \textit{loc. cit.} for relevant notation.

The local roots $q_v^\nu z_jz_k'$ in \eqref{RS-L-factor2} satisfy $|q_v^\nu z_jz_k'|=q_v^{\nu-\sigma_j-\sigma_k'-n_j/2-n_k'/2}$.  Under the Ramanujan conjecture, we have $\sigma_j=\sigma_k'=0$, so that 
\begin{equation}\label{RS-tempered}
|q_v^\nu z_jz_k'|=q_v^{\nu-n_j/2-n_k'/2}\leq 1.
\end{equation}
Unconditionally, the Jacquet-Shalika bounds \cite{JS} show that $0\leq |\sigma_j|,|\sigma_k'|<1/2$, so that 
\begin{equation}\label{JS-RS}
|q_v^\nu z_jz_k'| <q_v^{\nu+1-n_j/2-n_k'/2}\leq q_v.
\end{equation}
Rudnick-Sarnak \cite[Appendix]{RS} improved this to $q_v^{1-\delta}$, where $\delta=1/(n^2+1)+1/(m^2+1)$.\end{remark}

\subsubsection{General formula for Dirichlet coefficients}

We put together the descriptions of the prime-to-$S$ coefficients in \eqref{prime2S} with the ramified coefficients in \eqref{lambda-nu-a}. 

We continue to write $v$ for a finite place with associated prime ideal $\kp$. Recall the set $\mathrm{Comb}_v$ from \S\ref{RS-local-factor}, whose elements index the combinatorial data $\mathscr{C}_v=(M_v,\underline{J}_v, \mathbf{d}_v, \mathbf{e}_v)$ described in \S\ref{induction-data}. Via the expansion \eqref{RS-L-factor}, $\mathscr{C}_v$ gives rise to a set 
\[
\{A_\pi^{e_v}(\kp;a_v,\nu_v): 1\leq a_v\leq A_v, 1\leq \nu_v\leq n\},
\]
encoding the analytic data. We shall write ${\rm Index}(\mathscr{C}_v)$ for the indexing set of pairs $(a_v,\nu_v)$. 

Now let $S$ once again denote the prime support of the ideal $\kq$ and put $\mathrm{Comb}_S=\prod_{v\in S}\mathrm{Comb}_v$. For any $\mathscr{C}\in\mathrm{Comb}_S$, we let $\mathscr{A}_{\rm cusp}(\kq,\chi,\mathscr{C})$ denote the set of $\pi\in\mathscr{A}_{\rm cusp}(\kq,\chi)$ such that $\pi_S$ has combinatorial data $\mathscr{C}$. Let $\pi,\pi'\in\mathscr{A}_{\rm cusp}(\kq,\chi,\mathscr{C})$. Let $\mathfrak{n}$ be an integral ideal and $\bbmu\in\mathscr{P}_{n-1}(\mathfrak{n})$. Let $(a,\nu)\in {\rm Index}(\mathscr{C})$. Generalizing \eqref{prime2S}, we write
\[
a_\pi(\mathfrak{n},\bbmu;a,\nu)=\prod_{\kp\notin S}s_{\mu_{\kp}}(A_\pi(\kp))\prod_{\kp\in S}s_{\mu_{\kp}}(A_\pi^{e_{a_v}}(\kp;a_{\kp},\nu_{\kp})).
\]
Then the Dirichlet coefficients of $L(s,\pi\times\pi')$, denoted $a_{\pi\times\pi'}(\mathfrak{n})$, can be written as
\begin{equation}\label{full-coeff}
a_{\pi\times\tilde\pi'}(\mathfrak{n})=\sum_{\mathfrak{m}^n\mid\mathfrak{n}}\; \sum_{\bbmu\in\mathscr{P}_{n-1}(\mathfrak{n}/\mathfrak{m}^n)}\;\sum_{(a,\nu)\in \mathscr{C}}a_\pi(\mathfrak{n}/\mathfrak{m}^n,\bbmu;a,\nu)\overline{a_{\pi'}(\mathfrak{n}/\mathfrak{m}^n,\bbmu;a,\nu)},
\end{equation}
extending \eqref{n-prime2S} to all ideal $\mathfrak{n}\in I$.

\subsubsection{Global Rankin-Selberg estimates}\label{sec:RS}

We now recall a few basic analytic properties of the Rankin-Selberg $L$-function $L(s,\pi\times\tilde\pi')$ associated with a pair $(\pi,\pi')\in \mathscr{F}_{\kq,\chi}(Q)\times \mathscr{F}_{\kq,\chi}(Q)$. 

The convexity bound of Li \cite{Li} (see also \cite{Brumley_2} for the cases $n=3,4$) states that
\begin{equation}\label{RS-convex}
\left(\frac{s-1}{s-2}\right)L(s,\pi\times\tilde\pi')\ll C(\pi\times\tilde\pi',s)^{(1-\sigma)/2} \qquad ({\rm Re} (s)\leq 1).
\end{equation}
We have the factorization $C(\pi\times\tilde\pi',s)=D^{n^2}N_{\pi\times\tilde\pi'}K_{\pi\times\tilde\pi'}(s)$. For $\pi_f$ and $\pi_f'$ of conductor $\kq$, whose central characters are equal up to an unramified twist, Theorem \ref{thm:improvedBH} of Appendix \ref{AppendixB} implies that
\begin{equation}\label{conductor-pairs}
N_{\pi\times\tilde\pi'}\leq {\rm Norm}(\kq)^{2n-2}.
\end{equation}
Moreover, the bounds \cite[Lemma A.2]{B2} imply 
\begin{equation}\label{conductor-pairs2}
K_{\pi\times\tilde\pi'}(s)\ll (1+|s|)^{dn^2}(K_{\pi}K_{\tilde\pi'})^n.
\end{equation}
We deduce that, for $\pi,\pi'\in\mathscr{F}_{\kq,\chi}(Q)$, we have
\begin{equation}\label{eq:conv}
\left(\frac{s-1}{s-2}\right)L(s,\pi\times\tilde\pi')\ll (D^{-n^2}(1+|s|)^{dn^2}{\rm Norm}(\kq)^{-2}Q^{2n})^{(1-\sigma)/2} \qquad ({\rm Re} (s)\leq 1).
\end{equation}

The function $L(s,\pi\times\tilde\pi')$ is regular at $s=1$ if and only if $\pi'\neq\pi$. In the case where $\pi'=\pi$, we have a lower bound of polynomial type on the residue at $s=1$. Indeed, \cite[Theorem 3]{Brumley} establishes the existence of an $A>0$ such that
\begin{equation}\label{eq:res}
\Run L(s,\pi\times\tilde{\pi})\gg Q^{-A}.
\end{equation}

\begin{remark}\label{rem:residue}
In \cite{Brumley_2} it is shown that $\Run L(s,\pi\times\tilde\pi)\gg C(\pi\times\tilde\pi)^{-\frac78+\frac{5}{8n}-\varepsilon}$. From the upper bound \eqref{conductor-pairs}, one obtains $7n/4-5/4+\epsilon$ as an admissible value of $A$. This exponent will not play a role in Theorem \ref{thm:poly-growth}.
\end{remark}

\subsection{Refining the set-up in \S\ref{sec:outline}}\label{sec:comb-decomp}
We put $\mathscr{F}_{\kq,\chi,\mathscr{C}}(Q)=\mathscr{F}(Q)\cap \mathscr{A}_{\rm cusp}(\kq,\chi,\mathscr{C})$. Then 
\begin{equation}\label{combinatorial}
|\mathscr{F}_{\kq,\chi}(Q)|=\sum_{\mathscr{C}\in\mathrm{Comb}_S}|\mathscr{F}_{\kq,\chi,\mathscr{C}}(Q)|.
\end{equation}
We shall prove that
\begin{equation}\label{FqQ}
|\mathscr{F}_{\kq,\chi,\mathscr{C}}(Q)|\ll_\epsilon (D^{-n^2}{\rm Norm}(\kq)^{-2}Q^{2n})^{1+\epsilon},
\end{equation}
uniformly in $\kq$. Note that for every $v$ we have $|\mathscr{C}_v|=O_n(1)$. Thus the number of terms in \eqref{combinatorial} is $|\mathscr{C}|=O(|S|^{O_n(1)})=O(\log^{O_n(1)} Q)$. Inserting this into \eqref{combinatorial} and \eqref{LevelSum} will then prove Theorem \ref{thm:poly-growth}.

Recall the set $\mathscr{I}^S$ of $\mathcal{P}_{n-1}$-decorated prime-to-$S$ ideals from \S\ref{sec:outline}. We shall now enrich $\mathscr{I}^S$ at the places in $S$ to account for the combinatorial information $\mathscr{C}$. We shall call a $(\mathcal{P}_{n-1},\mathscr{C})$\textit{-decorated ideal} a triple $(\mathfrak{n},\bbmu, (a,\nu))$, where $\mathfrak{n}\in I$ is an integral ideal, $\bbmu\in\mathscr{P}_{n-1}(\mathfrak{n})$, and $(a,\nu)\in {\rm Index}(\mathscr{C})$. We shall generally write this as $(\mathfrak{n},\bbmu; a,\nu)$. The set of such triples shall be denoted $\mathscr{I}_S$. We have a map $\mathscr{I}_S\rightarrow I$, $(\mathfrak{n},\bbmu;a,\nu)\mapsto\mathfrak{n}$, where we forget the decorations and take the underlying ideal $\mathfrak{n}$. We sometimes write $\tilde{\mathfrak{n}}$ for an element in $\mathscr{I}_S$ with underlying ideal $\mathfrak{n}$. Let $\mathscr{I}_S(X)$ denote the set of $\tilde{\mathfrak{n}}\in\mathscr{I}_S$ with ${\rm Norm}(\mathfrak{n})\leq X$.

Let $V_S(X)$ be the vector space of complex valued functions on $\mathscr{I}_S$. Endow $V_S(X)$ with the standard scalar product
\[
\langle a,b\rangle=\sum_{\tilde{\mathfrak{n}}\in\mathscr{I}_S(X)} f(\tilde{\mathfrak{n}})\overline{g(\tilde{\mathfrak{n}})}.
\]
We shall map $\mathscr{A}_{\rm cusp}(\kq,\chi,\mathscr{C})$ to $V_S(X)$ by sending $\pi\in\mathscr{A}_{\rm cusp}(\kq,\chi,\mathscr{C})$ to the vector $\mathbf{v}_\pi\in V_S(X)$ given by the formula
\[
\mathbf{v}_\pi(\tilde{\mathfrak{n}})=\sqrt{F_X(\mathfrak{n})}a_\pi(\tilde{\mathfrak{n}}),
\]
where $f:\R\rightarrow\R$ is as in \S\ref{sec:outline} and 
\[
F_X(\mathfrak{n})=\sum_\mathfrak{m}f({\rm Norm}(\mathfrak{nm}^n)/X).
\]

The above enrichment allows us to identify the inner product $\langle \mathbf{v}_\pi, \mathbf{v}_{\pi'}\rangle$ in terms of the \textit{full} finite part Rankin-Selberg $L$-function. Indeed, by \eqref{full-coeff} and Mellin inversion we have
\begin{equation}\label{beef}
\langle \mathbf{v}_\pi,\mathbf{v}_{\pi'}\rangle=\sum_{\tilde{\mathfrak{n}}\in\mathscr{I}_S}F_X(\mathfrak{n})a_\pi(\tilde{\mathfrak{n}})\overline{a_{\pi'}(\tilde{\mathfrak{n}})}=\frac{1}{2\pi i}\int_{(2)}L(s,\pi\times\tilde{\pi}')\hat{f}(s)X^sds.
\end{equation}
The above formula is the culmination of the combinatorial explication of the Rankin-Selberg $L$-functions in \S\ref{sec:outline}-\ref{sec:comb-decomp}. It is the basis of the following section.

\subsection{Executing steps \eqref{step1} and \eqref{step2}}

We now execute the first two steps of the proof outline in \S\ref{sec:outline}, using the facts we collected from Rankin-Selberg theory in \S\ref{sec:RS}.

\subsubsection{First step} We begin by establishing the following result.

\begin{prop}\label{prop:mult-1} Let $\epsilon>0$ and $X\gg (D^{-n^2}{\rm Norm}(\kq)^{-2}Q^{2n})^{1+\epsilon}$. Then the map $\mathscr{F}_{\kq,\chi,\mathscr{C}}(Q)\rightarrow V_S(X)$ given by $\pi\mapsto \mathbf{v}_\pi$ is injective.
\end{prop}

\begin{proof} Indeed, \cite[Theorem 7]{Brumley} shows the existence of a $B>0$ such that when $X\gg Q^B$ any pair $(\pi,\pi')\in\mathscr{F}_\kq(Q)\times \mathscr{F}_\kq(Q)$ satisfying $a_\pi(\tilde{\mathfrak{n}})=a_{\pi'}(\tilde{\mathfrak{n}})$ for $\tilde{\mathfrak{n}}\in\mathscr{I}_S$ lies along the diagonal $\pi=\pi'$. It is shown in \cite{LW} that an admissible value for the exponent $B$ is $2n+\varepsilon$, for any $\varepsilon>0$. In fact, their result can be refined, under the assumption that $\pi_f$ and $\pi_f'$ have the same (finite) conductor $\kq$ and central character $\chi$. Indeed, in this case, the bounds of Theorem \ref{thm:improvedBH} of Appendix \ref{AppendixB} save ${\rm Norm}(\kq)^2$ off of this.\end{proof}

\subsubsection{Second step} As in \S\ref{sec:outline}, we let
\[
\mathbf{u}_\pi=\frac{\mathbf{v}_\pi}{\langle \mathbf{v}_\pi,\mathbf{v}_\pi\rangle^{1/2}}
\]
be the projection of the vector $\mathbf{v}_\pi$ to the unit sphere in $V$. We now proceed to show that the vectors $\mathbf{u}_\pi$ and $\mathbf{u}_{\pi'}$ (for $\pi\neq\pi'$) are quasi-orthogonal, in a quantifiable sense.

\begin{prop}\label{prop:IP}
Let $(\pi,\pi')\in\mathscr{F}_{\kq,\chi,\mathscr{C}}(Q)\times\mathscr{F}_{\kq,\chi,\mathscr{C}}(Q)$. For $\epsilon>0$ let
\[
X\gg (D^{-n^2}{\rm Norm}(\kq)^{-2}Q^{2n})^{1/2+\epsilon}.
\] 
If $\pi'\neq\pi$ then $\langle \mathbf{u}_\pi,\mathbf{u}_{\pi'}\rangle\ll_{\epsilon,r} Q^{-r}$ for all $r>0$.
\end{prop}

\begin{proof}
We shall show that there is $C>0$ such that for $X \gg (D^{-n^2}{\rm Norm}(\kq)^{-2}Q^{2n})^{1/2+\epsilon}$ any pair $(\pi,\pi')\in\mathscr{F}_{\kq,\chi,\mathscr{C}}(Q)\times\mathscr{F}_{\kq,\chi,\mathscr{C}}(Q)$ satisfies
\begin{equation}\label{eq:IP-est}
\begin{cases}
\langle \mathbf{v}_\pi,\mathbf{v}_\pi\rangle\gg Q^{-C},&\text{if $\pi=\pi'$;}\\ 
\langle \mathbf{v}_\pi,\mathbf{v}_{\pi'}\rangle\ll_{\epsilon,r} Q^{-r}, & \text{if $\pi\neq\pi'$.}
\end{cases}
\end{equation}
(We can take $C=3n/4-9/4-\epsilon$ for any $\epsilon>0$, but this value is irrelevant for the proof of this proposition; see Remark \ref{rem:residue}.) These two estimates imply the result.

Recall the identity \eqref{beef}. By hypothesis $\hat{f}(1)=1$, and since $f$ is of compact support, $\hat{f}(s)$ is entire. Using \eqref{eq:conv}, we shift the contour to $(-r)$ for $r>0$ to obtain
\[
\langle \mathbf{v}_\pi,\mathbf{v}_{\pi'}\rangle=\Run L(s,\pi\times\tilde\pi')X+O_r(( D^{-n^2}{\rm Norm}(\kq)^{-2}Q^{2n})^{(1+r)/2}X^{-r}).
\]
If $\pi\neq\pi'$ then the residual term vanishes, and hence
\[
\langle \mathbf{v}_\pi,\mathbf{v}_{\pi'}\rangle\ll_r ( D^{-n^2}{\rm Norm}(\kq)^{-2}Q^{2n})^{(1+r)/2}X^{-r}.
\]
If $\pi=\pi'$ we recall the lower bound \eqref{eq:res}. This produces
\[
\langle \mathbf{v}_\pi,\mathbf{v}_\pi\rangle\gg Q^{-A}X+O_r(( D^{-n^2}{\rm Norm}(\kq)^{-2}Q^{2n})^{(1+r)/2}X^{-r}).
\] 
Letting $X\gg (D^{-n^2}{\rm Norm}(\kq)^{-2}Q^{2n})^{1/2+\epsilon}$, we take $r$ sufficiently large (relative to $n$ and $\epsilon$) to arrive at the two estimates in \eqref{eq:IP-est}.
\end{proof}

\begin{remark}
We note that we could avoid quoting the convexity bound \eqref{RS-convex} of \cite{Li} by dualizing the $L$-function, as was done, for example, in \cite{LW}. This does not, however, lead to an improvement in the resulting bounds.
\end{remark}
\begin{remark}\label{rem:n-square}
The analog of Proposition \ref{prop:IP}, when stated with $\langle\mathbf{v}_\pi^S,\mathbf{v}_{\pi'}^S\rangle$, would incur a loss of $n^2$ in the power of $Q$. Indeed, with the set-up of \S\ref{sec:outline} one needs to bound 
\[
L^S(s,\pi\times\tilde\pi')=L(s,\pi\times\tilde\pi')L_S(s,\pi\times\tilde\pi')^{-1}
\]
for ${\rm Re}(s)=\sigma\rightarrow-\infty$. From \eqref{RS-L-factor}, each local correction factor $L_v(s, \pi\times\tilde\pi')^{-1}$ is the product of at most $n^2$ local factors of the form $1-\alpha_{\pi\times\pi'}(\kp;\nu,j,k)^eq_v^{-es}$. By Remark \ref{rem:local-roots}, and in particular the bound \eqref{JS-RS} on the Rankin-Selberg local roots, we deduce that $L_v(s, \pi\times\tilde\pi')^{-1}\ll (1+q_v^{1-\sigma})^{n^2}$. Thus, for ${\rm Re}(s)=\sigma<0$, we have $L_S(s,\pi\times\tilde\pi')^{-1}\ll (\prod_{v\mid\kq}q_v)^{n^2(1-\sigma)}$, which accounts for the weakened exponent.  Moreover, the same loss by $Q^{n^2}$ would arise in the proof of Proposition \ref{prop:mult-1}, were we only to assume that $\pi_{\kp}\simeq\pi_{\kp}'$ for $\kp\nmid\kq$.

Note that in the critical strip the correction factor $L_S(s,\pi\times\tilde\pi')^{-1}$ is uniformly bounded under the Ramanujan conjecture (see \eqref{RS-tempered}). Nevertheless, a contour shift to anywhere within the critical strip leads to insufficient correlation bounds relative to known sphere packing bounds. 
\end{remark}

\subsection{Executing step \eqref{step3}}
We finally come to the fact that a large dimensional sphere can only contain so many quasi-orthogonal vectors.

Let $N$ denote the cardinality of the set $\mathscr{I}$; this is the same as the dimension of $V_S(X)$. Denote by $K$ the cardinality of $\mathscr{F}_{\kq,\chi,\mathscr{C}}(Q)$. We shall show that if $X=(D^{-n^2}{\rm Norm}(\kq)^{-2}Q^{2n})^{1+\epsilon}$, then $K\leq N$. Since $N\asymp X=(D^{-n^2}{\rm Norm}(\kq)^{-2}Q^{2n})^{1+\epsilon}$, this will complete the proof of \eqref{FqQ}, and hence of Theorem \ref{thm:poly-growth}.

By our choice of $X$, we may apply both Propositions \ref{prop:mult-1} and \ref{prop:IP}, so that $\mathscr{F}_{\kq,\chi,\mathscr{C}}(Q)$ can be viewed as a finite system of unitary quasi-orthogonal vectors in $V_S(X)$. The following abstract result establishes the desired bound. To apply it to our situation, we identify $V_S(X)=\mathbb{C}^N=\R^M$, where $M=2N$.

\begin{prop}\label{lem:sphere}
Let $M\geq 2$ and put $V=\R^M$. Let $\mathbf{u}_1,\ldots ,\mathbf{u}_K\in V$ be unitary vectors such that $|\langle \mathbf{u}_i,\mathbf{u}_j\rangle| < M^{-1}$ for $i\neq j$. Then $K\leq M$.
\end{prop}

Before passing to the proof of Proposition \ref{lem:sphere}, we make several remarks.
\begin{remark}
The conclusion of the proposition is sharp, since one can certainly put $M$ orthonormal vectors (and no more) on the unit sphere in $\R^M$. The idea of the proof of Proposition \ref{lem:sphere} is that, in high dimensions, a $1/M$ error off of strict orthogonality is imperceptible. (In fact a $\frac12M^{-1/2}$ error is provably imperceptible: see Remark \ref{rem:codes}.)

Note that the quasi-orthogonality relations established in Proposition \ref{prop:IP} for the family $\{\mathbf{u}_\pi: \pi\in\mathscr{F}(Q)\}$ are much stronger (rapid decay) that the required bounds for Proposition \ref{lem:sphere}. However, it is of no advantage to have $O_r(M^{-r})$ correlation decay, instead of the required rate of $O(1/M)$, since in any case, strictly vanishing off-diagonal correlations (an orthonormal basis) still produces $K=M$.
\end{remark}

\begin{remark}\label{rem:codes}
Let $M\geq 2$ and $\theta\in[0,\pi)$. Denote by $A(M,\theta)$ the maximum cardinality of a subset $\{\mathbf{u}_1,\ldots ,\mathbf{u}_K\}$ of $S^{M-1}$ with $\max_{i\neq j}\langle \mathbf{u}_i,\mathbf{u}_j\rangle\leq \cos\theta$. Such a subset is called a \textit{spherical code}.

If $\theta>\pi/2$, then an elementary argument shows that $A(M,\theta)$ is bounded by an expression depending only on $\theta$. Indeed, as remarked in \cite[\S 3.2]{HV}, we have
\begin{equation}\label{elementary}
0\leq\langle \mathbf{u}_1+\cdots +\mathbf{u}_K,\mathbf{u}_1+\cdots +\mathbf{u}_K\rangle\leq K+K(K-1)\cos\theta.
\end{equation}
Thus, if $\cos\theta$ is strictly negative, this provides a bound for $A(M,\theta)$ which depends only on $\theta$. In particular, if $\theta>\pi/2$ is \textit{fixed}, then $A(M,\theta)$ is bounded uniformly in $M$.

If $\theta>\pi/2$ is now allowed to depend on $M$, the inequality \eqref{elementary} still yields an upper bound on $A(M,\theta)$. For example, if $\cos\theta=-M^{-\alpha}$, for $\alpha\geq 0$, we obtain $A(M,\theta)\leq M^\alpha$. As $\alpha$ varies through the interval $[0,1]$, this provides an interpolation of the uniformly bounded range (where $\theta>1/2$ is fixed) and the range treated by Proposition \ref{lem:sphere}.

On the other hand, when $\theta<\pi/2$ is fixed, then $A(M,\theta)$ grows exponentially in $M$. The work of Kabatjanskii-Levenshtein \cite{KL} provides upper bounds in this regime. It is known, however, that for $\theta=\pi/2-c/\sqrt{M}$, one still retains a polynomial upper bound. See \cite[Theorem 2.1]{FKM} and \cite{Tao}. Indeed, Lemma 2 of \textit{loc. cit} shows that one retains a \textit{linear bound} as long as $\cos\theta\leq \frac12M^{-1/2}$. This latter result would in fact be sufficient for our purposes.
\end{remark}

\begin{proof}
An elementary exercise establishes the result for $M=2$. Suppose the result is true in dimension $M-1$. We claim this implies the result in dimension $M$.

Let $W$ be the orthogonal complement to $\mathbf{u}_K$ in $V=\R^M$. For every $1\leq i\leq K-1$ let $\mathbf{w}_i$ be the projection of the vector $\mathbf{u}_i$ to $W$. We define $\lambda_i\in\R$ by the equality $\mathbf{w}_i=\mathbf{u}_i-\lambda_i \mathbf{u}_K$; then $\lambda_i= \langle\mathbf{u}_i,\mathbf{u}_K\rangle$. For $1\leq i, j\leq K-1$ we have $\langle \mathbf{w}_i,\mathbf{w}_j\rangle=\langle \mathbf{u}_i,\mathbf{u}_j\rangle+\lambda_i\lambda_j$.  If $i=j$ we obtain $\|\mathbf{w}_i\|^2=1+\lambda_i^2$. By hypothesis, $|\lambda_i|<1/M$ which implies $\|\mathbf{w}_i\|^2> 1-M^{-2}=M^{-2}(M^2-1)$.  Moreover, if $i\neq j$ we have $|\langle \mathbf{u}_i,\mathbf{u}_j\rangle|<1/M$; thus
\[
|\langle \mathbf{w}_i,\mathbf{w}_j\rangle|<1/M+1/M^2=M^{-2}(M+1).
\]

Now, consider the $K-1$ unitary vectors $\mathbf{u}'_j=\mathbf{w}_j/\|\mathbf{w}_j\|$ in the $(M-1)$-dimensional subspace $W$. For $1\leq i\neq j\leq K-1$ we have
\[
|\langle \mathbf{u}_i',\mathbf{u}_j'\rangle|=\|\mathbf{w}_i\|^{-1}\|\mathbf{w}_j\|^{-1}|\langle \mathbf{w}_i,\mathbf{w}_j\rangle|< \frac{M^2}{M^2-1}\cdot \frac{M+1}{M^2}=\frac{1}{M-1}.
\]
From our recurrence hypothesis, we deduce that $K-1\leq M-1$, as claimed.
\end{proof}

\begin{remark}
The above induction argument works under the more general hypothesis that $\max_{i\neq j}|\langle \mathbf{u}_i,\mathbf{u}_j\rangle|<f(M)$, for any function $f$ verifying $\frac{f(M)}{1-f(M)}\leq f(M-1)$. But if $f(M)=M^{-\alpha}$, this inequality reads $(1-1/M)^\alpha\leq 1-M^{-\alpha}$. The left-hand side is approximated by $1-\alpha/M$, and one sees that one can do no better than $\alpha=1$.
\end{remark}

\subsection{Proof of Theorem \ref{under-RR}}\label{sec:Riemann}

We now address the question of improving the upper bound on $|\mathscr{F}(Q)|$ in Theorem \ref{thm:poly-growth}, under the Riemann hypothesis for Rankin-Selberg $L$-functions as well as the Ramanujan conjecture at finite places for members of $\mathscr{A}_{\rm cusp}$.

It is easy to see that the exponent of $2n$ in Theorem \ref{thm:poly-growth} can be improved to $n+1$ under these assumptions, and that the discriminant dependence is as described there. This is due to the fact that, under Riemann and Ramanujan, the map $\pi\mapsto\mathbf{v}_\pi$ is injective as soon as $X\gg\log^2Q$ (see, for example, \cite[Proposition 5.22]{IK}). This replaces step \eqref{step1} in the proof of Theorem \ref{thm:poly-growth}. On the other hand, the proof of Proposition \ref{prop:IP} is insensitive to the Riemann hypothesis and the Ramanujan conjecture, despite the fact that the residue of the $L(s,\pi\times\tilde\pi)$ is bounded below by $1/\log Q$ under these assumptions (see \cite[Theorem 5.19]{IK}). In any case, with Theorem \ref{prop:mult-1} improved, we may take $X=(D^{-n^2}{\rm Norm}(\kq)^{-2}Q^{2n})^{1/2+\epsilon}$ in executing step \eqref{step3}. Indeed, the exponent of $Q$ required for the value of $X$ in step \eqref{step3} is the maximum of the exponents coming from Propositions \ref{prop:mult-1} and \ref{prop:IP}. Inserting this into \eqref{combinatorial} and \eqref{LevelSum} will then prove Theorem \ref{under-RR}.

\section{A bound for the Artin exponent of a pair,\\ by Colin J. Bushnell and Guy Henniart}\label{AppendixB}

Let $F$ be a locally compact non-Archimedean field, and $n,m$ two positive integers. Let $\pi$ be a smooth irreducible representation of $\GL_n(F)$, with central character $\omega_\pi$ and Artin conductor $\Ar(\pi)=a$, and let $\rho$ be a smooth irreducible representation of $\GL_m(F)$, with central character $\omega_\rho$ and Artin conductor $\Ar(\rho)=b$.

In \cite{BH} and \cite[Theorem C]{BH2017}, we proved that the pair $(\pi,\rho)$ satisfies
\begin{equation}\label{originalBH}
\Ar(\pi\times\rho)\leq ma +nb-\min (a,b).
\end{equation}
That bound cannot be improved in general but here, prompted by a query of F. Brumley, we improve \eqref{originalBH} under an additional hypothesis.

\begin{thm}\label{thm:improvedBH}
Assume that $\omega_\pi\omega_\rho$ is unramified. Then
\begin{equation}\label{improvedBH}
\Ar(\pi\times\rho)\leq ma+nb-2\min (a,b).
\end{equation}
\end{thm}
When $n=m$ and $a=b$, this gives $\Ar(\pi\times\rho)\leq (2n-2)a$, as used in the main text. Note also that when $n=m=1$ the hypothesis implies $a=b$ and $\Ar(\pi\times\rho)=0$, which is fortunate since the right hand side of \eqref{improvedBH} is also $0$!

Thanks to the Langlands correspondence, we may express the theorem in terms of Weil-Deligne representations, and we indeed use that language in the proofs. We fix a separable algebraic closure $F^{\rm sep}$ of $F$ and let $W_F$ be the Weil group of $F^{\rm sep}$ over $F$. We write $\sigma, \tau$ for the Weil-Deligne representations corresponding to $\pi, \rho$: they are directs sums of indecomposable Weil-Deligne representations. The theorem above is then equivalent to
\begin{thm}
Assume that $\det\sigma\det\tau$ is unramified. Then
\begin{equation}\label{improvedBH'}
\Ar(\sigma\otimes\tau)\leq ma+nb-2\min(a,b).
\end{equation}
\end{thm}

\begin{remark}
Assume that $\sigma$ is the direct sum of characters of $W_F$, all trivial but one, which then has to be $\det\sigma$. Take for $\tau$ the contragredient $\tilde\sigma$ of $\sigma$. Then $\Ar(\sigma\otimes\tilde\sigma)=(2n-2)a$, so one cannot improve \eqref{improvedBH'} or \eqref{improvedBH} in general, even assuming that $\tau=\tilde\sigma$.
\end{remark}

We now proceed to the proof, relying on the results and techniques of \cite{BH2017}.

\subsection{}\label{sec:B1}\hspace{-0,4cm} A basic point is a stronger inequality than \eqref{originalBH}, when $\sigma$ and $\tau$ are indecomposable.

\begin{lem}\label{indecomp1}
Assume $\sigma, \tau$ indecomposable. Then
\[
\Ar(\sigma\otimes\tau)/nm\leq \max(a/n, b/m),\quad\text{with equality if }\, a/n\neq b/m.
\]
\end{lem}

\begin{proof}
The case of inequality is \cite[Proposition 6.3]{BH2017}. The case of equality can be deduced from \cite[Proposition 5.5]{BH2017}.
\end{proof}
\begin{lem}\label{indecomp2}
Assume $\sigma$ indecomposable. Then $\Ar(\det\sigma)\leq a/n$.
\end{lem}

\begin{proof}
By \cite[Fact 2.1]{BH2017} and the notation there, we have $\sigma={\rm St}_r(\sigma')$, for some positive integer $r$ and some irreducible representation $\sigma'$ of $W_F$.

If $\sigma'$ is an unramified character of $W_F$ then $r=n$ and $a=n-1$, whereas $\det\sigma$ is unramified, so $\Ar(\det\sigma)=0\leq 	\frac{a}{n}$.

If $\sigma'$ is not an unramified character, then $a=r\Ar(\sigma')$ and $\det\sigma=(\det\sigma)^r$, so it is enough to treat the case where $\sigma=\sigma'$ is irreducible (and not unramified). But then $a-n$ is the Swan exponent of $\sigma$, so, using \cite[Fact 2.3]{BH2017},
\[
\frac{a}{n}-1=\inf\{\epsilon>0: \sigma (W_F^\epsilon)=1\}.
\]
Since $\det\sigma$ is certainly trivial on the ramification subgroup $W_F^\epsilon$ if $\sigma$ is, we see that the Swan exponent of $\det\sigma$ is at most $\frac{a}{n}-1$, so $\Ar(\det\sigma)\leq a/n$.
\end{proof}

Let us define the list of slopes $\sigma$. When indecomposable, $\sigma$ has a list of $n$ slopes, all equal to $a/n$. In general the list of slopes of $\sigma$ is obtained by gathering the lists of slopes of its indecomposable summands, in increasing order. We write $(a_1,\ldots , a_n)$ for the list of slopes of $\sigma$, and $(b_1,\ldots ,b_m)$ for the list of slopes of $\tau$; in particular, $a=a_1+\cdots + a_n$ and $b=b_1+\cdots +b_m$.

Applying Lemmas \ref{indecomp1} and \ref{indecomp2} to the indecomposable summands of $\sigma$ and $\tau$ we get
\begin{cor}\label{corollary} The following holds:
\begin{enumerate}
\item[(i)] $\Ar(\det\sigma)\leq a_n$
\item[(ii)] $\Ar(\sigma\otimes\tau)\leq n\Ar(\tau)$, if $a_n\leq b_1$, with equality if $a_n<b_1$.
\end{enumerate}
\end{cor}

\subsection{}\label{sec:B2}\hspace{-0,4cm} In this n$^{\rm o}$, we assume $n=1$. As the case $n=m=1$ is done, we also assume $m>1$.

We first deal with the case $b_{m-1}<b_m$. Then we can write $\tau=\tau'\oplus\eta$ for a character $\eta$ of $W_F$ with $\Ar(\eta)=b_m$. By the Corollary \ref{corollary} (i), $\Ar(\det\tau')\leq b_{m-1}$ and since $\det\tau=(\det\tau')\eta$ we get $\Ar(\det\tau)=b_m$. But $\sigma=\det\sigma$ and $\det\sigma\det\tau$ is unramified, so we have $a=b_m$.

By the Corollary \ref{corollary} (ii), we have $\Ar(\sigma\otimes\tau')=(m-1)a$ since $a=b_m>b_{m-1}$. We also have $\Ar(\sigma\otimes\eta)=\Ar((\det\sigma)\eta)=\Ar((\det\tau)^{-1}\eta)$, so $\Ar(\sigma\otimes\eta)=\Ar(\det\tau')\leq b_{m-1}$. Adding, we get $\Ar(\sigma\otimes\tau)\leq (m-1)a+b_{m-1}$.

On the other hand, $b\geq b_m=a$ hence $\min(a,b)=a$ and
\[
ma+b-2\min(a,b)=(m-2)a+b\geq (m-1)a+b_{m-1}
\]
because $b\geq b_{m-1}+b_m=a+b_{m-1}$. We have proved \eqref{improvedBH'} when $b_{m-1}< b_m$.

We now assume that $b_{m-1}=b_m$. By Corollary \ref{corollary} (i), $\Ar(\det\tau)\leq b_m$ and, reasoning as above, we now get $a\leq b_m$ from Corollary \ref{corollary} (ii). Write $\tau=\tau'\oplus\eta$, where $\eta$ is a Weil-Deligne representation with dimension $d\geq 2$ and all slopes equal to $b_m$. Let $b'=\Ar(\tau')$, so $b=b'+db_m$. We have $\Ar(\sigma\otimes\tau')\leq (m-d)a+b'-\min(a,b')$: this follows from \eqref{originalBH} if $\tau'\neq 0$ and $m=d$, $b'=0$, if $\tau'=0$.

On the other hand, $\Ar(\sigma\otimes\eta)\leq db_m$ by Corollary \ref{corollary} (ii), since $a\leq b_m$. Adding, we obtain
\[
\Ar(\sigma\otimes\tau)\leq (m-d)a+b'+db_m-\min (a,b'),
\]
so the result follows, provided $da+\min(a,b')\geq\min (a,b'+db_m)$, which is clear since $d\geq 2$. This again proves \eqref{improvedBH'}.

\subsection{}\label{sec:B3}\hspace{-0,4cm} From now on we assume $n>1$ and $m>1$.

We first deal with the situation where $a_{n-1}<a_n$ and $b_{m-1}<b_m$. Accordingly, we write $\sigma=\sigma'\oplus\chi$ for a character $\chi$ of $W_F$ with $\Ar(\chi)=a_n$, and $\tau=\tau'\oplus\eta$, for a character $\eta$ of $W_F$ with $\Ar(\eta)=b_m$. We put $a'=\Ar(\sigma')$, $b'=\Ar(\tau')$, so $a=a'+a_n$ and $b=b'+b_m$. Reasoning as above, we get $\Ar(\det\sigma')\leq a_{n-1}$, $\Ar(\det\tau')\leq b_{m-1}$, $a_n=b_m$, and $\Ar(\chi\eta)\leq \max(a_{n-1},b_{m-1})$.

On the other hand, we have by \eqref{originalBH}
\[
\Ar(\sigma'\otimes\tau')\leq (m-1)a'+(n-1)b'-\min(a',b'),
\]
and, by Corollary \ref{corollary} (ii) again, $\Ar(\sigma'\otimes\eta)=(n-1)b_m=(n-1)a_n$ and $\Ar(\chi\otimes\tau')=(n-1)a_n$. Adding, we get
\[
\Ar(\sigma\otimes\tau)\leq (m-1)a+(n-1)b-\min (a',b')+\max(a_{n-1},b_{m-1}).
\]
The result then follows provided that $a+b+\min(a',b')\geq 2\min(a,b)+\max(a_{n-1},b_{m-1})$, or, equivalently, 
\begin{equation}\label{suffices}
a'+b'+\min(a',b')\geq 2\min (a',b')+\max(a_{n-1},b_{m-1}).
\end{equation}
But $a'+b'=\min(a',b')+\max(a',b')$ and $\max(a',b')\geq (a_{n-1},b_{m-1})$ because $a'\geq a_{n-1}$ and $b'\geq b_{m-1}$, establishing \eqref{suffices}.

\subsection{}\label{sec:B4}\hspace{-0,4cm} We turn to the case where $a_{n-1}<a_n$ but $b_{m-1}=b_m$. Write $\sigma=\sigma'\oplus\chi$ for a character $\chi$ of $W_F$ with $\Ar(\chi)=a_n$, and $\tau=\tau'\oplus\eta$ for a Weil-Deligne representation $\eta$ with dimension $d\geq 2$ and all slopes equal to $b_m$. Put $a'=\Ar(\sigma')$, $b'=\Ar(\tau')$, so $a=a'+a_n$, $b=b'+d b_m$.

As in the second case of \S\ref{sec:B2}, we get $a_n\leq b_m$ and $\Ar(\chi\otimes\eta)\leq d b_m$ by Corollary \ref{corollary} (ii).

By \eqref{originalBH} (or because $\tau'=0$) we have
\[
\Ar(\sigma'\otimes\tau')\leq (m-d)a'+(n-1)b'-\min (a',b').
\]
Because $a_n\leq b_m$ we have $a_{n-1}<b_m$ so, by Corollary \ref{corollary} (ii), $\Ar(\sigma'\otimes\eta)=(n-1)d b_m$. Applying Lemma \ref{indecomp1} to $\chi\otimes\tau_j$ where $\tau_j$ is an indecomposable summand of $\tau'$, we obtain $\Ar(\chi\otimes\tau')\leq\sum_{i=1}^{m-d}\max(a_n,b_i)$. Adding gives
\begin{equation}\label{Artin-sigma-tau}
\Ar(\sigma\otimes\tau)\leq (m-d)a'+ndb_m+(n-1)b'-\min (a',b')+\sum_{i=1}^{m-d}\max(a_n,b_i).
\end{equation}
We claim that the right-hand side of \eqref{Artin-sigma-tau} is at most $m(a'+a_n)+n(b'+db_m)-2\min(a'+a_n,b'+db_m)$, or, equivalently, that
\begin{equation}\label{claim}
2\min (a'+a_n,b'+db_m)+\sum_{i=1}^{m-d}\max(a_n,b_i)\leq ma_n+da'+b'+\min(a',b').
\end{equation}
Indeed, since $\sum_{i=1}^{m-d}\max(a_n,b_i)\leq (m-d)a_n+b'$ and $d\geq 2$, we have
\[
2\min(a'+a_n,b'+db_m)+\sum_{i=1}^{m-d}\max(a_n,b_i)\leq 2(a'+a_n)+(m-d)a_n+b'\leq ma_n+da'+b',
\]
establishing \eqref{claim}. By symmetry, the case where $a_{n-1}=a_n$ but $b_{m-1}<b_m$ also holds.

\subsection{}\label{sec:B5}\hspace{-0,4cm} The final case is when $a_{n-1}=a_n$ and $b_{m-1}=b_m$. Here the hypothesis that $\det\sigma\det\tau$ is unramified plays no role. By symmetry we may and do assume $a_n\leq b_m$.

We write $\sigma=\sigma'\oplus\chi$ for a Weil-Deligne representation $\chi$ with dimension $e\geq 2$ and all slopes equal to $a_n$, and $\tau=\tau'\oplus\eta$ as in \S\ref{sec:B4}. We put $a'=\Ar(\sigma')$, $b'=\Ar(\tau')$, so that $a=a'+ea_n$, $b=b'+db_m$.

By \eqref{originalBH} (or because $\sigma'$ or $\tau'$ is 0) we have $\Ar(\sigma'\otimes\tau')\leq (m-d)a'+(n-e)b'-\min(a',b')$.  Since $a_n\leq b_p$ by hypothesis, we get from Corollary \ref{corollary} $\Ar(\sigma'\otimes\eta)\leq d(n-e)b_m$ and $\Ar(\chi\otimes\eta)\leq db_m$. As in \S\ref{sec:B4} we get, from Lemma \ref{indecomp1}, $\Ar(\chi\otimes\tau')\leq \sum_{i=1}^{m-d}e\max(a_n,b_i)$. Adding, this gives
\begin{equation}\label{adding}
\Ar(\sigma\otimes\tau)\leq (m-d)a'+(n-e)b'+ndb_m-\min(a',b')+e\sum_{i=1}^{m-d}\max(a_n,b_i).
\end{equation}

We claim that the right-hand side of \eqref{adding} is at most $ma'+mea_n+nb'+ndb_m-\min(a'+ea_n,b'+db_m)$. This is equivalent to the inequality
\begin{equation}\label{claim2}
da'+mea_n+eb'+\min(a',b')\geq \min(a'+ea_n,b'+db_m)+e\sum_{i=1}^{m-d}\max(a_n,b_i).
\end{equation}
Indeed using $\sum_{i=1}^{m-d}\max(a_n,b_i)\leq (m-d)a_n+b'$ and $d\geq 2$, as in \S\ref{sec:B4}, we deduce \eqref{claim2}.

\subsection{}\hspace{-0,4cm} With an entirely similar reasoning, but replacing Artin exponents $\Ar$ with Swan exponents $\Sw$, we get the following result, improving \cite[Theorem CS]{BH2017} in a special case.

\begin{thm}
Let $\sigma$, $\tau$ be semisimple representations of $W_F$. Assume that $\Sw(\det\sigma\det\tau)=0$. Then 
\[
\Sw(\sigma\otimes\tau)\leq (\dim\tau)\Sw(\sigma)+(\dim\sigma)\Sw(\tau)-2\min(\Sw(\sigma),\Sw(\tau)).
\]
\end{thm}

\section*{Acknowledgements}

The authors thank Dimitris Koukoulopoulos, James Maynard, Djordje Mili{\'c}evi{\'c}, Paul Nelson, Lillian Pierce, Maksym Radziwi{\l}{\l}, Abhishek Saha, and Kannan Soundararajan for several insightful discussions.  We also thank Colin Bushnell and Guy Henniart for kindly allowing us to include their answers to our questions as an appendix.  Farrell Brumley is partially supported by ANR grant 14-CE25.  Jesse Thorner is partially supported by a NSF Mathematical Sciences Postdoctoral Fellowship, and Asif Zaman is partial supported by a NSERC Postdoctoral Fellowship.  Part of this work was carried out at MSRI, Berkeley during the spring semester of 2017, supported in part by NSF grant DMS 1440140.

\bibliographystyle{abbrv}
\bibliography{GeneralizedLinnik}

\end{document}